\documentclass[preprint,11pt,3p,authoryear]{elsarticle}

\usepackage{lmodern}
\usepackage[english]{babel}
\usepackage[utf8]{inputenc} %caractère accentués
\usepackage{amsmath}
\usepackage{amsthm}
\usepackage[bookmarks=false,colorlinks]{hyperref}
\AtBeginDocument{%
	\hypersetup{
		linkcolor=blue,   
		citecolor=red,
		urlcolor=magenta}}
\usepackage{enumitem}
\usepackage{xparse}%
\usepackage{xkeyval}%
\usepackage{subfloat}
\usepackage{setspace}
\usepackage{tabularx,multirow}
\usepackage[algo2e,noend]{algorithm2e} 
\usepackage{amsfonts}
\usepackage{graphicx}
\usepackage{subfig}
\usepackage{xcolor}
\usepackage{breqn}
\usepackage{float}
\usepackage{bm}
\usepackage{natbib}
\usepackage{amssymb}
\usepackage{xspace}
\usepackage[tableposition=top]{caption}
\usepackage{booktabs}
\usepackage{etoolbox}
\usepackage{textcomp}

\newtheorem{prop}{Proposition}
\newtheorem{definition}{Definition}
\newtheorem{modl}{Model}
\newenvironment{model}{\begin{samepage}\begin{modl}}{\end{modl}\end{samepage}}

\newcommand{\xiz}{\widehat{x}_{i}}
\newcommand{\yiz}{\widehat{y}_{i}}
\newcommand{\xjz}{\widehat{x}_{j}}
\newcommand{\yjz}{\widehat{y}_{j}}
\newcommand{\viz}{\widehat{v}_{i}}
\newcommand{\vjz}{\widehat{v}_{j}}
\newcommand{\xijz}{\widehat{x}_{ij}}
\newcommand{\yijz}{\widehat{y}_{ij}}
\newcommand{\qi}{q_{i}}
\newcommand{\qj}{q_{j}}
\newcommand{\ti}{\theta_{i}}
\newcommand{\tj}{\theta_{j}}
\newcommand{\qmax}{\overline{q}_i}
\newcommand{\qmin}{\underline{q}_i}

\newcommand{\tiz}{\widehat{\theta}_{i}}
\newcommand{\tjz}{\widehat{\theta}_{j}}

\newcommand{\tmin}{\underline{\theta}_i}
\newcommand{\tmax}{\overline{\theta}_i}
\newcommand{\dx}{\delta_{i,x}}
\newcommand{\dy}{\delta_{i,y}}

\newcommand{\z}{z_{ij}}

\newcommand{\tm}{{t}^{\text{min}}_{ij}}

\newcommand{\vijx}{v_{ij,x}}
\newcommand{\vijy}{v_{ij,y}}
\newcommand{\vix}{v_{i,x}}
\newcommand{\viy}{v_{i,y}}
\newcommand{\vjx}{v_{j,x}}
\newcommand{\vjy}{v_{j,y}}

\newcommand{\R}{\mathbb{R}}
\newcommand{\C}{\mathcal{C}_{ij}}
\newcommand{\A}{\mathcal{A}}

\renewcommand{\P}{\mathcal{P}}

\newcommand{\eix}{\epsilon_{i,x}}
\newcommand{\ejx}{\epsilon_{j,x}}
\newcommand{\eiy}{\epsilon_{i,y}}
\newcommand{\ejy}{\epsilon_{j,y}}
\newcommand{\ei}{\epsilon_{i}}
\newcommand{\ej}{\epsilon_{j}}
\newcommand{\heix}{\overline{\epsilon}_{i,x}}

\newcommand{\heiy}{\overline{\epsilon}_{i,y}}

\newcommand{\tvy}{\tilde{v}_{ij,y}}
\newcommand{\tvx}{\tilde{v}_{ij,x}}
\newcommand{\tviy}{\tilde{v}_{i,y}}
\newcommand{\tvix}{\tilde{v}_{i,x}}

\newcommand{\B}{\mathcal{B}_{ij}}
\newcommand{\U}{\mathcal{U}}
\newcommand{\Br}{\tilde{\mathcal{B}}_{ij}}

\newcommand{\Python}{\textsc{Python}\xspace}
\newcommand{\Cplex}{\textsc{Cplex}\xspace}
\newcommand{\RS}{\mathcal{RS}}

\makeatletter
\def\munderbar#1{\underline{\sbox\tw@{$#1$}\dp\tw@\z@\box\tw@}}
\makeatother

\begin{document}
	
\begin{frontmatter}
	
\title{Robust aircraft conflict resolution under trajectory prediction uncertainty} 

\author[UNSW]{Fernando H. C. Dias}
\author[UNSW]{David Rey\corref{mycorrespondingauthor}}
\ead{drey@unsw.edu.au}
\cortext[mycorrespondingauthor]{Corresponding author}
\address[UNSW]{School of Civil and Environmental Engineering, UNSW Sydney, NSW, 2052, Australia}

\begin{abstract}
We address the aircraft conflict resolution problem under trajectory prediction uncertainty. We consider that aircraft velocity vectors may be perturbed due to weather effects, such as wind, or measurement errors. Such perturbations may affect aircraft trajectory prediction which plays a key role in ensuring flight safety in air traffic control. Our goal is to solve the aircraft conflict resolution problem in the presence of such perturbations and guarantee that aircraft are separated for any realization of the uncertain data. We propose an uncertainty model wherein aircraft velocities are represented as random variables and the uncertainty set is assumed to be polyhedral. We consider a robust optimization approach and embed the proposed uncertainty model within state-of-the-art mathematical programming formulations for aircraft conflict resolution. We then adopt the approach of \cite{bertsimas2004price} to formulate the robust counterpart. We use the complex number reformulation and the constraint generation algorithm proposed by \cite{dias2020disjunctive} to solve the resulting nonconvex optimization problem on benchmarking instances of the literature. Our numerical experiments reveal that perturbations of the order of $\pm 5\%$ on aircraft velocities can be accounted for without significantly impacting the objective function compared to the deterministic case. Our tests also show that for greater levels of uncertainty, several instances fail to admit conflict-free solutions, thus highlighting existing risk factors in aircraft conflict resolution. We attempt to explain this behavior by further analyzing pre- and post-optimization aircraft trajectories. Our findings show that most infeasible instances have both a relatively low total aircraft pairwise minimal distance and a high number of conflicts. 
\end{abstract}

\end{frontmatter}

\section{Introduction}

The usage of airspace has been increasing steadily in the last decades and even with recent challenges, current air traffic management (ATM) systems are working under pressure. Safety stands at the heart of Air traffic control (ATC) which operates within a constrained environment to ensure air space and flight safety. One of the missions of ATC is to guarantee that aircraft respect separation standards at all times to ensure that risks of collision are avoided. This requires the prediction of aircraft trajectories in the near future, typically within the next 30 minutes \citep{rey2016subliminal}. This limitation is largely due to the uncertainty inherent to aircraft trajectory prediction. Weather conditions such as wind, pressure and temperature, and its associated measurement errors may affect the velocity of aircraft thus potentially compromising the accuracy of the trajectory prediction task. The vast majority of exact approaches for aircraft conflict resolution methods do not take into account this uncertainty. This makes their solutions, i.e. conflict-free aircraft trajectories, potentially sensitive to perturbations and an increased risk of violating separation standards.
%In 2005, the International Civil Aviation Organization (ICAO) released an operational map establishing its global program \citep{international2005global}. In this report, the main goals were to balance air traffic operations, to reduce the workload of air traffic controllers and the overall improvement in safety, capacity, efficiency and environmental impact. 

ATC procedures are largely operated by Air Traffic Controllers which face an increasing workload in dense traffic conditions. Decision-support and automation in conflict resolution procedures are emerging alternatives to traditional methods that may help to balance controller workload and increase airspace capacity. The regulations related to civil aviation and the usage of airspace are controlled by the ICAO \citep{icao2010convention}. The ICAO is responsible for defining the standards for commercial and civil aviation and the minimum required separation between aircraft. During cruise stage, which is the focus of this paper, the minimum separation is established as 5 NM horizontally and 1000 ft vertically. Any pair of aircraft within less than those requirements is assumed to be in a conflict. Airspace is typically configured such that aircraft are assigned to flight levels separated by at least 1000 ft, hence vertical separation is guaranteed by design. Congested air traffic networks can lead to loss of separation between aircraft which impairs flight safety and may result in eventual collisions. The aircraft conflict resolution problem (ACRP) is typically formulated as an optimization model in which the objective is to find conflict-free trajectories for a set of aircraft with conflicting trajectories. Different strategies have been used to address this problem based on the type of deconfliction manoeuvres available, i.e. speed control (acceleration or deceleration), heading control, vertical control (flight level reassignment) or a combination of these manoeuvres. 
%For example, in \cite{rey2017complex}, the authors opt for using continuous variables and manoeuvres but with the risk of violation in dense scenario while \cite{omer2015space} proposes a model based smaller sets and focusing on trajectories points to provide violation free solutions.\\

In this paper, we consider the ACRP in the presence of uncertainty onto aircraft velocities. We propose a robust optimization approach which aims to guarantee that aircraft are separated for any realization of the random data and to identify minimum-deviation trajectories. We adopt the robust optimization approach proposed by \cite{bertsimas2004price} to control the level of robustness in the formulation and adapt state-of-the-art solution methods to the deterministic ACRP to solve the resulting robust ACRP. We make the following methodological contributions to the field: i) we propose a model of aircraft trajectory prediction uncertainty based on aircraft velocity components; ii) we show that the proposed uncertainty model can be incorporated in a robust optimization formulation for the ACRP; iii) we adapt the complex number formulation of \cite{rey2017complex} and the proposed exact algorithm of \cite{dias2020two} to solve the robust ACRP; iv) we conduct numerical experiments on benchmarking instances to test the proposed robust ACRP formulation. Our tests reveal that increasing the level of robustness or the size of the uncertainty set rapidly increase the likelihood of infeasibility; and, v) re- and post-optimization analyzes reveal that the number of conflicts and the total minimal pairwise distance between aircraft trajectories can explain the behaviour of the model.\\

The paper is organized as follows. We review the state-of-the-art on aircraft conflict resolution, highlight existing research gaps and position our contributions to the field in Section \ref{lit}. We present a deterministic nonconvex formulation for the ACRP in Section \ref{model}. We then introduce the uncertainty model and present the robust ACRP in Section \ref{robs}. Numerical results are provided in Section \ref{num} and concluding remarks are discussed in Section \ref{con}.

\section{Literature Review}
\label{lit}

\cite{shone2020applications} provides a recent and comprehensive review of stochastic modeling in ATM. We next focus on works which have addressed the problems of aircraft trajectory prediction under uncertainty and approaches for conflict resolution under uncertainty.

\subsection{Aircraft trajectory prediction under uncertainty}

Most of the efforts for trajectory prediction under uncertainty assume that weather is the primary source of randomness, considering that wind, rain and fog may affect aircraft trajectory prediction. More generally, adverse weather conditions can cause delays due to lower visibility, loss of friction in take-off and arrivals, etc, and has been the focus of several works in ATM. In \cite{nilim2001trajectory}, a dynamic routing mechanism was proposed to account for expected delays if the nominal trajectory is inaccessible due to weather conditions. By modelling the uncertainty using a statistical analysis of forecast data, \cite{hentzen2018maximizing} calculated the probability of an aircraft reaching its destination given that some action is taken to avoid the area affected by adverse weather. Similarly, \cite{clarke2009air} used available stochastic weather information into a dynamic model to determine the route capacity for each aircraft. \cite{pepper2003predictability} presented a Bayesian model to incorporate uncertainty from weather into air traffic flow to understand capacity behaviour under weather conditions. A sequential optimization approach was developed by \cite{grabbe2009sequential} to adapt route capacity to account for varying weather conditions. \cite{zheng2011modeling} developed a statistical model of wind uncertainty and trajectory prediction using constant speed. The authors used curated data from previous years, to estimate and calibrate the parameters in their statistical model and used that knowledge to determine optimal manoeuvres. In a more generic uncertainty context, \cite{murcca2017robust} presented a robust approach for optimizing runway usage and taxi-out time and \cite{radmanesh2018grey} solved the problem of path planning for unmanned air vehicles under random circumstances. Different sources of uncertainty were explored by \cite{kim2009air} who discretized flight speed uncertainty using a white Gaussian function and removing the crosswind effect to assess the efficiency of traffic flow. \cite{gonzalez2016wind} analyzed  flight paths determined via pseudo-spectral methods by assessing wind-optimal trajectories. \cite{franco2017optimal} created a structured space and applied the Dijkstra algorithm to obtain optimal paths under wind uncertainty. \cite{rivas2017analysis} and \cite{valenzuela2017sector} analyzed the effects of wind uncertainty in fuel consumption and demand. %\cite{hassan2020mixed} proposed a mixed-integer linear programming formulation for multi-sector planning that combines speed and heading changes to provide conflict-free aircraft trajectories at a more tactical level than traditional conflict resolution approaches.

\subsection{Aircraft conflict resolution under uncertainty}

At an operational level, stochastic approaches to air traffic modeling have mainly focused on calculating the probability of conflict or collision. \cite{paielli1997conflict} proposed an approach to calculate the probability of conflict in cross-track and along-track trajectories. \cite{krozel1997strategic} presented a stochastic model assuming that velocity, speed and heading follow Gaussian distributions. \cite{blin2000stochastic} discussed the role of accounting for uncertainty by evaluating its impact on the error on position, speed and acceleration of aircraft. \cite{prandini2000probabilistic} proposed two models to track perturbations on aircraft trajectories: one for calculating the probability of conflict, and a second which combines a deterministic motion with a Brownian motion perturbation. A similar perturbation model was proposed by \cite{jacquemart2016adaptive} where the authors calculated the probability of conflict and, using rare event probabilities, the probability of collision. \cite{irvine2002geometrical} proposed a geometrical approach for estimating conflict probabilities. \cite{lygeros2002aircraft} incorporated Gaussian winds into a discrete-time probabilistic model while \cite{hu2005aircraft} applied the same technique by assuming they are spatially correlated. \cite{chaloulos2007effect} analyzed how wind affects conflict probability through spatial and temporal correlations. \cite{matsuno2016near} and \cite{matsuno2015stochastic} developed probabilistic models under the assumption of spatial uncertainty while \cite{rodionova2016conflict} evaluated strategic maneuvers for conflict resolution over the North Atlantic oceanic airspace. \\

Initial approaches for solving conflict resolution problems under uncertainty focused on meta-heuristics. \cite{durand1996automatic} and \cite{durand1997optimal} proposed genetic algorithms to solve aircraft conflict resolution problems with heading control and trajectory recovery under speed uncertainties. In the latter, a conflict solver is proposed and its performance is reported realistic air traffic scenarios in France. In addition, \citet{durand2009ant} presented an ant colony optimization approach as a solution for larger conflict resolution problems under uncertainty involving up to 30 aircraft using a limited number of pheromone trails by assigning multiple ants for each aircraft. Recently, \cite{wang2020cooperation} proposed a generic framework based on altitude, speed and heading control and accounting for uncertainties on aircraft positions. The proposed framework combines a memetic algorithm with integer-linear programming formulations. While several works explored the impact of accounting for uncertainty on aircraft trajectories and the associated probability of a loss of separation, only a few studies proposed methods which accounts for aircraft trajectory prediction uncertainty within mathematical programming approaches for the ACRP, as highlighted in a recent survey \citep{pelegrin2020airspace}. \cite{vela2009two} proposed a two-stage stochastic optimization approach to identify the optimal trade-off between deviation costs and conflict probabilities under wind uncertainty. In a first stage, aircraft speeds are determined based on expected costs, and in a second stage manoeuvres necessary to compensate wind modelling errors are identified. A multi-aircraft model for conflict detection based on aircraft dynamics, flight plan and flight management system were proposed by \cite{glover2004multi}, where statistics properties from wind field were extracted from weather data. \cite{rey2016subliminal} examined the behaviour of a deterministic conflict resolution optimization formulation in the presence of random perturbations onto aircraft velocity but did not incorporate uncertainty within the optimization formulation.\\

As highlighted in this review of the literature, most efforts are either focused on probabilistic models to calculate conflict rates and to improve the accuracy of decision-support systems for ATC by providing a more comprehensive model of reality; or on (meta-)heuristic approaches that aim to solve the ACRP under uncertainty. To the best of knowledge, no attempts have been made to solve the ACRP to optimality without any form of discretization while accounting for uncertainty in a robust optimization framework. We address this research gap in this work by considering a simple uncertainty model on aircraft velocities and proposing a scalable robust counterpart formulation.

\section{Aircraft Conflict Resolution Problem}
\label{model}

In this section, we recall the premises of the aircraft conflict resolution problem (ACRP) and present a non-convex formulation for the ACRP. We consider a set of aircraft and assume that all aircraft are initially separated and that their initial positions and velocities are known. This sets the context of the optimization problem of interest: given a set of aircraft in their initial configuration, find least-deviating conflict-free trajectories for all aircraft. To address this problem, we propose a two-dimensional conflict resolution formulation based on heading changes and speed control which are the most widely separation maneuvers in for aircraft conflict resolution \citep{omer2015space,hassan2020mixed}. This formulation is based on the formulations proposed by \cite{rey2017complex} and \cite{dias2020disjunctive}.

\subsection{Separation Conditions}

Consider a set of aircraft $\A$ sharing the same flight level. For each aircraft $i \in \mathcal{A}$, assuming uniform motion laws apply, its position is: $p_i(t) = [x_i(t) = \widehat{x}_i + q_i\hat{v}_i\cos (\widehat{\theta}_{i} + \theta_i)t, y_i(t) = \widehat{y}_i + q_i\hat{v}_i\sin(\widehat{\theta}_{i} + \theta_i)t]^\top$ in which $v_i$ is the speed, $\widehat{x}_i$ and $\widehat{y}_i$ are the initial coordinates of $i$ at the time of optimization, $\widehat{\theta}_{i}$ is its initial heading angle, $\theta_i$ is its deviation angle and $q_i$ is the speed deviation. The relative velocity vector of $i$ and $j$, denoted $v_{ij}$, can be expressed as:
\begin{equation}
    v_{ij} = [\vijx,\vijy]^\top = [\vix - \vjx,\viy - \vjy]^\top,
\end{equation}

where:
\begin{align}
    \vix = q_i\hat{v}_i\sin(\widehat{\theta}_{i} + \theta_i), \\
    \viy = q_i\hat{v}_i\cos(\widehat{\theta}_{i} + \theta_i). 
\end{align}

The relative position of aircraft $i$ and $j$ at time $t$ can be represented as $p_{ij}(t) = p_{i}(t) - p_{j}(t)$. Let $d$ be the horizontal separation norm, typically $d = 5$ NM. Two aircraft $i,j \in \mathcal{A}$ are horizontally separated if and only if: $||p_{ij}(t)|| \ge d$, for all $t \ge 0$. Incorporating these elements into the equation of motion gives:
\begin{subequations}
\begin{align}
& {x}_{i}(t) = \xiz + \vix t,\\
& {y}_{i}(t) = \yiz + \viy t.
\end{align}\label{randommotion}
\end{subequations}

Let $\P$ be the set by the pairs of aircraft, i.e. $\P = \{i \in \A, j \in \A, i <j\}$. For each pair $(i,j) \in \P$, the relative position vector ${p}_{ij}(t)$ is:
\begin{equation}
p_{ij}(t) = [x_{ij}(t),y_{ij}(t)]^\top,
\end{equation}

and the relative velocity vector ${v}_{ij} = [\vijx,\vijy]^\top,$ is:
\begin{subequations}\label{eq:v}
\begin{align}
& \vijx = \qi \viz\cos(\tiz + \ti) - \qj \vjz\cos(\tjz + \tj), \\
& \vijy = \qi \viz\sin(\tiz + \ti) - \qj \vjz\sin(\tjz + \tj).
\end{align}
\end{subequations}

Imposing the separation condition, gives for each pair $(i,j) \in \P$:
\begin{equation}\label{eq:sepcond}
||{p}_{ij}(t)|| \geq d \Rightarrow \sqrt{({x}_i(t) - y_j(t))^2 + ({y}_i(t)- {y}_j(t))^2} \geq d, \quad \forall t \ge 0.
\end{equation}

Let $\xijz = \xiz - \xjz$ and $\yijz = \yiz - \yjz$. Squaring both sides in Eq. \eqref{eq:sepcond}, we obtain:
\begin{equation}
{f}_{ij}(t) \equiv ((\vijx)^2 + (\vijy)^2)t^2 + (2\vijx\xijz + 2\vijy\yijz)t + \xijz^2 + \yijz^2- d^2 \geq 0.
\end{equation}

The function ${f}_{ij}(t)$ is a second-order polynomial in $t$ which is minimal for ${f}^\prime_{ij}(t) = 0$:
\begin{equation}\label{trobust}
{f}^\prime_{ij}(t) = 0 \Rightarrow \tm \equiv -\frac{\xijz\vijx + \yijz\vijy}{(\vijx)^2 + (\vijy)^2}.
\end{equation}

The time instant $\tm$ represents the time of minimal separation of aircraft $i$ and $j$. As noted in several studies \citep{cafieri2017maximizing,cafieri2017mixed,rey2017complex}, if $\tm \leq 0$ then aircraft $i$ and $j$ are diverging and, assuming aircraft are separated at $t=0$, they are thus separating for any $t \geq 0$. Further, substituting $\tm$ in ${f}_{ij}(t)$ yields:
\begin{equation}\label{grobust}
{g}_{ij} \equiv {f}_{ij}(\tm) = (\vijy)^2(\xijz^2 - d^2) + (\vijx)^2(\yijz^2 - d^2) - (2 \xijz\yijz\vijx\vijy).
\end{equation}

Hence, if $g_{ij} \geq 0$ then ${f}_{ij}(t)$ for all $t \geq 0$ and aircraft $i$ and $j$ are separated. Writing  the terms $g_{ij}$ and $\tm$ as functions of the relative velocity variables $\vijx$ and $\vijy$, we obtain the following disjunctive pairwise aircraft separation conditions:
\begin{equation}\label{eq:sepconditions}
||{p}_{ij}(t)|| \geq d, \forall t \geq 0 \Leftrightarrow {g}_{ij}(\vijx,\vijy) \geq 0 \vee \tm(\vijx,\vijy) \leq 0.
\end{equation}

The separation condition \eqref{eq:sepconditions} can be further linearized following the approach described by \cite{rey2017complex} and \cite{dias2020disjunctive}. By alternatively fixing variables $\vijx$ and $\vijy$ and solving the resulting quadratic equations, we can obtain the solution for $g_{ij}(\vijx,\vijy) = 0$. By isolating each variable, we obtain the discriminants:
\begin{equation}\label{discriminant}
\begin{cases}
\Delta_{\vijx} = 4d^2(\vijy)^2(\xijz^2 + \yijz^2 - d^2),\\
\Delta_{\vijy} = 4d^2(\vijx)^2(\xijz^2 + \yijz^2 - d^2).
\end{cases}
\end{equation}

Assuming aircraft are initially separated, then $\xijz^2 + \yijz^2 - d^2 \geq 0$ holds and thus the discriminants are positive and the roots of equation $g(\vijx,\vijy) = 0$ are the lines defined by the system of equations:
\begin{subequations}
\begin{align}
(\hat{y}_{ij}^2 -d^2)\vijx - (\hat{x}_{ij}\hat{y}_{ij} + d\sqrt{\hat{x}_{ij}^2 + \hat{y}_{ij}^2 - d^2})\vijy = 0, \\
(\hat{y}_{ij}^2 -d^2)\vijx - (\hat{x}_{ij}\hat{y}_{ij} - d\sqrt{\hat{x}_{ij}^2 + \hat{y}_{ij}^2 - d^2})\vijy = 0, \\ 
(\hat{x}_{ij}^2 -d^2)\vijy - (\hat{x}_{ij}\hat{y}_{ij} + d\sqrt{\hat{x}_{ij}^2 + \hat{y}_{ij}^2 - d^2})\vijx = 0, \\
(\hat{x}_{ij}^2 -d^2)\vijy- (\hat{x}_{ij}\hat{y}_{ij} - d\sqrt{\hat{x}_{ij}^2 + \hat{y}_{ij}^2 - d^2})\vijx = 0. 
\end{align}
\label{robustmotion}
\end{subequations}

\begin{figure}[t]
\centering
\includegraphics[width=1.0\linewidth]{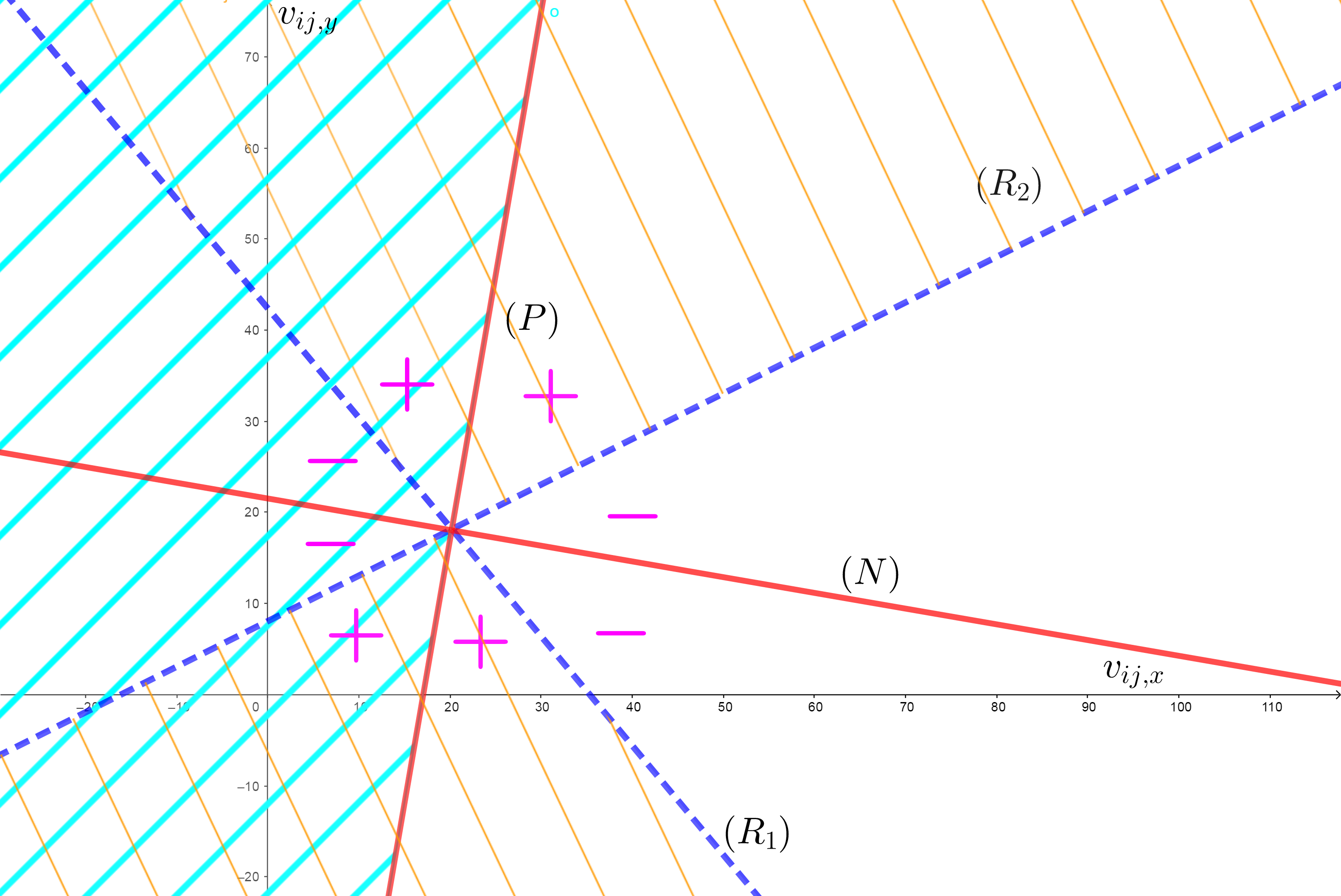}\label{fig:gplotrobust}
\caption{Illustration of a two-aircraft conflict in the plane $\{(\vijx,\vijy) \in \mathbb{R}^2\}$. The red lines represent the lines $P$ and $N$. The dashed blue lines correspond to the linear equations $R_1$ and $R_2$ that are the roots of $g(\vijx,\vijy) = 0$. The sign of $g(\vijx,\vijy)$ is shown by the + and - pink symbols. The hashed orange region represents $g(\vijx,\vijy) \geq 0$. The hashed blue half-plane represents diverging trajectories, i.e. $\tm(\vijx,\vijy) \leq 0$.}
\end{figure}

Considering that if all coefficient in Eq. \eqref{robustmotion} are non-zero, they define two lines, denoted $R_1$ and $R_2$, in the plane $\{(\vijx,\vijy) \in \mathbb{R}^2\}$ and the sign of $g_{ij}(\vijx,\vijy)$ can be characterized based on the position of $(\vijx,\vijy)$ relative to these lines (see Figure \eqref{fig:gplotrobust}). Recall that according to Eq. \eqref{trobust}, the sign of the dot product~$\hat{p}_{ij} \cdot v_{ij}$ indicates aircraft convergence or divergence. Let \eqref{plane} be the equation of the line corresponding to the dot product~$\hat{p}_{ij} \cdot v_{ij}$. 
\begin{equation}\label{plane}\tag{$P$}
\vijx\xijz + \vijy\yijz = 0.
\end{equation}

The line defined by \eqref{plane} splits the plane $\{(\vijx,\vijy) \in \mathbb{R}^2\}$ in two half-planes, each of which representing converging and diverging trajectories, respectively. Consider the line normal to \eqref{plane}, denoted \eqref{normalplane}:
\begin{equation}\label{normalplane}\tag{$N$}
\vijy\hat{x}_{ij} - \vijx\hat{y}_{ij} = 0.
\end{equation}

Recall that any point $(\vijx,\vijy)$ such that $\tm \leq 0$ or $g_{ij}(\vijx,\vijy) \geq 0$ corresponds to a pair of conflict-free trajectories. Hence, the conflict-free region is nonconvex. According to Figure \ref{fig:gznp}, we can observe that the dashed purple lines divide the solution space into two regions: the region hashed in red corresponds to the conflict region $\C$ while the complement of this region represents conflict-free trajectories. We recall the formal definition of the conflict region. 

\begin{definition}[Conflict region, \cite{dias2020disjunctive}]
\label{cr}
Consider a pair of aircraft $(i,j) \in \P$. Let $\C$ be the subset of $\R^2$ defined as:
\begin{equation}
\C = \left\{(\vijx,\vijy) \in \R^2 : \vijx \gamma_{ij}^l - \vijy \phi_{ij}^l \geq 0 \land \vijx \gamma_{ij}^u - \vijy \phi_{ij}^u \leq 0 \right\}.
\end{equation}
$\C$ is the conflict region of $(i,j) \in \P$.
\end{definition}

The conflict region of pair of aircraft represents the set of relative velocity vectors $(\vijx,\vijy)$ which corresponds to conflicts. As depicted in Figure \ref{fig:gznp}, the feasible region is non-convex which that it can be solved via a disjunctive formulation. However, as presented in \cite{rey2017complex} and \cite{dias2020disjunctive}, the set of linear equations described by Eqs. \ref{eq:zd} is equivalent to \eqref{grobust} as detailed by Theorem 1 in \cite{dias2020disjunctive}. In each convex sub-region, the lines delineate the conflict-free region. The expressions of these lines depend on aircraft initial positions, i.e. $\xijz$, $\yijz$. Integer-linear separation conditions with regards to aircraft velocity components can be derived as follows and we model this disjunction using the variable $z_{ij} \in \{0,1\}$ defined as:
\begin{subequations}
\label{eq:zd}
\begin{align}
\vijy \xijz - \vijx \yijz \leq 0, &\quad \text{ if } \z=1, \quad \forall (i,j) \in \P,\\
\vijy \xijz - \vijy \yijz \geq 0, &\quad \text{ if } \z=0, \quad \forall (i,j) \in \P,\\
\vijy \gamma_{ij}^l - \vijx \phi_{ij}^l \leq 0, &\quad \text{ if } \z=1, \quad \forall (i,j) \in \P, \label{eq:sep1}\\
\vijy \gamma_{ij}^u - \vijx \phi_{ij}^u \geq 0, &\quad \text{ if } \z=0, \quad \forall (i,j) \in \P. \label{eq:sep2}
\end{align}
\end{subequations}

To characterize the set of 2D conflict-free trajectories, we examine the relative velocity vector $v_{ij}$ as a function of trajectory control bounds. For each aircraft $i \in \A$, we assume that the speed rate variable is lower bounded by $\qmin$ and upper bounded by $\qmax$, i.e.:
\begin{equation}\label{eq:speedbound}
\qmin \leq \qi \leq \qmax, \qquad \forall i \in \A.
\end{equation}

We assume that the heading deviation is lower bounded by $\tmin$ and upper bounded by $\tmax$, i.e.:
\begin{equation}\label{eq:thetabound}
\tmin \leq \ti \leq \tmax,  \qquad \forall i \in \A.
\end{equation}

To derive lower and upper bounds on relative velocity components $\vijx$ and $\vijy$, we re-arrange Eq. \eqref{eq:v} using trigonometric identities:
\begin{subequations}\label{eq:v2}
\begin{align}
& \vijx = \qi \viz\cos(\tiz)\cos(\ti) - \qi \viz\sin(\tiz)\sin(\ti) - \qj \vjz\cos(\tjz)\cos(\tj) + \qj \vjz\sin(\tjz)\sin(\tj), \\
& \vijy = \qi \viz\sin(\tiz)\cos(\ti) + \qi \viz\cos(\tiz)\sin(\ti) - \qj \vjz\sin(\tjz)\cos(\tj) - \qj \vjz\cos(\tjz)\sin(\tj).
\end{align}
\end{subequations}

Let $\underline{v}_{ij,x},\overline{v}_{ij,x}$ and $\underline{v}_{ij,y},\overline{v}_{ij,y}$ be the lower and upper bounds for $\vijx$ and $\vijy$, respectively. These bounds can be determined using Eq. \eqref{eq:v2} and the bounds on speed and heading control provided in Eqs. \eqref{eq:speedbound} and \eqref{eq:thetabound}. The derived bounds on the relative velocity components can be used to define a box in the plane $\{(\vijx,\vijy) \in \mathbb{R}^2\}$, as noted by \cite{dias2020disjunctive}.

\begin{definition}[Relative velocity box, \cite{dias2020disjunctive}]
\label{rbox}
Consider a pair of aircraft $(i,j) \in \P$. Let $\B$ be the  subset of $\R^2$ defined as
\begin{equation}
\B = \left\{(\vijx,\vijy) \in \mathbb{R}^2 : \underline{{v}}_{ij,x}\leq \vijx \leq \overline{{v}}_{ij,x}, \underline{{v}}_{ij,y} \leq \vijy \leq \overline{{v}}_{ij,y}\right\}.
\end{equation}
$\B$ is the relative velocity box of $(i,j) \in \P$.
\end{definition}

The relative velocity box $\B$ characterize all possible trajectories for the pair $(i,j) \in \P$ based on the available 2D deconfliction resources, i.e. speed and heading controls. To characterize the set of conflict-free trajectories of a pair of aircraft $(i,j) \in \P$, we compare the relative position of the relative velocity box $\B$ with the conflict region of this pair of aircraft. Observe that the conflict region is convex and can be defined by reversing the inequalities \eqref{eq:sep1}-\eqref{eq:sep2} and omitting the disjunction $\z \in \{0,1\}$. 

The objective function is based on the deviation angle and the speed variation. To design the objective function, we introduce a preference weight $w \in \ ]0,1[$ to balance the trade-offs among velocity controls, i.e. speed and heading. Using the definitions above, the nonconvex formulation the ARCP is summarized in Model \ref{mod:robustnonconvex}.

\begin{model}[Nonconvex Formulation]
\label{mod:robustnonconvex}
\begin{subequations}
\begin{align}
&\emph{Minimize} \quad \sum_{i \in \A}(1 - w)(1 - q_i)^2 + w\theta_i^2, \label{eq:obj} \nonumber \\
&\emph{Subject to:}  && \nonumber \\
& \vijx = q_i\hat{v}_i\sin(\hat{\ti} + \ti) - q_j\hat{v}_j\sin(\hat{\tj} + \tj), && \forall (i,j) \in \P,\\
& \vijy = q_i\hat{v}_i\cos(\hat{\ti} + \ti) - q_j\hat{v}_j\cos(\hat{\tj} + \tj), && \forall (i,j) \in \P,\\
&\vijy \hat{x}_{ij} - \vijx \hat{y}_{ij} \leq 0, \quad \text{ if } \z=1, && \forall (i,j) \in \P,\\ \label{c1}
&\vijy \hat{x}_{ij} - \vijx \hat{y}_{ij} \geq 0, \quad \text{ if } \z=0, && \forall (i,j) \in \P, \\ \label{c2}
&\vijy \gamma_{ij}^l - \vijx \phi_{ij}^l \leq 0, \quad \text{ if } \z=1, && \forall (i,j) \in \P, \\ \label{c3}
&\vijx \gamma_{ij}^u - \vijy \phi_{ij}^u \geq 0, \quad \text{ if } \z=0, && \forall (i,j) \in \P, \\ \label{c4}
& \vijx,\vijy \in \B, && \forall (i,j) \in \P,\\
& \z \in \{0,1\}, &&\forall (i,j) \in \P, \\
& \qmin \leq q_i \leq \qmax, && \forall i \in \A,\\
& \tmin \leq \theta_i \leq \tmax, && \forall i \in \A.
\end{align}
\end{subequations}
\end{model}

Model \ref{mod:robustnonconvex} provides a compact formulation for the ACRP with speed and heading control which requires a single binary variable per pair of aircraft. This reformulation does not require the introduction of auxiliary variables and is proved to provide the tightest continuous relaxation for each On/Off constraint. Note that coefficients $\gamma_{ij}^l$, $\phi_{ij}^l$, $\gamma_{ij}^u$ and $\phi_{ij}^u$ (present in \eqref{eq:zd}) can be pre-processed based on the sign of $\hat{x}_{ij}$ and $\hat{y}_{ij}$. For implementation details, a fully reproducible formulation can be found at: \small \url{https://github.com/acrp-lib/acrp-lib}\normalsize.\newline

Model \ref{mod:robustnonconvex} is a deterministic formulation of the ACRP with speed and heading control. We next extend this formulation using robust optimization to account for the impact of aircraft trajectory prediction with the conflict resolution procedure. 

\section{Robust Aircraft Conflict Resolution Problem}
\label{robs}

In this section, we introduce formulations for the robust aircraft conflict resolution problem. We first define the uncertainty model before discussing how it can be incorporated within the ACRP. We then propose a tractacble robust counterpart formulation for the robust ACRP.

\subsection{Uncertainty Model}

We assume that each aircraft has a source of randomness and this affects its current velocity and position. Let $\ei = [\eix,\eiy]^\top$ be a vector of random variables representing the uncertainty on the velocity components of aircraft $i \in \A$. We next use this vector of random variables to define aircraft-based uncertainty sets.

\begin{definition}[Uncertainty set of aircraft]
The uncertainty set of aircraft $i \in \A$, is defined as:
\begin{equation}
\U_i \equiv \{\ei \in \mathbb{R}^2\ |\  -\heix \leq \eix \leq \heix, -\heiy \leq \eiy \leq \heiy \},
\end{equation}

where $\heix \geq 0$ and $\heiy \geq 0$ represent the maximum perturbations on the velocity components $\vix$ and $\viy$, respectively, of aircraft $i$.
\end{definition}

We denote $[\tvix,\tviy]^\top$ the vector of random aircraft velocity components where the random variables $\tvix$ and $\tviy$ take values in $\tvix \in [-\vix(1 + \heix), \vix(1 + \heix)]$ and $\tviy \in [-\viy(1 + \heiy), \viy(1 + \heiy)]$, respectively. Accordingly, for each pair of aircraft $(i,j) \in \P$, the random relative velocity components $\tvx$ and $\tvy$ are:
\begin{subequations}
\label{tvxtvy}
\begin{align}
& \tvx = \vix(1+\eix) - \vjx(1 + \ejx) = \vijx + \vix\eix - \vjx\ejx,  \\
& \tvy = \viy(1+\eiy) - \vjy(1 + \ejy) = \vijy + \viy\eiy - \vjy\eiy. 
\end{align}
\end{subequations}

Let $\underline{\tilde{v}}_{ij,x},\overline{\tilde{v}}_{ij,x}$ and $\underline{\tilde{v}}_{ij,y},\overline{\tilde{v}}_{ij,y}$ be the lower and upper bounds for $\tvx$ and $\tvy$, respectively. These bounds can be determined using Eq. \eqref{tvxtvy} and the bounds on speed and heading control provided in Eqs. \eqref{eq:speedbound} and \eqref{eq:thetabound}. The derived bounds on the random relative velocity components can be used to define the random relative velocity box which characterizes all possible trajectories for the pair $(i,j) \in \P$ under uncertainty.

\begin{definition}[Random relative velocity box]
\label{box}
Consider a pair of aircraft $(i,j) \in \P$. Let $\U_i$ and $\U_j$ be the uncertainty sets of aircraft $i$ and $j$, respectively. Let $\Br(\U_i,\U_j)$ be the parametric subset of $\R^2$ defined as
\begin{equation}
\Br(\U_i,\U_j) \equiv \left\{(\tvx,\tvy) \in \mathbb{R}^2 : \underline{\tilde{v}}_{ij,x}\leq \tvx \leq \overline{\tilde{v}}_{ij,x}, \underline{\tilde{v}}_{ij,y} \leq \tvy \leq \overline{\tilde{v}}_{ij,y}\right\}.
\end{equation}
$\Br(\U_i,\U_j)$ is the random relative velocity box of $(i,j) \in \P$ under the uncertainty sets $\U_i$ and $\U_j$.
\end{definition}

\begin{figure}[!ht]
\centering
\includegraphics[width=0.8\linewidth]{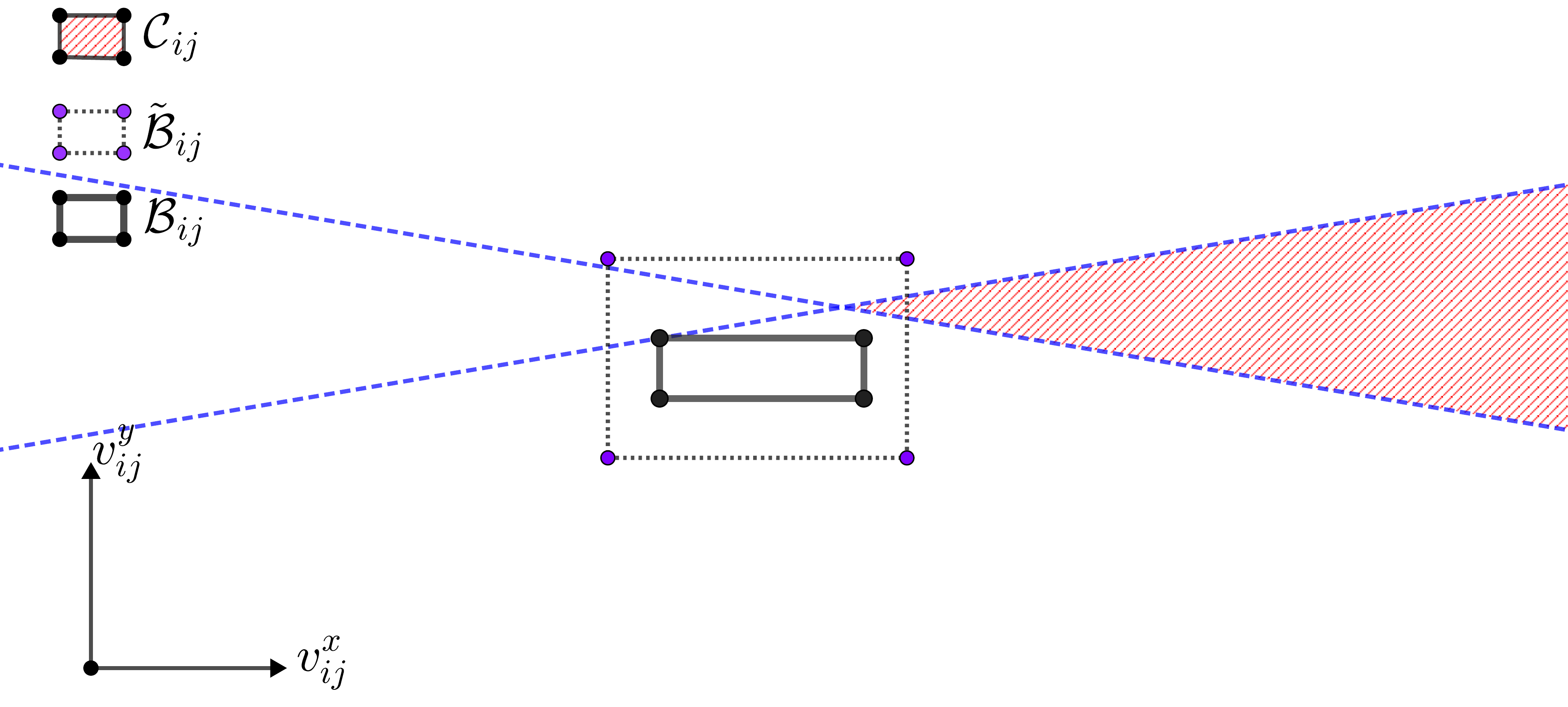}
\caption{Illustration of a two-aircraft conflict in the plane $\{(\vijx,\vijy) \in \mathbb{R}^2\}$. The inner box with black lines corresponds to the velocity bounds $\B$ in the deterministic scenario while the box with purple dots $\Br(\U_i,\U_j)$ is the random relative velocity box. The region is hashed in red corresponds to the conflict region $\C$. If the uncertainty sets $\U_i$ and $\U_j$ of aircraft $i$ and $j$ are empty, i.e. aircraft trajectories are deterministic, then aircraft $i$ and $j$ are conflict-free. In turn, if the uncertainty sets $\U_i$ and $\U_j$ are such that the random relative velocity $\Br$ intersects with the conflict region $\C$, then there exists a risk of conflict.}
\label{fig:gznp}
\end{figure}

The random relative velocity box is illustrated in Figure \ref{fig:gznp} for a two-aircraft conflict. To ensure that a pair of aircraft $(i,j) \in \P$ is separated for any realization of the random variables $\ei \in \U_i$ and $\ej \in \U_j$, we redefine \eqref{c1}-\eqref{c4} using the random velocity variables $\tvx$ and $\tvy$:
\begin{subequations}
\label{eq:zdrs}
\begin{align}
\tvy \xijz - \tvx \yijz \leq 0, &\quad \text{ if } \z=1, \quad \forall (i,j) \in \P, \tag{$N_1$}\label{eq:n1}\\
\tvy \xijz - \tvx \yijz \geq 0, &\quad \text{ if } \z=0, \quad \forall (i,j) \in \P, \tag{$N_0$}\label{eq:n0}\\
\tvy \gamma_{ij}^l - \tvx \phi_{ij}^l \leq 0, &\quad \text{ if } \z=1, \quad \forall (i,j) \in \P, \tag{$S_1$}\label{eq:s1}\\
\tvy \gamma_{ij}^u - \tvx \phi_{ij}^u \geq 0, &\quad \text{ if } \z=0, \quad \forall (i,j) \in \P. \tag{$S_0$}\label{eq:s0}
\end{align}
\end{subequations}

% \begin{subequations}
% \label{eq:zdrs}
% \begin{align}
% \tvy \xijz - \tvx \yijz \leq 0, &\quad \text{ if } \z=1, \quad \forall (i,j) \in \P, \forall \ei \in \U_i, \ej \in \U_j, \tag{$N_1$}\\
% \tvy \xijz - \tvx \yijz \geq 0, &\quad \text{ if } \z=0, \quad \forall (i,j) \in \P, \forall \ei \in \U_i, \ej \in \U_j, \tag{$N_0$}\\
% \tvy \gamma_{ij}^l - \tvx \phi_{ij}^l \leq 0, &\quad \text{ if } \z=1, \quad \forall (i,j) \in \P, \forall \ei \in \U_i, \ej \in \U_j, \tag{$S_1$}\label{eq:r1}\\
% \tvy \gamma_{ij}^u - \tvx \phi_{ij}^u \geq 0, &\quad \text{ if } \z=0, \quad \forall (i,j) \in \P, \forall \ei \in \U_i, \ej \in \U_j. \tag{$S_0$}\label{eq:r2}
% \end{align}
% \end{subequations}

To find robust aircraft trajectories under uncertainty sets $\U_i$ and $\U_j$, we require that the pairwise separation constraints \eqref{eq:n0}, \eqref{eq:n1}, \eqref{eq:s1} and \eqref{eq:s0} hold for any $(\tvx,\tvy) \in \Br(\U_i,\U_j)$. We next use state-of-the-art approaches in robust optimization to integrate these constraints into a robust counterpart formulation.

\subsection{Robust Counterpart Formulation}
\label{robus}

%The state-of-art approaches to handle models with uncertainty sets and robustness are highlighted by \cite{gorissen2015practical} where the authors describe different uncertainty sets and how to apply them directly into our formulation. This strategy provides a general framework based on the definition of the dual problem and the consequent equivalence between dual and primal. On the other hand, \cite{bertsimas2004price}, which is also based on primal-dual equivalence, proposed an alternative approach that allows controlling the robustness level via a single parameter.

The constraints \eqref{eq:n0}, \eqref{eq:n1}, \eqref{eq:s1} and \eqref{eq:s0} are function of the random variables $\tvx$ and $\tvy$. We use the approach of \cite{bertsimas2004price} to reformulate these robust separation constraints as integer-linear constraints with regards to deterministic relative  velocity variables $\vijx$ and $\vijy$. The constraints \eqref{eq:n0}, \eqref{eq:n1}, \eqref{eq:s1} and \eqref{eq:s0} are of the form $a\tvx + b\tvy \leq 0$ (omitting the disjunction), and can be rearranged by separating deterministic and random elements as follows:
\begin{equation}
a(\vix - \vjx) - b (\viy - \vjy) + a \vix\eix - a \vjx\ejx - b \viy\eiy + b \vjy\ejy \leq 0.
\label{zd1exp}    
\end{equation}

We introduce new variables $\nu_{i,x} \geq 0$ and $\nu_{i,y} \geq 0$ for each $i \in \A$ through constraints \eqref{new_mu} to impose artificial bounds on aircraft velocity components $\vix$ and $\viy$:
\begin{subequations}\label{new_mu}
\begin{align}
-\nu_{i,x} \leq \vix \leq \nu_{i,x}, && \forall i \in \A,\\
-\nu_{i,y} \leq \viy \leq \nu_{i,y}, && \forall i \in \A.
\end{align}
\end{subequations}

Let $\RS \equiv \{N_1, N_0, S_1, S_0\}$ be a set of indices corresponding to constraints \eqref{eq:n0}, \eqref{eq:n1}, \eqref{eq:s1}, \eqref{eq:s0}, respectively. Further, let $\alpha_k = 1$ if $k=N_1$ or $k=S_1$, and let $\alpha_k = 0$ if $k=N_0$ or $k=S_0$. Let $\Gamma$ be a real parameter that takes values in the range [0,4], where the upper bound is given by the number of decision variables in constraint \eqref{zd1exp}. The parameter $\Gamma$ determines the level of robustness for each robust separation constraint. To link the level of robustness $\Gamma$ with each robust separation constraint $k \in \RS$ and each aircraft pair $(i,j) \in \P$, we introduce real decision variables $\psi_{ij}^k \geq 0$. Further, each constraint of the form \eqref{zd1exp} involves four decision variables, hence for each constraint $k \in \RS$ and $(i,j) \in \P$, we introduce associated variables $\rho^{l,k}_{ij,x} \geq 0$ and $\rho^{l,k}_{ij,y} \geq 0$ for $l \in \{i,j\}$. These artificial variables reflect the consumption of robustness resources of each velocity variable in \eqref{zd1exp}. The following constraints link variables $\psi_{ij}^k$ $\rho^{l,k}_{ij,x}$ and $\rho^{l,k}_{ij,x}$ with $\nu^x_i$ and $\nu^x_j$:
\begin{subequations}
    \label{new_rho}
	\begin{align}
	\psi_{ij}^k +\rho^{l,k}_{ij,x} \geq \nu_{i,x}\heix, \quad \text{ if } \z=\alpha_k, && \forall (i,j) \in \P, \forall l \in \{i,j\}, \forall k \in \RS,\\
	\psi_{ij}^k +\rho^{l,k}_{ij,y} \geq \nu_{i,y}\heiy, \quad \text{ if } \z=\alpha_k, && \forall (i,j) \in \P, \forall l \in \{i,j\}, \forall k \in \RS.
	%\rho^{l,k}_{i,j,x},\rho^{l,k}_{i,j,x} \geq 0, \quad \text{ if } \z=\alpha_k, && \forall (i,j) \in \P, \forall l \in \{i,j\},\forall k \in RS \\
	%\psi_{i,j}^k \geq 0,\quad \text{ if } \z=\alpha_k, && \forall (i,j) \in \P,\forall k \in RS
	\end{align}
\end{subequations}

Constraints of the form  \eqref{zd1exp} can then be rewritten as:
\begin{equation}
a(\vix - \vjx) - b (\viy - \vjy) + \psi_{ij}^k\Gamma +  \sum_{l \in \{i,j\}} (\rho^{l,k}_{ij,x} + \rho^{l,k}_{ij,y}) \leq 0, \text{ if } \z=\alpha_k,  \forall (i,j) \in \P,\forall k \in \RS.    
\end{equation}

Combining the above constraints and artificial real variables, we obtain a tractable formulation of the robust separation constraints.

\begin{prop}
The set robust separation constraints \eqref{eq:n0}, \eqref{eq:n1}, \eqref{eq:s1} and \eqref{eq:s0} are equivalent to the following set of integer-linear constraints and real variables:
\begin{subequations}
\label{vrobust2}
\begin{align}
& a(\vix - \vjx) - b (\viy - \vjy) \notag\\
& + \psi_{ij}^k\Gamma + \sum_{l \in \{i,j\}} (\rho^{l,k}_{ij,x} + \rho^{l,k}_{ij,y}) \leq 0,\quad &&\text{ if } \z=\alpha_k, && \forall (i,j) \in \P,\forall k \in \RS,\\    
&\psi_{ij}^k +\rho^{l,k}_{ij,x} \geq \nu_{i,x}\heix, \quad &&\text{ if } \z=\alpha_k, && \forall (i,j) \in \P, \forall l \in \{i,j\}, \forall k \in \RS,\\
&\psi_{ij}^k +\rho^{l,k}_{ij,y} \geq \nu_{i,y}\heiy, \quad &&\text{ if } \z=\alpha_k, && \forall (i,j) \in \P, \forall l \in \{i,j\}, \forall k \in \RS,\\
&\rho^{l,k}_{ij,x},\rho^{l,k}_{ij,x} \geq 0, \quad &&\text{ if } \z=\alpha_k, && \forall (i,j) \in \P, \forall l \in \{i,j\},\forall k \in \RS, \\
&\psi_{ij}^k \geq 0,\quad &&\text{ if } \z=\alpha_k, && \forall (i,j) \in \P,\forall k \in \RS,\\
&-\nu_{i,x} \leq \vix \leq \nu_{i,x}, && && \forall i \in \A,\\
&-\nu_{i,y} \leq \viy \leq \nu_{i,y}, && && \forall i \in \A,\\
&\nu_{i,x},\nu_{i,y} \geq 0, && && \forall i \in \A.
\end{align}
\end{subequations}
\label{robust}
\end{prop}

The proof of Proposition \ref{robust} follows from Theorem 1 of \cite{bertsimas2004price}. Proposition \ref{robust} establishes the equivalency between constraints \eqref{eq:n0}, \eqref{eq:n1}, \eqref{eq:s1} and \eqref{eq:s0}, which are expressed in terms of random relative velocity variables $\tvx$ and $\tvy$, and a tractable integer-linear reformulation of these constraints which uses additional real variables. Using this reformulation, we propose the following robust counterpart formulation of the ACRP.

\begin{model}[Robust Counterpart Formulation of the ACRP]
\label{mod:gammalinear}
\begin{subequations}
\begin{align}
&\emph{Minimize} \quad \sum_{i \in \mathcal{A}} (1 - w)(1 - q_i)^2 + w\theta_i^2, \nonumber \\
&\emph{Subject to:}  && \nonumber \\
&\emph{Motion Equations} \quad \eqref{eq:v}, \nonumber\\
&\emph{Robust Separation Constraints and Variables} \quad \eqref{vrobust2}, \nonumber \\
&\emph{Speed and Heading Control Constraints} \quad  \eqref{eq:speedbound}, \eqref{eq:thetabound}, \nonumber\\
& \vijx,\vijy \in \B, &&\forall (i,j) \in \P, \\
& \z \in \{0,1\}, &&\forall (i,j) \in \P, \\
& \qmin \leq q_i \leq \qmax, && \forall i \in \A,\\
& \tmin \leq \theta_i \leq \tmax, && \forall i \in \A.
%& \nu^k_{i,x},\nu^k_{i,y} \geq 0, && \forall i \in \A, && \forall k \in RS\\
%& \rho^k_{i,x},\rho^k_{i,y} \geq 0, && \forall i \in \A, && \forall k \in RS.
\end{align}
\end{subequations}
\end{model}

\subsection{Solution Method for the Robust ACRP}
\label{rcnf}

The formulation presented in Model \ref{mod:gammalinear} is nonconvex due to trigonometric functions and non-linear components. To solve the robust optimization problem represented by Model \ref{mod:gammalinear}, we adopt the approach proposed by \cite{dias2020disjunctive}. This approach solves the deterministic ACRP by using the so-called complex number formulation of the ACRP \citep{rey2017complex} and embedding it within a cut generation algorithm. We next recall the main elements of the complex number formulation for the ACRP and outline how we adapt this formulation for the robust ACRP.

We use artificial variables $\dx$ and $\dy$ defined as: $\dx \equiv q_i \cos(\theta_i)$ and $\dy \equiv q_i\sin(\theta_i)$ to reformulate the expression of aircraft velocities. Specifically, we obtain:
\begin{subequations}
\begin{align}
&\vix =  \dx\viz\cos(\tiz) - \dy\viz\sin(\tiz), \\
&\viy =  \dy\viz\cos(\tiz) - \dx\viz\sin(\tiz).
\end{align}\label{motiondxdy}
\end{subequations}

Bounds on variables $\dx$ and $\dy$ can be derived as:
\begin{subequations}\label{eq:boundsdxdy}
\begin{align}
& \qmin\cos(\max\{|\tmin|,|\tmax|\}) \leq \dx \leq \qmax, && \forall i \in \A, \label{eq:bounddx}\\
& \qmax\sin(\tmin) \leq \dy \leq \qmax\sin(\tmax), && \forall i \in \A. \label{eq:bounddy}
\end{align}
\end{subequations}

The speed control constraint \eqref{eq:speedbound} can be reformulated in quadratic form as:
\begin{subequations}
\label{eq:boundq}
\begin{align}
& \qmin^2 \leq \dx^2 + \dy^2, && \forall i \in \A,\label{eq:lboundq}\\
& \qmax^2 \geq \dx^2 + \dy^2, && \forall i \in \A.\label{eq:uboundq}
\end{align}
\end{subequations}

The heading control constraint \eqref{eq:thetabound} can be reformulated in linear form as:
\begin{subequations}
\label{eq:boundtheta}
\begin{align}
& \dy \geq \dx\tan(\tmin), && \forall i \in \A,\\
& \dy \leq \dx\tan(\tmax), && \forall i \in \A.
\end{align}
\end{subequations}

The objective function of the complex number reformulation of the ACRP is:
\begin{equation}\label{eq:robustobj}
\text{Minimize} \sum_{i \in \A} (1 - w)(1 - \dx)^2 + w\dy^2.
\end{equation}

Objective function \eqref{eq:robustobj} is an approximation of the one used in Model \eqref{mod:gammalinear} which is minimal for $q_i=1$ and $\theta_i = 0$, thus preserving the optimality conditions of the original objective function \eqref{eq:obj}, as detailed in \cite{dias2020disjunctive}. 

The complex number reformulation of the robust is then constructed by substituting the original motion equations \eqref{eq:v} with the reformulated aircraft velocities constraints \eqref{motiondxdy}; replacing the speed and heading control constraints and variables bounds with \eqref{eq:boundq}, \eqref{eq:boundtheta} and \eqref{eq:boundsdxdy}, respectively; and substituting the original objective function \eqref{eq:obj} with \eqref{eq:robustobj}. Observe that the robust separation constraints \eqref{vrobust2} remain unchanged since aircraft velocity variables $\vix$ and $\viy$ are linked to variables $\dx$ and $\dy$ via constraints \eqref{motiondxdy}. We refer to the resulting formulation as the robust complex number formulation and use Algorithm 2 from \cite{dias2020disjunctive} to solve this problem. This solution method initially relaxes the nonconvex quadratic constraint in \eqref{eq:boundq} and uses a cut generation procedure to iteratively solve the relaxed problems until global optimality is reached.

\section{Numerical Experiments}
\label{num}

We first introduce the experimental framework used to test the proposed mixed-integer formulation for the robust problem in Section \ref{fixed}. We then present a detailed  analysis of four instances of the ACRP in Section \ref{ed}. We explore the computational performance of the proposed approach and conduct sensitivity analyzes on the level of robustness and the size of aircraft uncertainty sets in Section \ref{fixedsets} and Section \ref{fixedgamma}, respectively. We attempt to further analyze the behavior of the robust ACRP by examining specific sets of instances in Section \ref{fease}.

%%%%%%%%%%%%%%%%%%%%%%%%%%%%%%%
\subsection{Experiments Design}
\label{ed}

\begin{figure}[ht]	
	\centering	
		\subfloat[CP-7: circle problem with 7 aircraft]{\includegraphics[width=0.4\textwidth]{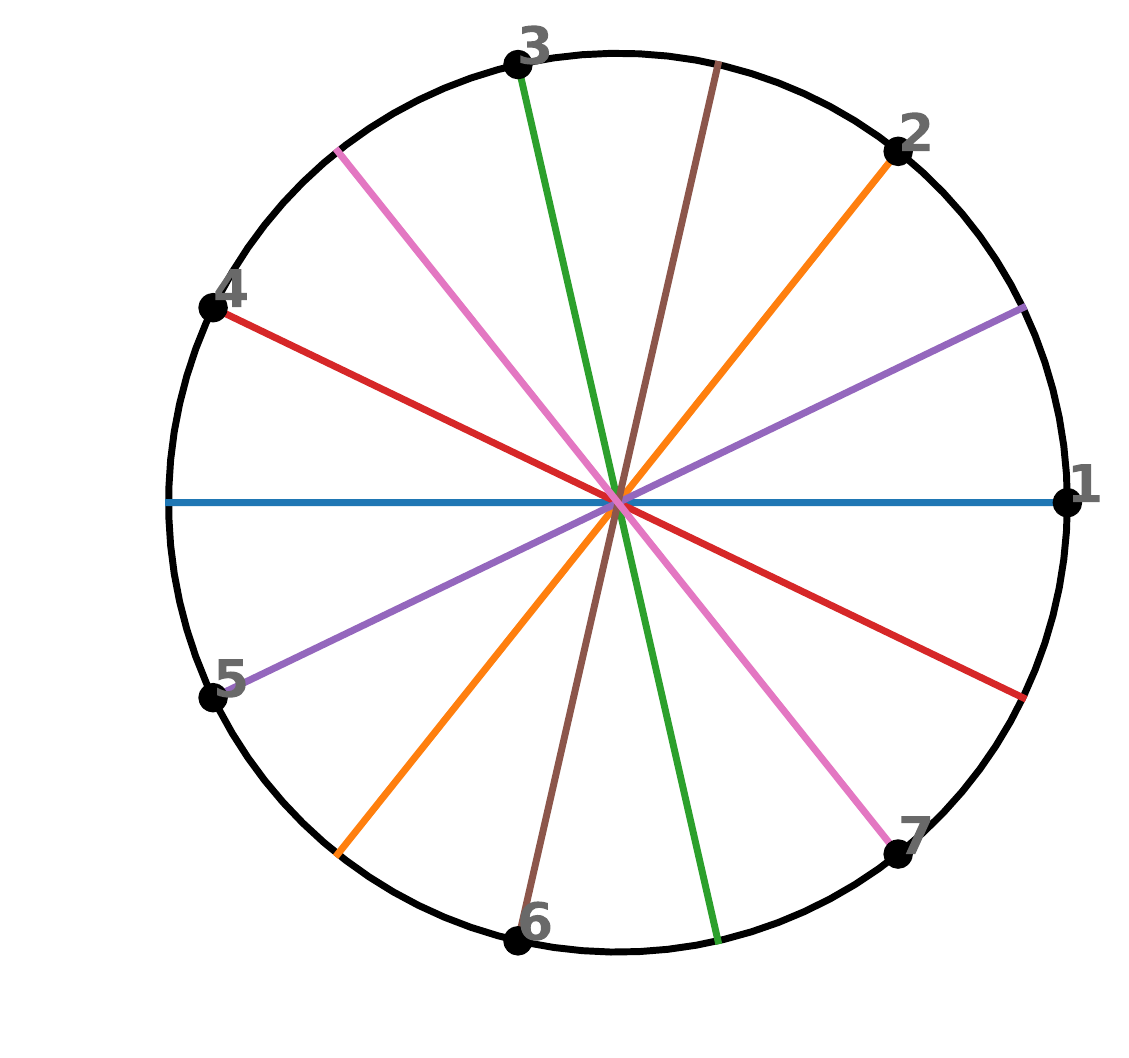}\label{fig:cp}} %\hspace*{2cm}
		\hspace{2cm}
		\subfloat[RCP-10: a random circle problem with 10 aircraft]{\includegraphics[width=0.4\textwidth]{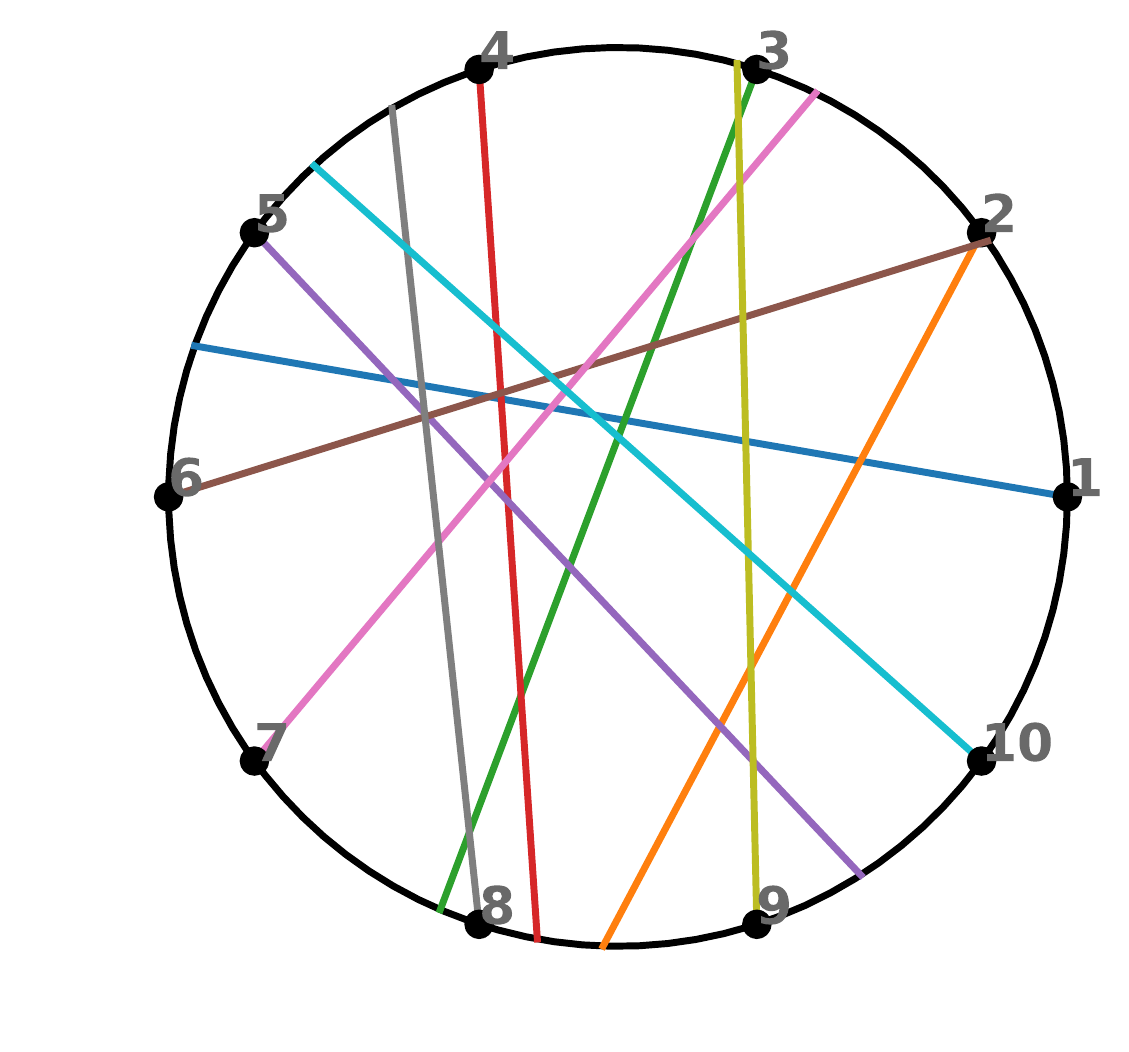} \label{fig:rcp10-5}}\\
	\caption{Example of benchmarking instances for the Circle Problem (CP) and Random Circle Problem (RCP). Aircraft initial positions are represented with black dots.}
	\label{case1}
\end{figure}

We test the performance of the proposed approach using two benchmarking problems from the literature: the Circle Problem (CP) and the Random Circle Problem (RCP). The two types of benchmark instances are illustrated in Figure \ref{case1}. The CP instances consist of a set of aircraft uniformly positioned on a circle and heading towards its centre. Aircraft speeds are assumed to be identical, hence the problem is highly symmetric (see Fig. \ref{fig:cp}). The CP instances are notoriously difficult due to the geometry of aircraft initial configuration. To break the symmetry of CP instances, \citet{vanaret2012benchmarking} introduced the RCP which builds on the same framework, but aircraft initial speeds and headings are randomly deviated within specified ranges to create randomized instances (see Fig. \ref{fig:rcp10-5}). Both benchmarking instances have been widely used in CD\&R algorithms in the literature \citep{durand2009ant,cafieri2017mixed,cafieri2017maximizing,rey2017complex}. CP and RCP instances are named CP-N and RCP-N-ID, respectively, where N is the total number of aircraft and ID is the identifier of the corresponding RCP instance with N aircraft. 
 
In all experiments, we use a circle of radius 200NM. For CP instances, all aircraft have an initial speed of 500 NM/h. For RCP instances, aircraft initial speeds are randomly chosen in the range 486-594 NM/h and their initial headings are deviated from a radial trajectory (i.e. towards the centre of the circle) by adding a randomly chosen angle between $-\frac{\pi}{6}$ and $+\frac{\pi}{6}$. 

We report numerical results for problems with a subliminal speed control range of $[-6\%, +3\%]$ \citep{bonini2009erasmus}. We consider a heading control range of $[-30^\circ,+30^\circ]$ as commonly used in the literature \citep{cafieri2017mixed,rey2017complex}. We assume that the maximum uncertainty on aircraft velocity components is uniform across directions and aircraft, and define $\bar{\epsilon} = \bar{\epsilon}^x_i = \bar{\epsilon}^y_i$ for all aircraft $i \in \A$. In our experiments we consider the following maximum levels of uncertainty: $\bar{\epsilon} = 2.5\%$, $5\%$, $7.5\%$ and $10\%$, and we consider varying levels of robustness $\Gamma = 0$, $1$, $2$, $3$ and $4$. The proposed robust formulations are compared against its deterministic counterpart, which corresponds to $\Gamma=0$ and/or $\bar{\epsilon} = 0$.

All problems are solved with a relative optimality gap of $1\%$ and a time limit of 10 minutes. All models are implemented using \Python on a personal computer with 16 GB of RAM and an Intel i7 processor at 2.9GHz, solved with \Cplex v12.10 \citep{cplex2009v12} API for \Python using default options. For reproducibility purposes, all instances used are available at \small \url{https://github.com/acrp-lib/acrp-lib}.\normalsize

%%%%%%%%%%%%%%%%%%%%%%%%%%%%%%%%%%%%
\subsection{Analysis of the robust ACRP}
\label{fixed}

\begin{figure}
	\centering\vspace{-2pt}
	\subfloat[CP-10]{\includegraphics[width=0.4\textwidth]{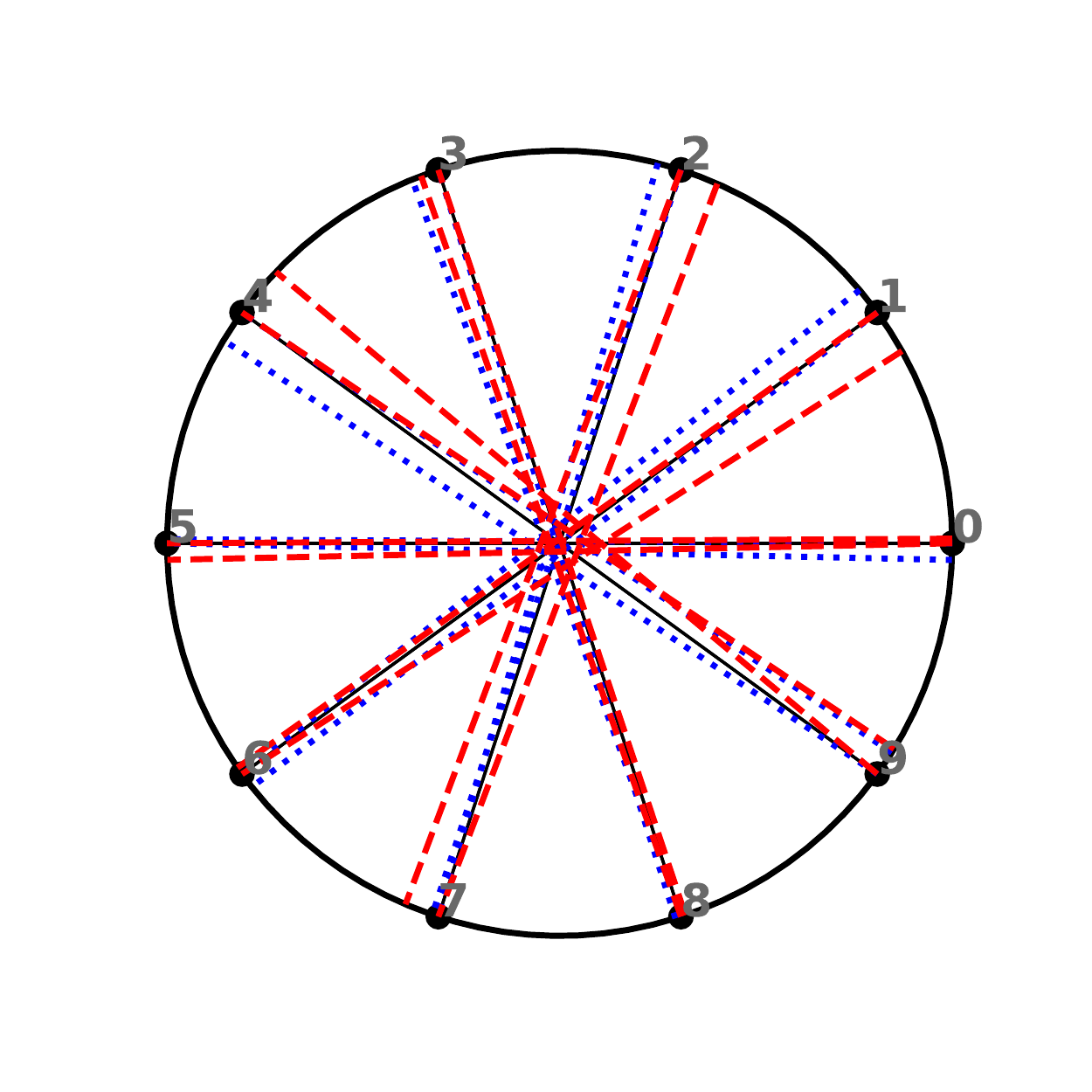}\label{cp10}} 
	\subfloat[RCP-10-1]{\includegraphics[width=0.4\textwidth]{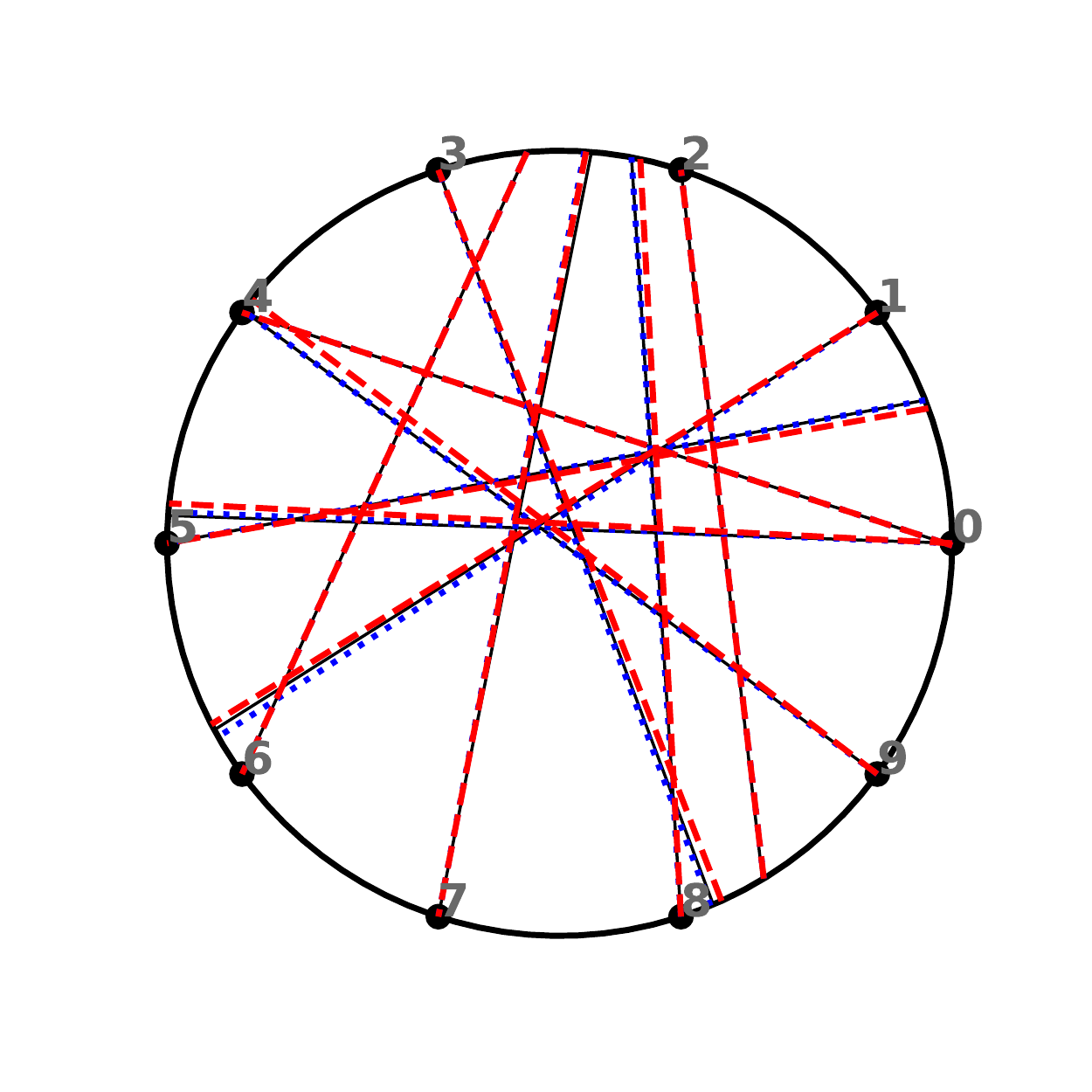}\label{rcp1}}\\
	\subfloat[RCP-20-1]{\includegraphics[width=0.4\textwidth]{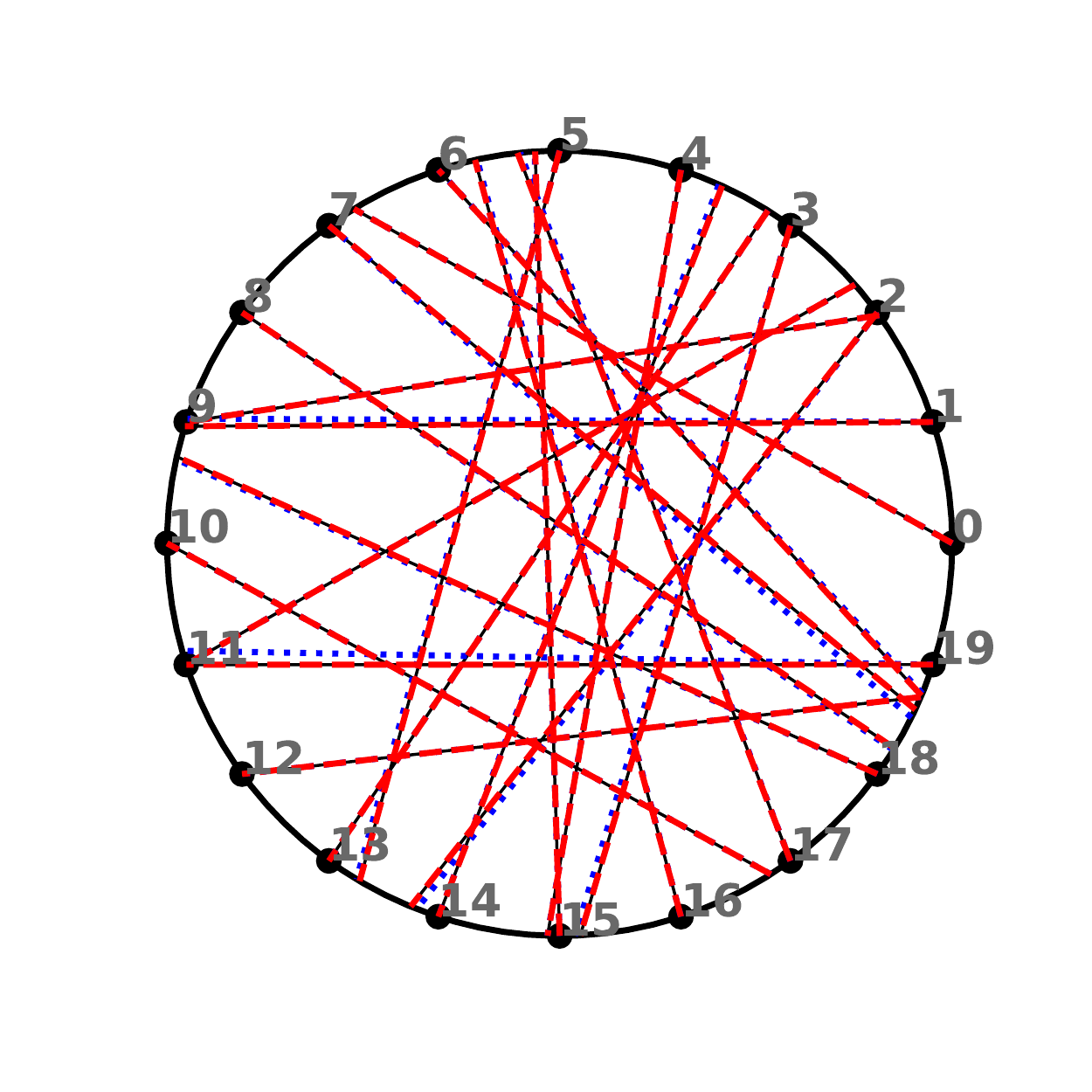}\label{rcp21}} 
	\subfloat[RCP-30-1]{\includegraphics[width=0.4\textwidth]{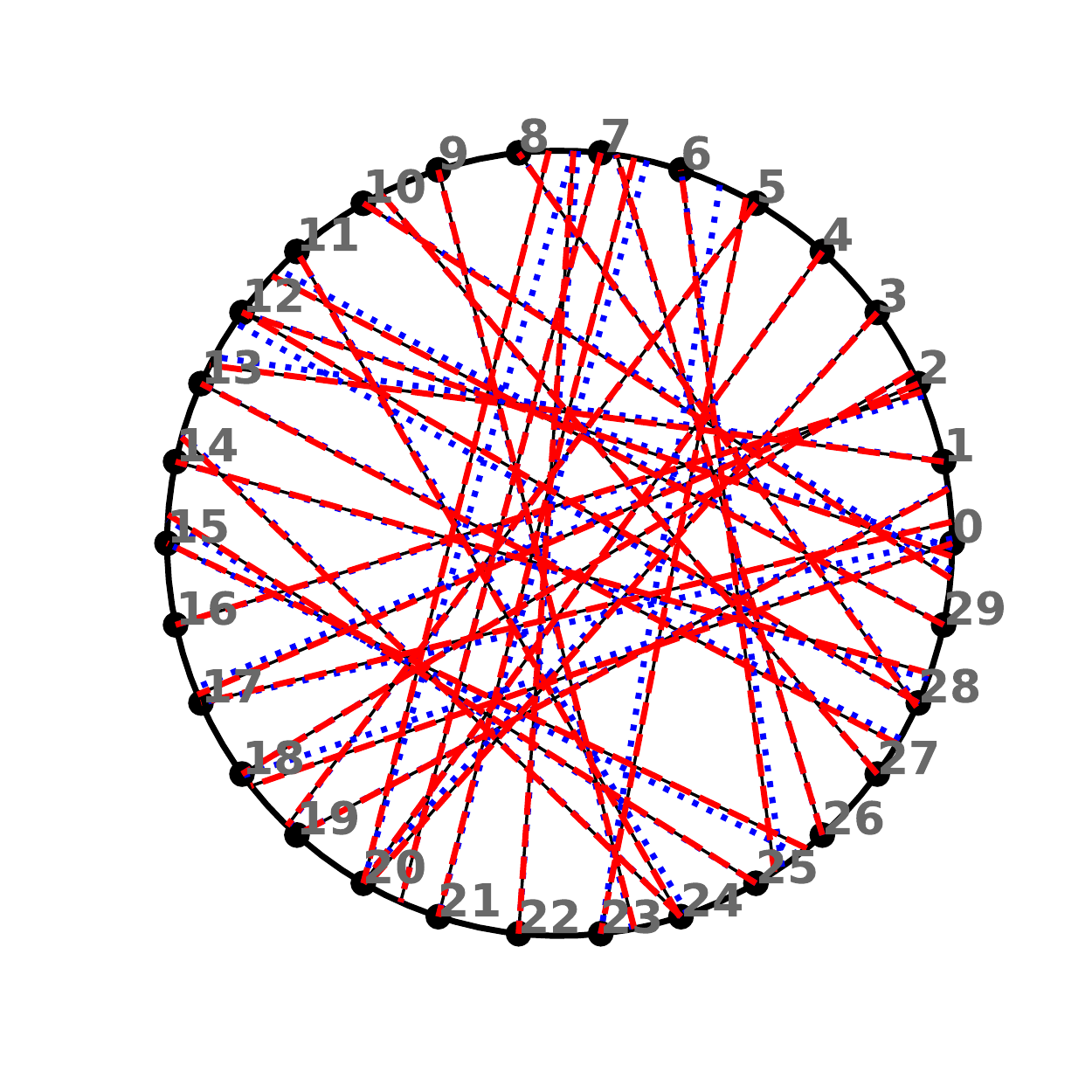}\label{rcp31}}\\
	\caption{Illustration of aircraft optimal trajectories using the robust complex number formulation with a maximum uncertainty of $\bar{\epsilon} = 5\%$. \textcolor{black}{Black} lines represent aircraft initial trajectories. \textcolor{blue}{Blue} lines represent optimal trajectories obtained using $\Gamma = 0$. \textcolor{red}{Red} lines represent avoidance trajectories obtained using $\Gamma$ = 4.}
	\label{ilus} 
\end{figure}

To analyze and illustrate the behavior of the proposed formulation for the robust ACRP, we focus on four instances: CP-1, RCP-10-1, RCP-20-1 and RCP-30-1. We set the maximum uncertainty to $\bar{\epsilon} = 5\%$ and we plot the optimal solutions obtained using the robust complex number formulation as described in section \ref{rcnf} in Figure \ref{ilus}: black lines represent aircraft nominal trajectories, blue lines represent aircraft optimal trajectories using $\Gamma=0$ (which corresponds to the deterministic case), and green lines represent optimal trajectories using $\Gamma = 4$ (which corresponds to the maximum level of robustness).

\begin{figure}[H]
	\centering\vspace{-2pt}
	\subfloat[CP-10]{\includegraphics[width=0.4\textwidth]{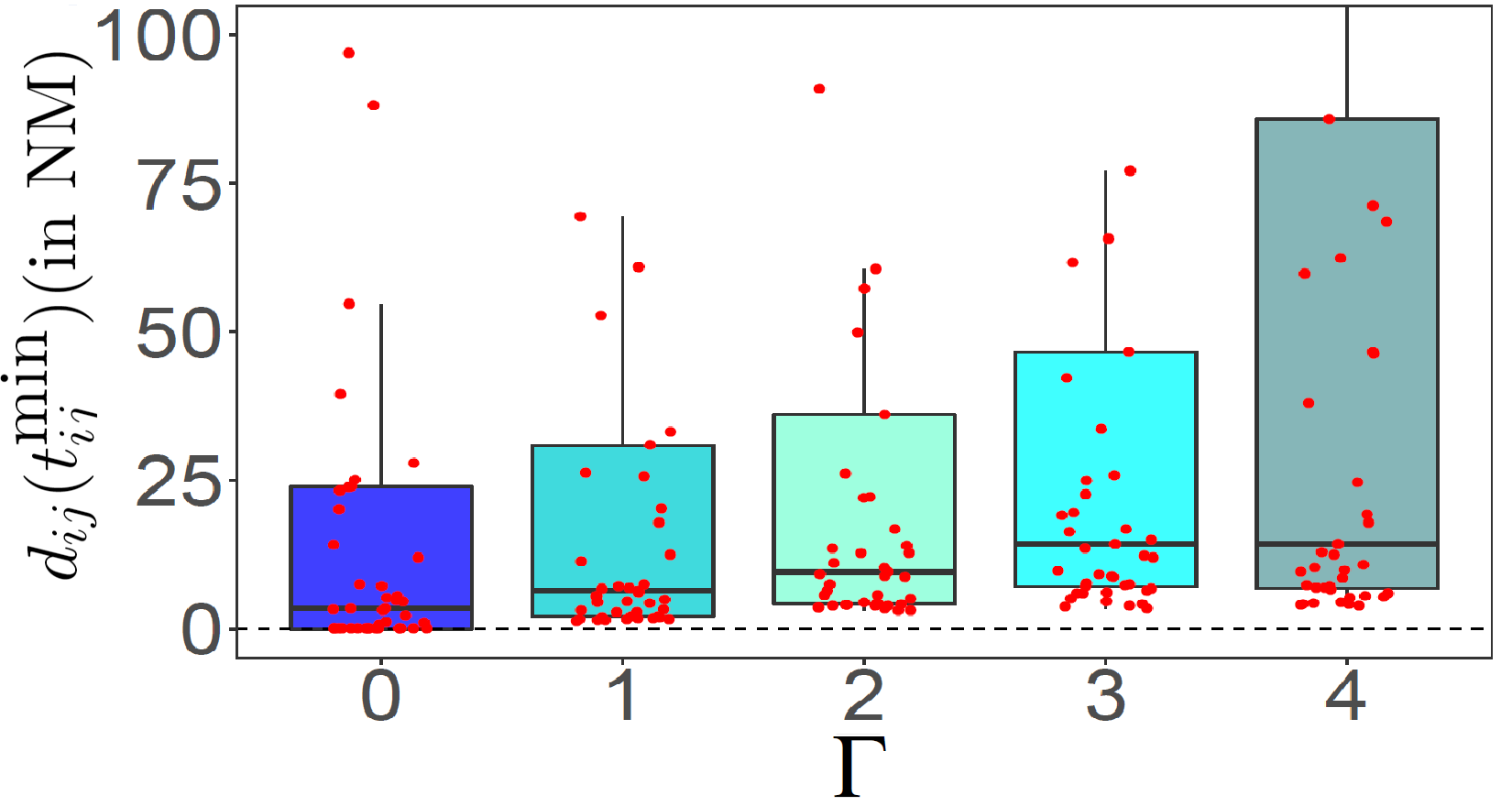}\label{cp_gamma}} \qquad
	\subfloat[RCP-10-1]{\includegraphics[width=0.4\textwidth]{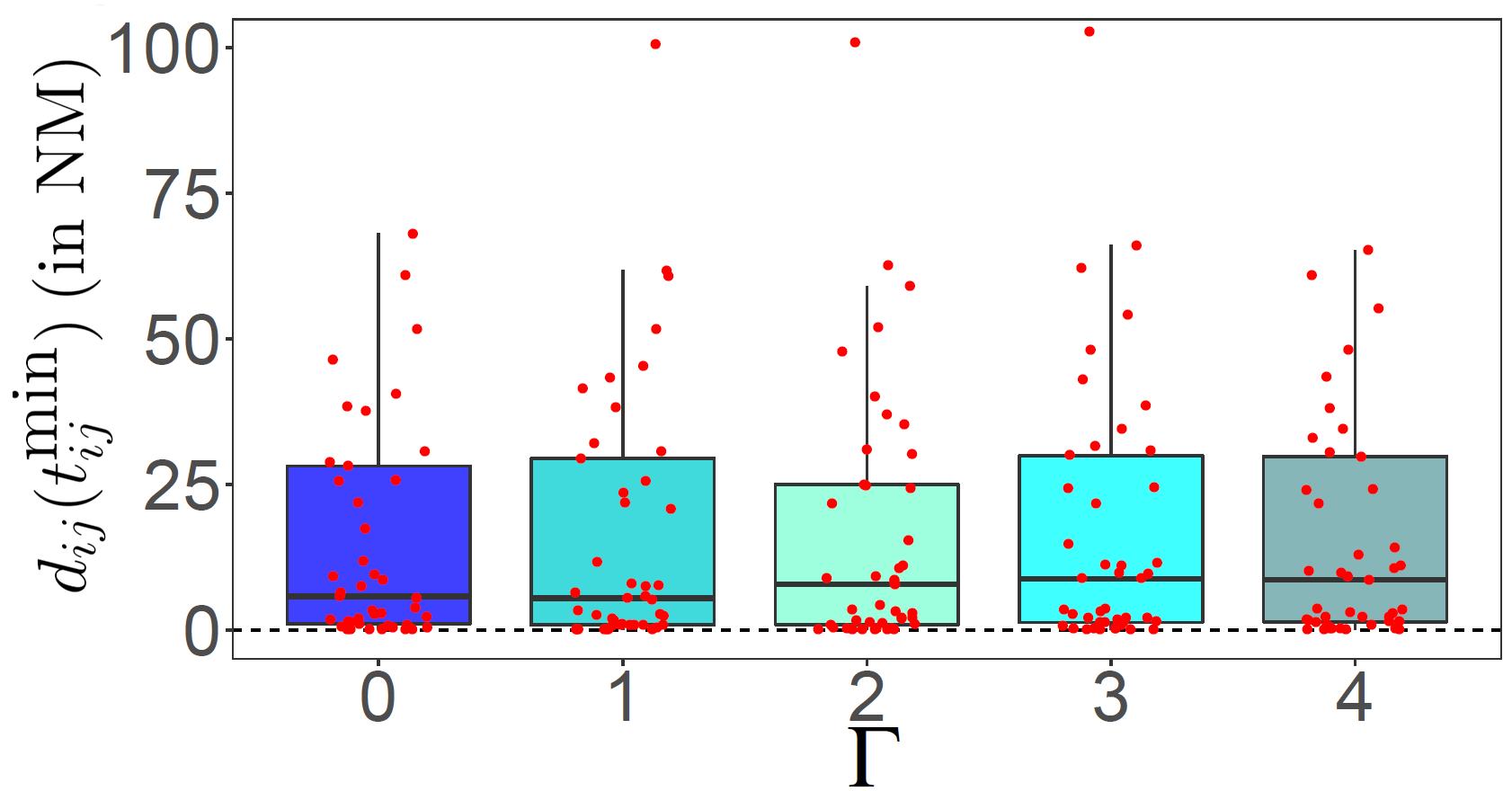}\label{rcp10_gamma}}\\
	\subfloat[RCP-20-1]{\includegraphics[width=0.4\textwidth]{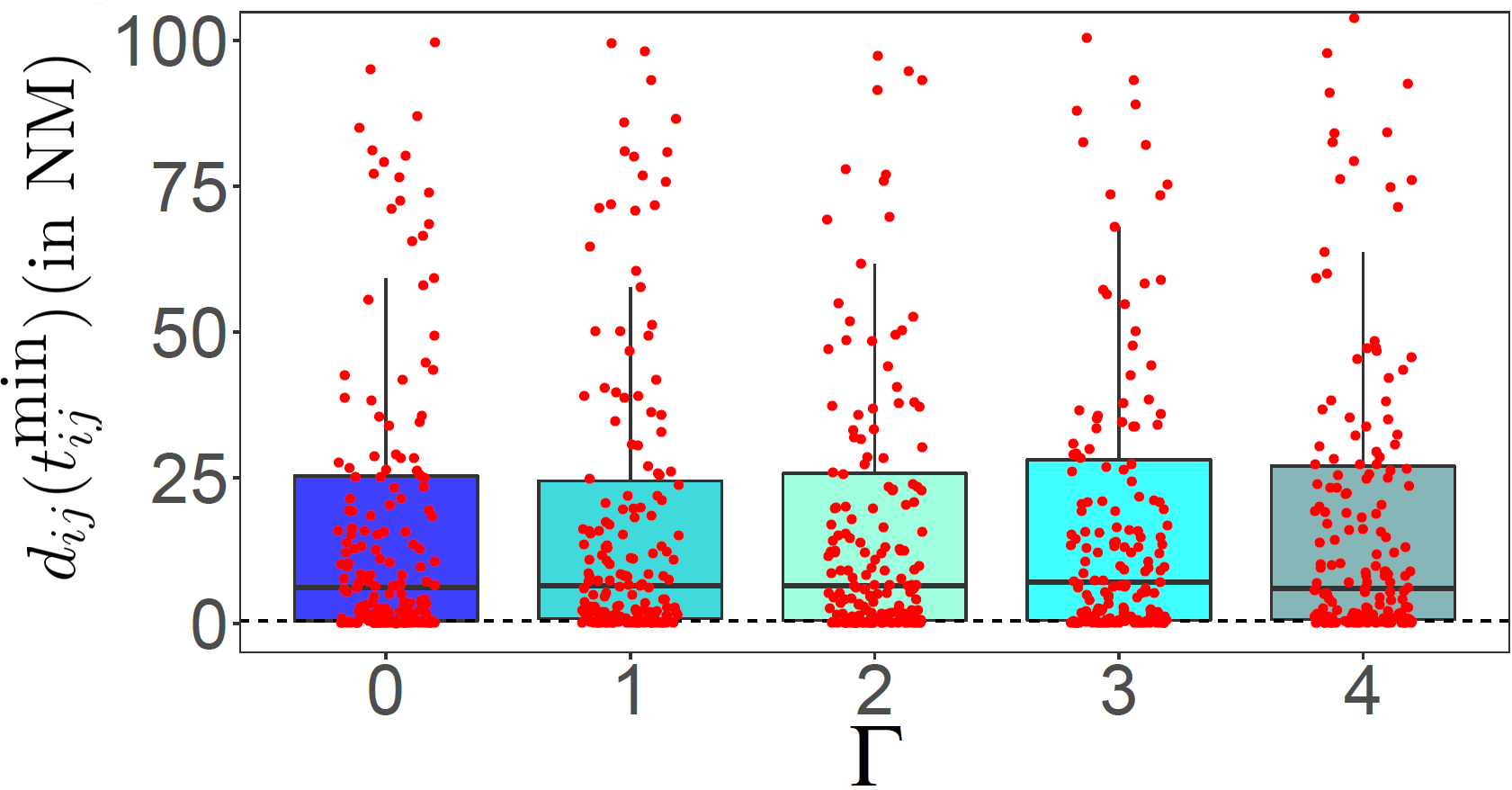}\label{rcp20_gamma}} \qquad
	\subfloat[RCP-30-1]{\includegraphics[width=0.4\textwidth]{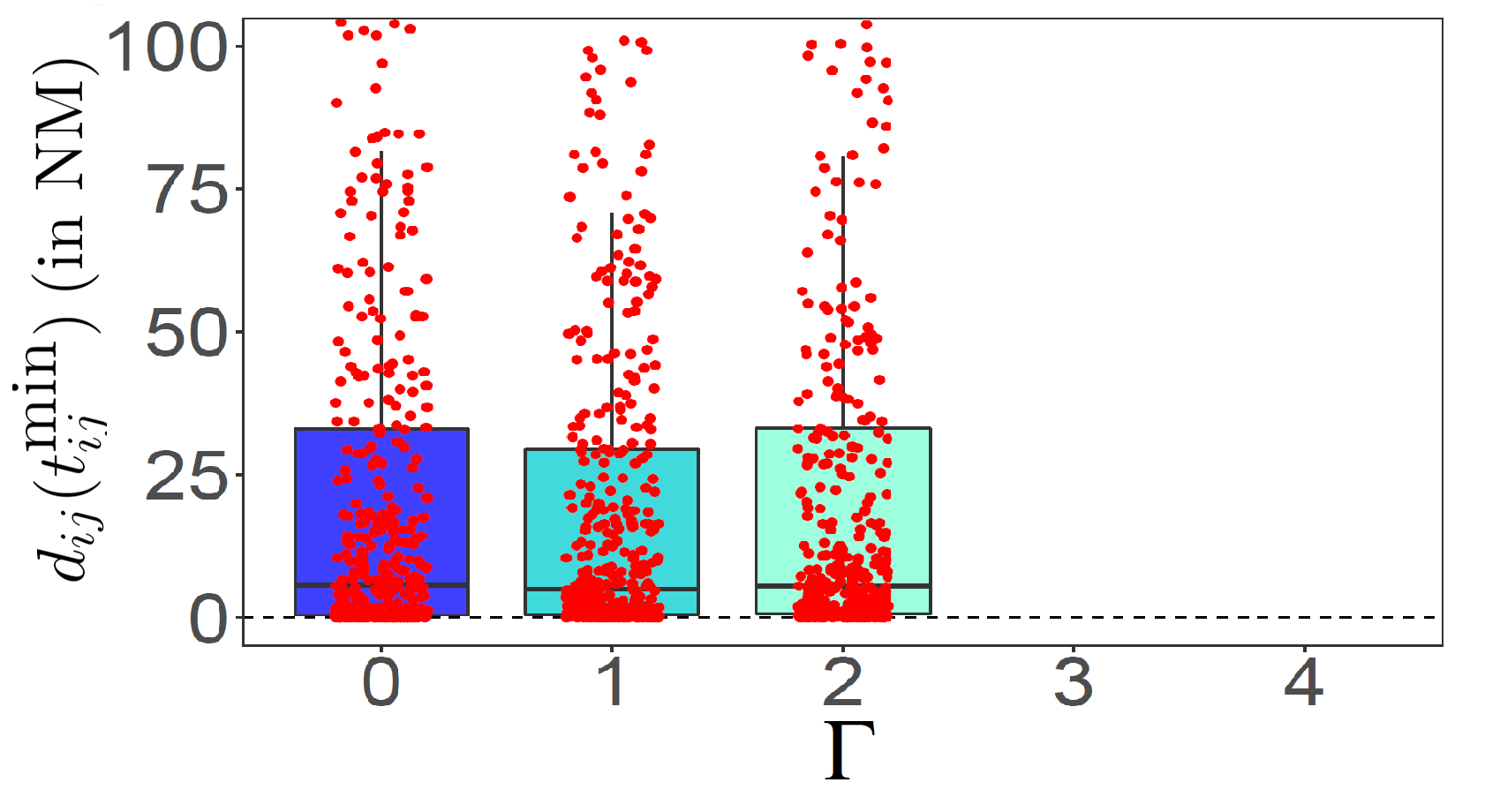}\label{rcp30_gamma}}\\
	\caption{Distribution of aircraft pairwise minimal separation distances after optimization for varying level of robustness $\Gamma$. The maximum uncertainty is set to $\bar{\epsilon} = 5\%$. The \textcolor{red}{red dots} correspond to individual values for each pair $(i,j) \in \P$. The \textcolor{black}{black} dashed line correspond to the minimal separation requirement of 5 NM.}
	\label{gammailus} 
\end{figure}

\begin{figure}[H]
	\centering\vspace{-2pt}
	\subfloat[CP-10]{\includegraphics[width=0.4\textwidth]{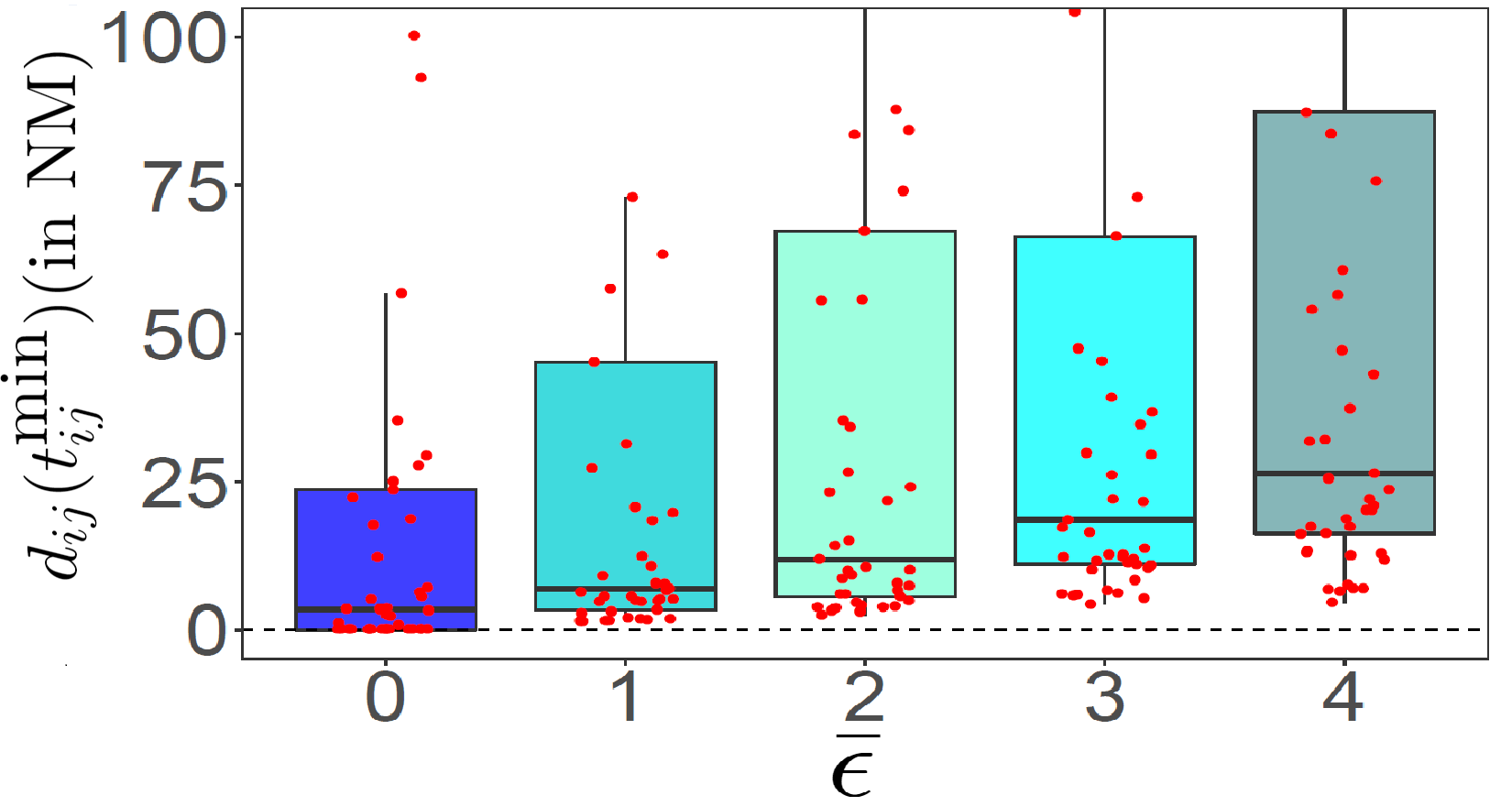}\label{cp_epsilon}} \qquad
	\subfloat[RCP-10-1]{\includegraphics[width=0.4\textwidth]{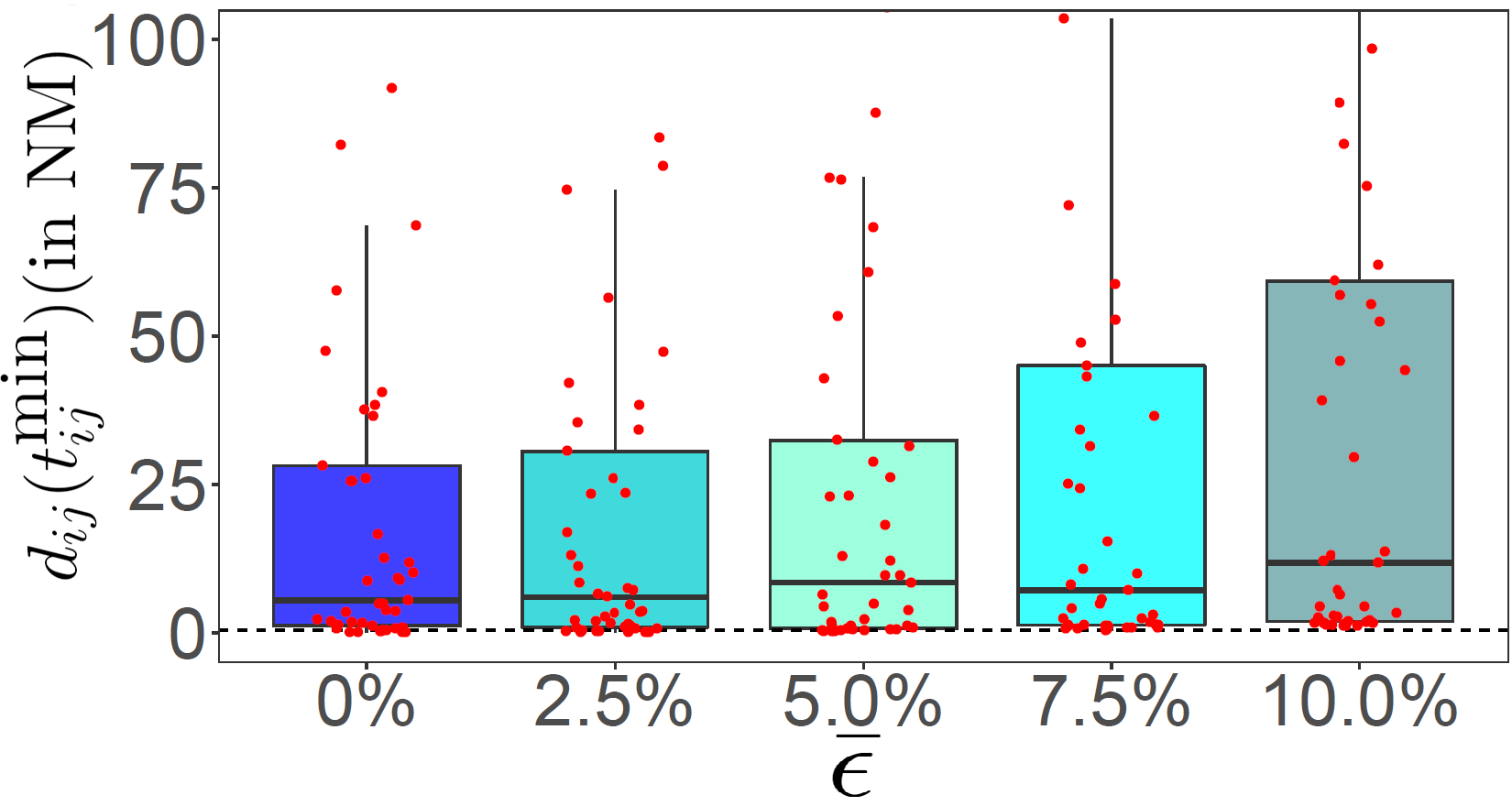}\label{rcp_10_epsilon}}\\
	\subfloat[RCP-20-1]{\includegraphics[width=0.4\textwidth]{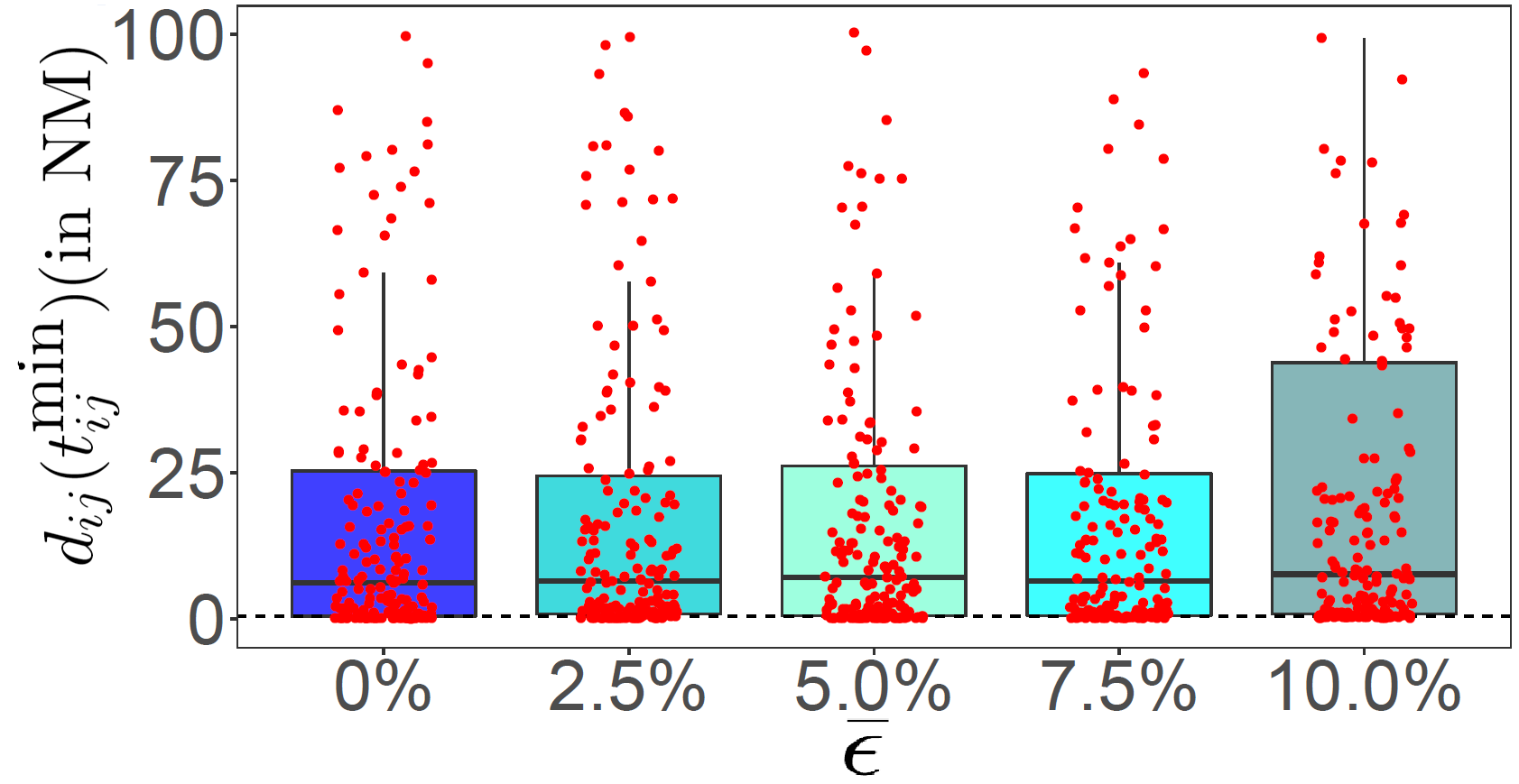}\label{rcp_20_epsilon}} \qquad
	\subfloat[RCP-30-1]{\includegraphics[width=0.4\textwidth]{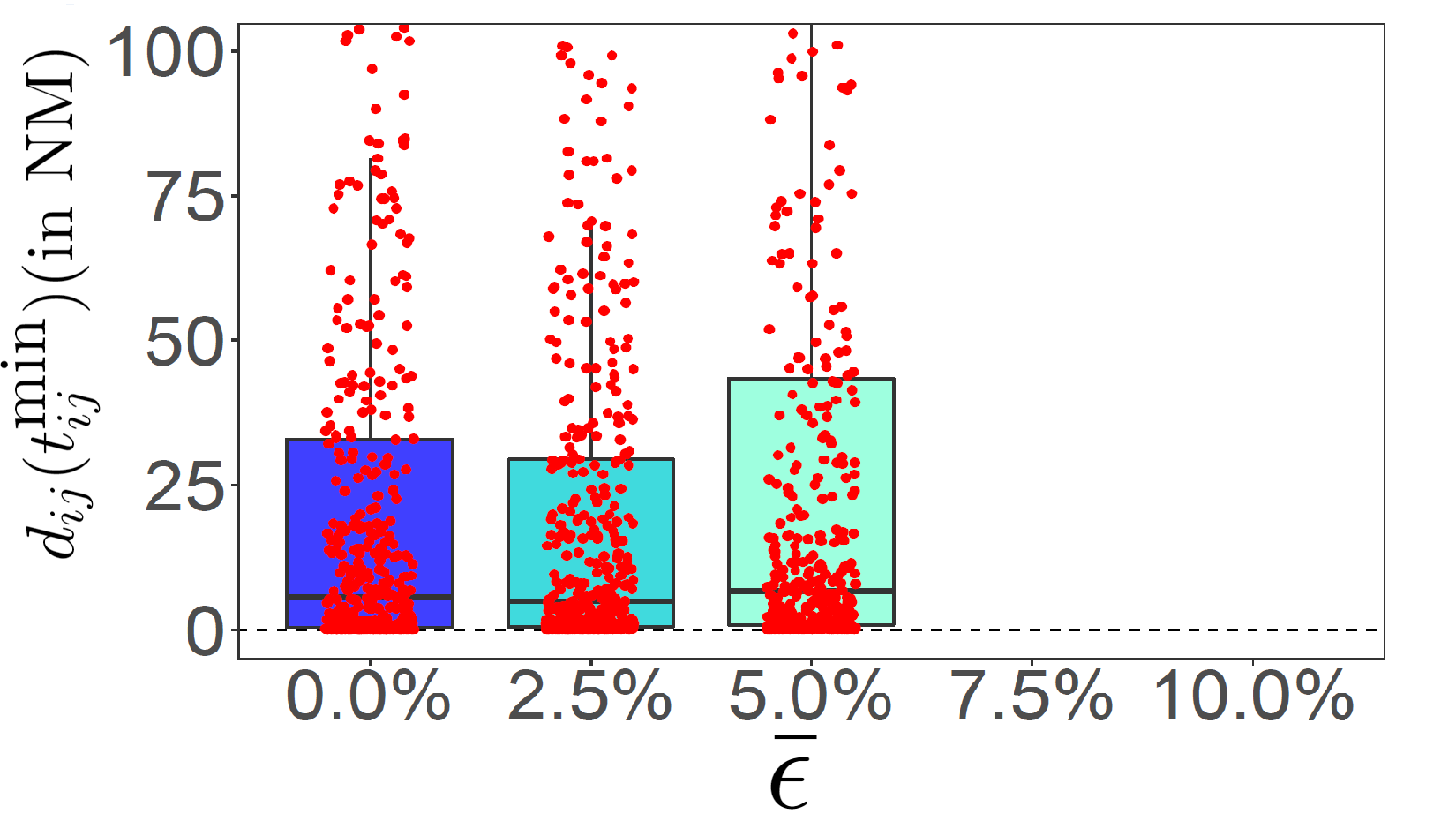}\label{rcp_30_epsilon}}\\
	\caption{Distribution of aircraft pairwise minimal separation distances after optimization for varying maximum uncertainty $\bar{\epsilon}$. The level of robustness is set to $\Gamma=4$. The \textcolor{red}{red dots} correspond to individual values for each pair $(i,j) \in \P$. The \textcolor{black}{black} dashed line correspond to the minimal separation requirement of 5 NM.}
	\label{epsilonilus} 
\end{figure}

Figure \ref{ilus} shows that there are significant differences between deterministic and robust aircraft trajectories, especially for instances with a high number of conflicts, such as CP-10 (Figure \ref{cp10}) and RCP-30-1 (Figure \ref{rcp31}). For RCP instances, this difference increases with the number of aircraft. To further analyze the behavior of the robust ACRP, we examine the distribution of post-optimization aircraft minimal separation distances, which is the distance $d_{ij}(\tm)$ between each pair of aircraft $(i,j) \in \P$ after optimization, against the level of robustness ($\Gamma$) and size of the uncertinty set ($\bar{\epsilon}$), in Figures \ref{gammailus} and \ref{epsilonilus}, respectively. For CP-10 (see Figure \ref{cp_gamma}), we observe that the average post-optimization minimal separation distance increases with $\Gamma$. Because this instance is symmetrical and all pairs are in conflict, the minimal separation to avoid conflict tends to grow rapidly with $\Gamma$. For RCP-10-1, RCP-20-1 and RCP-30-1 (see Figures \ref{rcp10_gamma}, \ref{rcp20_gamma}, \ref{rcp30_gamma}, which have a higher density of pairs of aircraft and less conflicts, we find that increasing the level of robustness does not significantly affect the distribution of aircraft minimal separation which remains near the minimal 5 NM mark. However, RCP-30-1 is found to be infeasible for $\Gamma = 3$ and $\Gamma = 4$, which suggests that even though trajectories may not substantially affected overall, such highly robust configurations fail to admit feasible solutions. Examing the same instances for a varying maximum uncertainty reveals an overall similar pattern although with more significant changes. Figure \ref{epsilonilus} shows that the average minimum separation distance increases in all instances, particularly for CP-10. This puts in evidence that increasing the value of $\bar{\epsilon}$ has a higher impact on the behaviour of the solution obtained. For RCP-10-1, RCP-20-1 and RCP-30-1 (see Figures \ref{rcp_10_epsilon}, \ref{rcp_20_epsilon} and \ref{rcp_30_epsilon}), we find that the distribution of aircraft pairwise minimal distances tends to increase for $\bar{\epsilon} \geq 5\%$. Further, we find that RCP-30-1 is infeasible for $\bar{\epsilon} = 7.5\%$ and $10\%$.

\begin{figure}[H]
	\centering\vspace{-2pt}
	\subfloat[RCP-10-1]{\includegraphics[width=0.33\textwidth]{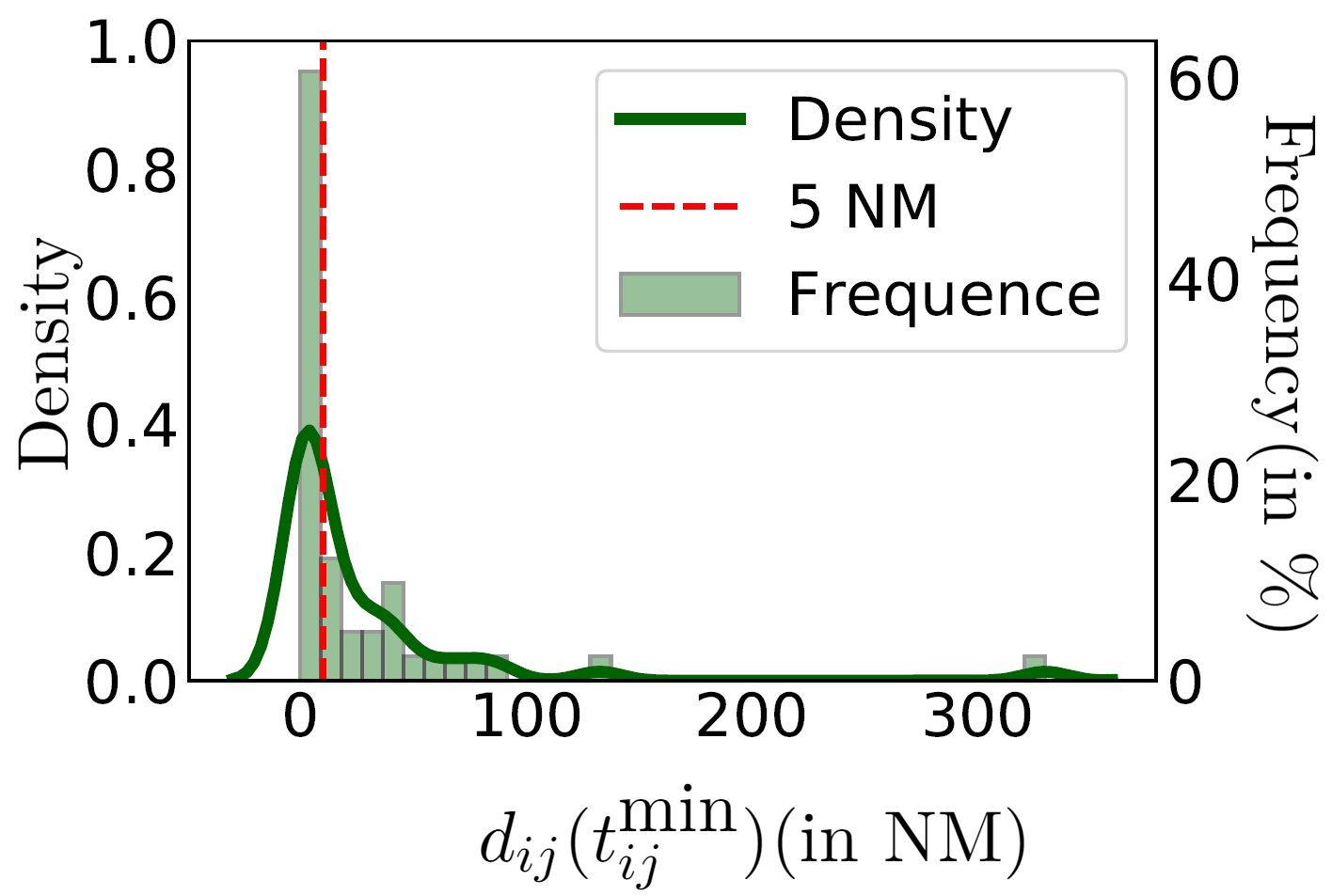}\label{frcp10}}
	\hfill
	\subfloat[RCP-20-1]{\includegraphics[width=0.33\textwidth]{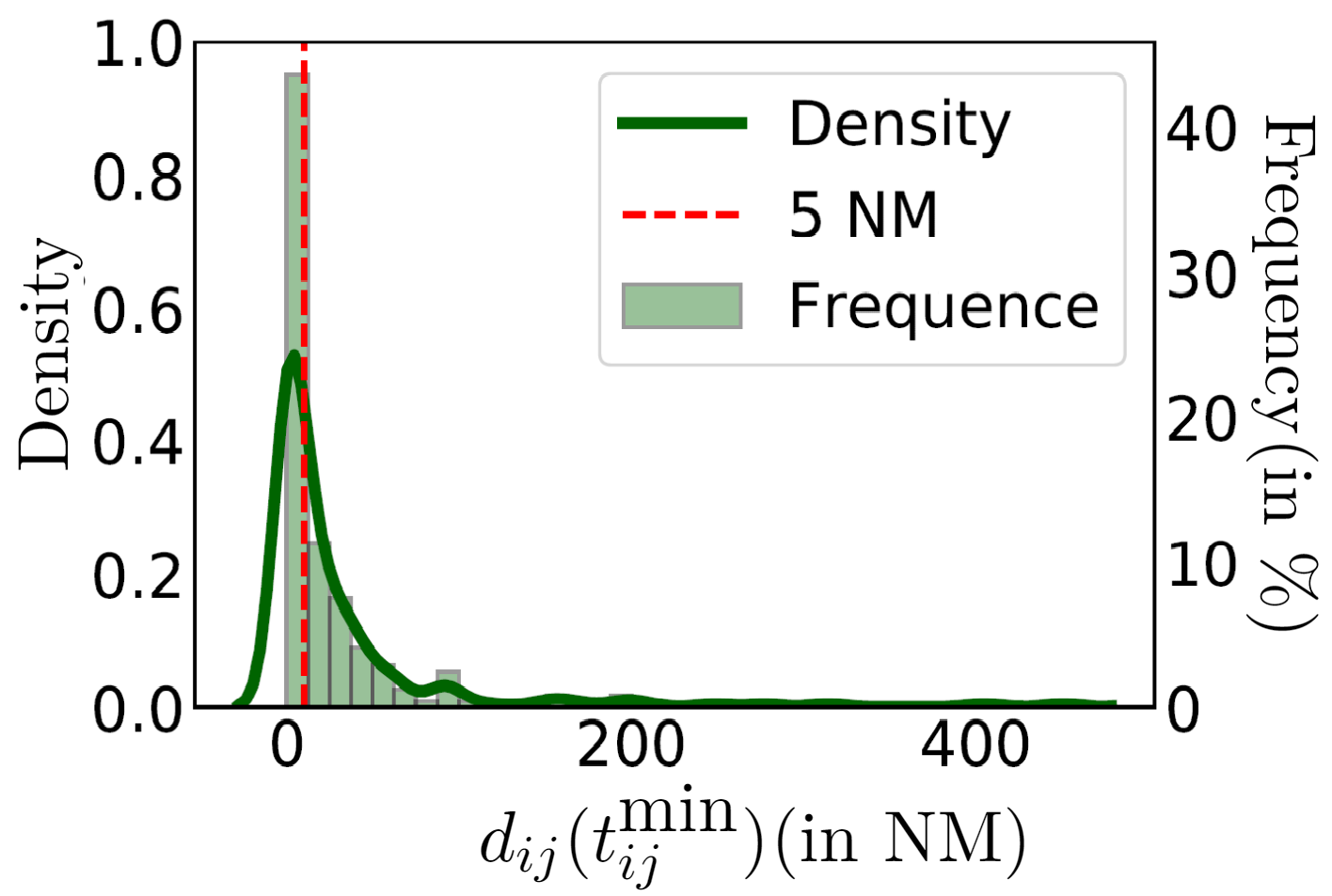}\label{frcp20}}
	\hfill
	\subfloat[RCP-30-1]{\includegraphics[width=0.33\textwidth]{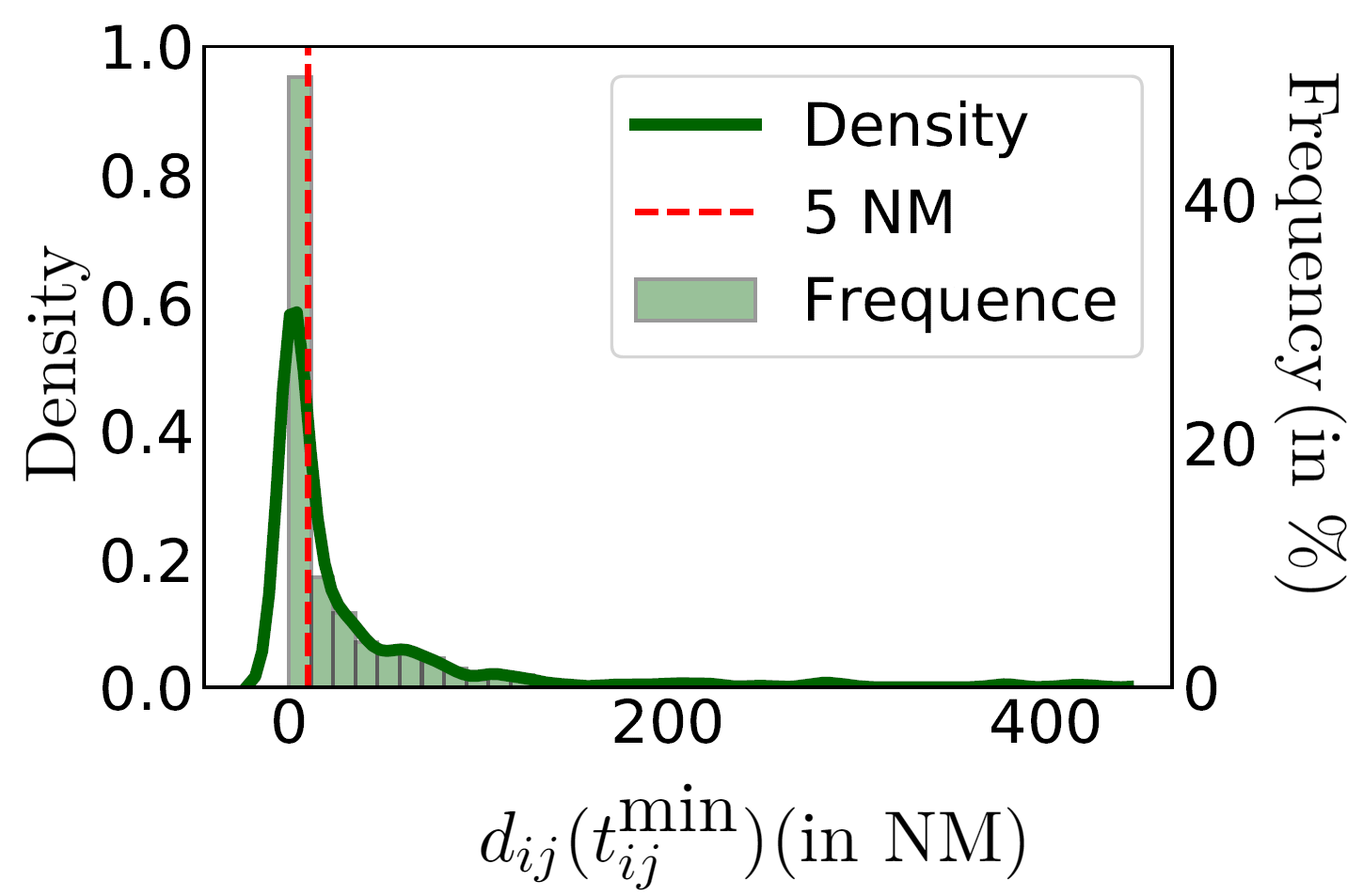}\label{frcp30}}
	\caption{Histogram and density plot for instances RCP-10-1, RCP-20-1 and RCP-30-1. The red dashed line represents the minimum separation distance of 5 NM.}
	\label{ilus2} 
\end{figure}

For RCP instances, we further examine the minimum separation distance between aircraft before optimization in Figure \ref{ilus2}. Figures \ref{frcp10}, \ref{frcp20} and \ref{frcp30} correspond to instances RCP-10-1, RCP-20-1 and RCP-30-1, respectively. To be more specific, around 55\%, 50\% and 45\% of the pair of aircraft are less than 5 NM of minimum separation (represented by the red dashed line on the plots) for RCP-10, RCP-20 and RCP-30, respectively. Most of the remaining pairs are within close distance to the minimum separation distance. This suggests that  introducing robustness can make such instances harder to solved. For RCP-10-1, RCP-20-1 and RCP-30-1, we can extrapolate that the distribution of aircraft initial minimal separation distance follows a chi-square distribution with small mean value. This confirms that in most instances, the mean value is close to 5 NM and a substantial number of aircraft pairs are in conflict or are only separated by a marginal amount in addition to the minimum separation distance.

We next conduct a comprehensive computational benchmarking of the proposed solution method for the robust ACRP.

\subsection{Sensitivity analysis on the level of robustness}
\label{fixedsets}

In this section, we analyze the performance of the robust complex number formulation for a varying level of robustness ($\Gamma$) for a fixed maximum uncertinaty of $\bar{\epsilon} = 5\%$. We present the results for 7 CP instances ranging from 4 to 10 aircraft in Table \ref{CP_g}. Results for RCP instances are reported in Table \ref{RCP_g} for 3 instance sizes with 10, 20 and 30 aircraft per group. For each instance group, 100 RCP instances are randomly generated and we report the average performance along with the standard deviation in parenthesis. For all instances, we compare the performance of the proposed formulations for varying level of robustness where $\Gamma=0$ corresponds to the deterministic case and $\Gamma=4$ corresponds to the most robust configuration. 

Each row in the results tables represents an instance (CP) or a group of instances (RCP). The header of the results tables is presented from left to right: $\Gamma$ is the level of robustness; $UB$ is the objective function value; Gap is the optimality gap in percent; Time is the total runtime in seconds, $n_i$ is the number of MIQCP iterations where $n_i$ = 0 means that an optimal solution was found or the time limit was reached during the initial MIQP solve, $n_t$ represents the proportion of instances that could not be solved within the time limit (10 minutes), i.e. the number of time-outs and $n_\emptyset$ indicates the number/proportion of infeasible instances.

\begin{table}[H]
\centering
\resizebox{0.65\columnwidth}{!}{%
\begin{tabular}{llllllllll}
\toprule
$\Gamma$ & Instance & UB & Gap (\%) & Time (s) & $n_i$ & $n_t$ & $n_\emptyset$ \\
\midrule
\multirow{2}{*}[-3em]{$0$} 
& CP-4	&	$6.25e^{-4}$	&	0.00	&	0.02	&	0	&	0 & 0	\\
& CP-5	&	$1.14e^{-3}$	&	0.00	&	0.05	&	0	&	0 & 0	\\
& CP-6	&	$1.81e^{-3}$	&	0.00	&	0.05	&	0	&	0 & 0	\\
& CP-7	&	$2.37e^{-3}$	&	0.04	&	0.17	&	0	&	0 & 0	\\
& CP-8	&	$3.46e^{-3}$	&	0.03	&	1.05	&	0	&	0 & 0	\\
& CP-9	&	$4.31e^{-3}$	&	0.02	&	32.1	&	0	&	0 & 0	\\
& CP-10	&	$5.55e^{-3}$	&	0.02	&	354	&	0	&	0 & 0	\\
\midrule
\multirow{2}{*}[-3em]{$1$} 
&	CP-4	&	$1.07e^{-3}$ &	0.03	&	0.16	&	0		&	0		&	0	\\
&	CP-5	&	$1.33e^{-3}$	&	0.00	&	0.05	&	0		&	0		&	0	\\
&	CP-6	&	$2.63e^{-3}$	&	0.00	&	0.13	&	0		&	0		&	0	\\
&	CP-7	&	$3.37e^{-3}$	&	0.00	&	0.45	&	0		&	0		&	0	\\
&	CP-8	&	$5.44e^{-3}$	&	0.01	&	3.73	&	1		&	0		&	0	\\
&	CP-9	&	$6.93e^{-3}$	&	0.01	&	68.9	&	4		&	0		&	0	\\
&	CP-10	&	$8.81e^{-3}$	&	0.01	&	51.9	&	1		&	0		&	0	\\
\midrule
\multirow{2}{*}[-3em]{$2$} 
&	CP-4	&	$1.21e^{-3}$	&	0.00	&	0.03	&	0		&	0		&	0	\\
&	CP-5	&	$1.49e^{-3}$	&	0.00	&	0.14	&	0		&	0		&	0	\\
&	CP-6	&	$3.13e^{-3}$	&	0.01	&	0.75	&	2		&	0		&	0	\\
&	CP-7	&	$4.18e^{-3}$	&	0.00	&	0.38	&	1		&	0		&	0	\\
&	CP-8	&	$6.93e^{-3}$	&	0.01	&	7.72	&	4		&	0		&	0	\\
&	CP-9	&	$8.92e^{-3}$	&	0.00	&	54.6	&	6		&	0		&	0	\\
&	CP-10 &	$1.39e^{-2}$	&	0.03	&	103	    &	6		&	0		&	0	\\
\midrule
\multirow{2}{*}[-3em]{$3$} 
&	CP-4	&	$1.22e^{-3}$	&	0.00	&	0.02	&	0		&	0		&	0	\\
&	CP-5	&	$1.59e^{-3}$	&	0.00	&	0.05	&	0		&	0		&	0	\\
&	CP-6	&	$3.38e^{-3}$	&	0.07	&	0.92	&	2		&	0		&	0	\\
&	CP-7	&	$4.64e^{-3}$	&	0.01	&	0.67	&	1		&	0		&	0	\\
&	CP-8	&	$8.98e^{-3}$	&	0.03	&	32.2	&	10		&	0		&	0	\\
&	CP-9	&	$1.09e^{-2}$	&	0.01	&	54.1	&	6		&	0		&	0	\\
&	CP-10&	$2.08e^{-2}$	&	0.00	&	174	    &	8   	&	0		&	0	\\
\midrule
\multirow{2}{*}[-3em]{$4$} 
&	CP-4	&	$1.22e^{-3}$	&	0.00	&	0.03	&	0		&	0		&	0	\\
&	CP-5	&	$1.63e^{-3}$	&	0.00	&	0.05	&	0		&	0		&	0	\\
&	CP-6	&	$3.49e^{-3}$	&	0.04	&	0.80	&	0		&	0		&	0	\\
&	CP-7	&	$4.77e^{-3}$	&	0.01	&	3.87	&	3		&	0		&	0	\\
&	CP-8	&	$9.20e^{-3}$	&	0.00	&	26.7	&	9		&	0		&	0	\\
&	CP-9	&	$1.24e^{-2}$	&	0.00	&	142	&	10		&	0		&	0	\\
&	CP-10&	$2.38e^{-2}$	&	0.02	&	215	&	9		&	0		&	0	\\
\bottomrule
\end{tabular}
}
\caption{Summary of results for CP instances with a speed control range of $[-6\%,+3\%]$ and a heading control range of $[-30^\circ,+30^\circ]$ using a maximal uncertainty of $\bar{\epsilon} = 5\%$.}
\label{CP_g}
\end{table}

The experiments on the CP instances (Table \ref{CP_g}) show that, as expected, the best objective value (UB) increases with the number of aircraft. Further, we find that UB increases with the level of robustness. Specifically, for CP-4, the objective value  increases over 90\% from $\Gamma=0$ to $\Gamma=4$; while for CP-10 this figure is greater than 300\%. On average, we find that compared to the deterministic case, introducing a maximum level of robustness in CP instances yields an increase of 143.62\% in the objective value. While the runtime increases exponentially with the number of aircraft, increasing $\Gamma$ does not significantly impact runtime, even when compared to the deterministic case of $\Gamma=0$. We find that increasing $\Gamma$ tends to increase the number of iterations of the solution algorithm, which suggests that the higher levels of robustness require more deconfliction resources, i.e. speed or heading change.

%The implementation of the proposed approach on CP instances reveals that while all instances up to CP-8 can be solved via the MIQP iteration when the value for $\Gamma = 0$ or $\Gamma =1$, for any higher values may require additional MIQCP iterations. We find that for increasing $\Gamma$, the time limit is enough to solve all instances, but with a higher level of randomness, more speed deviation is observed and therefore, more MIQCP iterations are necessary to obtain a conflict-free and violation-free optimal solution. This shows an expected behaviour such as that with higher randomness effect granted by increasing $\Gamma$ leads to a more complex problem and therefore a more difficult to solve. In addition to this, with a higher value of $\Gamma$ which represents a higher level of robustness, it is expected that we would observe more violations in speed to compensate for a wider range of speed and angle manoeuvres. With more violations, the solution will consequently require more iteration to be solved to optimality. By the time limit imposed, some instances could not be solved for a higher $\Gamma$. However, as observed, the gamma formulation is enough to provide robust solutions given more runtime.

\begin{table}[H]
\centering
\resizebox{0.8\columnwidth}{!}{%
\begin{tabular}{llllllllll}
\toprule
$\Gamma$ & Instance & UB & Gap & Time & $n_i$ & $n_t$ & $n_\emptyset$ \\
\midrule
\multirow{2}{*}[-0.65em]{$0$} 
& RCP-10 & $2.2e^{-4}$ ($2e^{-4}$) & 0.00 (0.0) & 0.05 (0.01) & 0.0 (0.0) & 0 & 0\\
& RCP-20 & $1.7e^{-3}$ ($9e^{-4}$) & 0.00 (0.0) & 0.26 (0.10) & 0.0 (0.0) & 0 & 0\\
& RCP-30 & $7.1e^{-3}$ ($2e^{-3}$) & 0.17 (0.7) & 135 (239) &  0.4 (0.6) & 3 & 0\\
\midrule
\multirow{2}{*}[-0.65em]{$1$} 
& RCP-10 & $1.0e^{-2}$ ($7e^{-6}$) & 0.00 (0.0) & 0.24 (0.06) & 0.0 (0.0) & 0 & 0\\
& RCP-20 & $9.8e^{-2}$ ($4e^{-2}$) & 0.00 (0.0) & 10.2 (4.84) & 0.0 (0.0) & 0 & 0\\
& RCP-30 & $7.5e^{-1}$ ($2e^{-1}$) & 0.17 (0.7) & 287 (116) &  0.4 (0.0) & 10 & 60 \\
\midrule
\multirow{2}{*}[-0.65em]{$2$} 
& RCP-10 & $2.7e^{-2}$ ($2e^{-2}$) & 0.00 (0.0) & 0.34 (0.13) & 0.0 (0.0) & 0 & 0\\
& RCP-20 & $2.8e^{-2}$ ($9e^{-3}$) & 0.00 (0.0) & 25.1 (16.8) & 0.0 (0.0) & 8 & 34\\
& RCP-30 & - & - & - & - & - & 100\\
\midrule
\multirow{2}{*}[-0.65em]{$3$} 
& RCP-10 & $4.2e^{-2}$ ($2e^{-2}$) & 0.00 (0.0) & 0.54 (0.18) & 0.0 (0.0) & 0 & 0\\
& RCP-20 & $7.4e^{-1}$ ($2e^{-1}$) & 0.00 (0.0) & 134 (55.6) & 0.0 (0.0) & 10 & 46\\
& RCP-30 & - & - & - & - & - & 100\\
\midrule
\multirow{2}{*}[-0.65em]{$4$} 
& RCP-10 & $5.3e^{-2}$ ($3e^{-2}$) & 0.00 (0.0) & 0.59 (0.19) & 0.0 (0.0) & 0 & 0\\
& RCP-20 & $8.8e^{-1}$ ($2e^{-1}$) & 0.00 (0.0) & 201 (144) & 0.0 (0.0) & 15 & 67\\
& RCP-30 & - & - & - & - & - & 100 \\
\bottomrule
\end{tabular}
}
\caption{Summary of results for RCP instances with a speed control range of $[-6\%,+3\%]$ and a heading control range of $[-30^\circ,+30^\circ]$ using a maximal uncertainty of $\bar{\epsilon} = 5\%$.}
\label{RCP_g}
\end{table}
 
The results of the experiments on RCP instances are summarized in Table \ref{RCP_g}. Increasing the level of robustness in RCP instances also tends to increase the objective value (UB), and the effect is on average magnified on instances with a large number of aircraft. For RCP-10 instances, the average objective value for $\Gamma=4$ is two orders of magnitude greater than that obtained using a deterministic configuration. For RCP-20 instances, the average objective values obtained are one order of magnitude greater than those obtained for RCP-10. We observe that RCP-10 instances can be solved within less than a second using for any level of robustness. RCP-20 instances require less than a minute for $\Gamma \leq 2$, on average, but using $\Gamma$ equal to 3 and to  4, requires around 2 and 3 minutes, respectively. We also find that 34\%, 46\% and 67\% of RCP-20 instances cannot be solved when $\Gamma$ is equal to 2, 3 and 4, respectively. The results for RCP-30 instances reveal that all problems with $\Gamma \geq 2$ are infeasible, while only 40\% of these instances can be solved with $\Gamma  = 1$ and all 100 RCP-30 instances are feasible in the deterministic case.

\subsection{Sensitivity analysis on the size of the uncertainty set}
\label{fixedgamma}

\begin{table}
\centering
\resizebox{0.65\columnwidth}{!}{%
\begin{tabular}{llllllllll}
\toprule
$\bar{\epsilon}$ (\%) & Instance & UB & Gap (\%) & Time (s) & $n_i$ & $n_t$ & $n_\emptyset$ \\
\midrule
\multirow{2}{*}[-3em]{$0.0$} 
& CP-4	&	$6.25e^{-4}$	&		0.00	&	0.20	&	0	&	0 &	0		\\
& CP-5	&	$1.14e^{-3}$	&		0.00	&	0.40	&	0	&	0 &	0		\\
& CP-6	&	$1.81e^{-3}$	&		0.00	&	0.64	&	0	&	0 &	0		\\
& CP-7	&   $2.37e^{-3}$	&		0.00	&	0.39	&	0	&	0 &	0		\\
& CP-8	&	$3.46e^{-3}$	&		0.02	&	3.04	&	0	&	0 &	0		\\
& CP-9	&	$4.31e^{-3}$	&		0.02	&	7.97	&	0	&	0 &	0		\\
& CP-10	&	$5.55e^{-3}$	&		0.02	&	72.3	&	0	&	0 &	0		\\
\midrule
\multirow{2}{*}[-3em]{$2.5$} 
&	CP-4	&	$1.22e^{-3}$	&	0.00	&	0.03	&	0		&	0		&	0	\\
&	CP-5	&	$1.63e^{-3}$	&	0.00	&	0.05	&	0		&	0		&	0	\\
&	CP-6	&	$3.49e^{-3}$	&	0.04	&	0.80	&	2		&	0		&	0	\\
&	CP-7	&	$4.77e^{-3}$	&	0.01	&	3.87	&	4		&	0		&	0	\\
&	CP-8	&	$9.20e^{-3}$	&	0.00	&	26.7	&	9		&	0		&	0	\\
&	CP-9	&	$1.24e^{-2}$	&	0.00	&	142	    &	10		&	0		&	0	\\
&	CP-10&	$2.38e^{-2}$	&	0.02	&	215 	&	8		&	0		&	0	\\
\midrule
\multirow{2}{*}[-3em]{$5.0$} 
&	CP-4	&	$1.70e^{-3}$	&	0.00	&	0.03	&	0		&	0		&	0	\\
&	CP-5	&	$3.94e^{-3}$	&	0.02	&	0.65	&	2		&	0		&	0	\\
&	CP-6	&	$1.00e^{-2}$	&	0.00	&	4.12	&	6		&	0		&	0	\\
&	CP-7	&	$2.12e^{-2}$	&	0.00	&	16.2	&	9		&	0		&	0	\\
&	CP-8	&	$4.01e^{-2}$	&	0.53	&	12.9	&	5		&	0		&	0	\\
&	CP-9	&	$6.35e^{-2}$	&	2.13	&	31.9	&	5		&	0		&	0	\\
&	CP-10&	$8.45e^{-2}$	&	3.23	&	20.4	&	2		&	0		&	0	\\
\midrule
\multirow{2}{*}[-3em]{$7.5$} 
&	CP-4	&	$3.80e^{-3}$	&	0.00	&	1.28	&	4		&	0		&	0	\\
&	CP-5	&	$1.14e^{-2}$	&	0.00	&	4.56	&	11		&	0		&	0	\\
&	CP-6	&	$3.02e^{-2}$	&	0.00	&	3.14	&	5		&	0		&	0	\\
&	CP-7	&	$6.12e^{-2}$	&	0.46	&	9.44	&	8		&	0		&	0	\\
&	CP-8	&	$9.53e^{-2}$	&	7.55	&	5.88	&	3		&	0		&	0	\\
&	CP-9	&	$1.52e^{-1}$	&	1.99	&	9.36	&	3		&	0		&	0	\\
&	CP-10&	-	&-	& -	&	-		& -		&	100	\\
\midrule
\multirow{2}{*}[-3em]{$10.0$} 
&	CP-4	&	$7.45e^{-3}$	&	0.00	&	4.15	&	12		&	0		&	0	\\
&	CP-5	&	$2.14e^{-2}$	&	0.00	&	2.81	&	8		&	0		&	0	\\
&	CP-6	&	$6.24e^{-2}$	&	0.00	&	3.92	&	6		&	0		&	0	\\
&	CP-7	&	$1.09e^{-1}$	&	0.25	&	4.68	&	5		&	0		&	0	\\
&	CP-8	&	$1.98e^{-1}$	&	3.52	&	4.14	&	3		&	0		&	0	\\
&	CP-9	&	$3.48e^{-1}$	&	18.8	&	5.57	&	3		&	0		&	0	\\
&	CP-10   &	-	            &	-	&	-	&	-		&	-		& 100	\\
\bottomrule
\end{tabular}
}
\caption{Summary of results for CP instances with a speed control range of $[-6\%,+3\%]$ and a heading control range of $[-30^\circ,+30^\circ]$ using a level of robustness of $\Gamma = 4$.}
\label{CP_u}
\end{table}

\begin{table}
\centering
\resizebox{0.8\columnwidth}{!}{%
\begin{tabular}{llllllllll}
\toprule
$\bar{\epsilon}$ (\%) & Instance & UB & Gap (\%) & Time (s) & $n_i$ & $n_t$ & $n_\emptyset$ \\
\midrule
\multirow{2}{*}[-0.65em]{$0.0$} 
& RCP-10 & $2.2e^{-4}$ ($2e^{-4}$) & 0.00 (0.0) & 0.05 (0.01) & 0.0 (0.0) & 0 & 0\\
& RCP-20 & $1.7e^{-3}$ ($9e^{-4}$) & 0.01 (0.0) & 0.24 (0.1) & 0.0 (0.0) & 0 & 0\\
& RCP-30 & $7.2e^{-3}$ ($2e^{-3}$) & 0.01 (0.1) & 66.3 (135) &  1.4 (0.8) & 3 & 0\\
\midrule
\multirow{2}{*}[-0.65em]{$2.5$} 
& RCP-10 & $1.6e^{-2}$ ($1e^{-2}$) & 0.00 (0.0) & 0.47 (0.15) & 0.0 (0.0) & 0 & 1\\
& RCP-20 & $8.9e^{-1}$ ($2e^{-1}$) & 0.00 (0.0) & 193 (138) & 0.0 (0.0) & 12 & 4\\
& RCP-30 & - & - & - & - & - & 100\\
\midrule
\multirow{2}{*}[-0.65em]{$5.0$} 
& RCP-10 & $5.3e^{-2}$ ($3e^{-2}$) & 0.00 (0.0) & 0.60 (0.18) & 0.0 (0.0) & 0 & 1\\
& RCP-20 & $8.8e^{-1}$ ($2e^{-1}$) & 0.00 (0.0) & 200 (144) & 0.0 (0.0) & 2 & 93\\
& RCP-30 & - & - & - & - & - & 100\\
\midrule
\multirow{2}{*}[-0.65em]{$7.5$} 
& RCP-10 & $1.3e^{-1}$ ($6e^{-2}$) & 0.00 (0.0) & 0.76 (0.24) & 0.0 (0.0) & 0 & 9\\
& RCP-20 & - & - & - & - & - & 100\\
& RCP-30 & - & - & - & - & - & 100\\
\midrule
\multirow{2}{*}[-0.65em]{$10.0$} 
& RCP-10 & $3.0e^{-1}$ ($1e^{-1}$) & 0.00 (0.0) & 0.89 (0.20) & 0.0 (0.0) & 0 & 26\\
& RCP-20 & - & - & - & - & - & 100\\
& RCP-30 & - & - & - & - & - & 100\\
\bottomrule
\end{tabular}
}
\caption{Summary of results for RCP instances with a speed control range of $[-6\%,+3\%]$ and a heading control range of $[-30^\circ,+30^\circ]$ using a level of robustness of $\Gamma = 4$.}
\label{RCP_u}
\end{table}

For this experiment, we compare the performance of the proposed formulation for a maximum level of robustness, i.e. $\Gamma = 4$ under varying maximum uncertainty $\bar{\epsilon}$: 0\% (for the deterministic counterpart), 2.5\%, 5\%, 7.5\% and 10\%. The performance is reported similarly as in the previous section. We present the results for 7 CP instances ranging from 4 to 10 aircraft in Table \ref{CP_u} and results for RCP instances are reported in Table \ref{RCP_u}.  

The experiments on the CP instances (Table \ref{CP_u}) show that, as expected, the objective value (UB) increases with the maximum uncertainty. For CP-4, we observe a growth in 1000\% of the objective value when using a maximal uncertainty of $\bar{\epsilon} = 10\%$ compared to the deterministic case. For CP-10, this growth is of the order of 7000\%. We find that CP instances with up to 9 aircraft can be solved with up $\bar{\epsilon}=10\%$, however CP-10 is infeasible for $\bar{\epsilon}=7.5\%$ and $\bar{\epsilon}=10\%$. All deterministic CP problems can be solved without any iteration ($n_i = 0$) of the algorithm. In turn, introducing uncertainty in the robust ACRP triggers several iterations of the solution algorithm which suggests that the relaxed mixed-integer convex programs are generating more violations due to the required robustness. Using $\bar{\epsilon} = 2.5\%$, the instances can be solved within 0.03 s for CP-4 and up to 215 seconds for CP-10. Increasing $\bar{\epsilon}$ results in a reduced runtime for all instances compared to the deterministic case.

The results obtained using RCP instances (Table \ref{RCP_u}) reveal that while all problems are feasible under deterministic conditions, the proportion of infeasible problems increase with the maximum uncertainty and the number of aircraft. Notably, all RCP-30-1 instances are found to be infeasible for $\bar{\epsilon} \geq 2.5\%$ (recall that a maximal level of robustness is used). Examining the change in the objective value, for RCP-10 instances we observe that on average the total deviation of aircraft increases by three orders of magnitude. In terms of runtime, RCP-10 instances, are solved within 1 s using any value of $\bar{\epsilon}$. RCP-20 instances, require on average 0.24 s in the deterministic case. Using $\bar{\epsilon} = 5\%$, 93 instances are infeasible, 2 timed-out and among the remaining 5 instances solved to optimality the average runtime is 200 s.\\

The results of the sensitivity analyses conducted highlight that while the robust ACRP can be solved without significantly increasing the level of computational resources required, the impact of on system costs (total deviation) increase rapidly with the level of robustness and/or uncertainty, and the number of aircraft. Further, these analyses revealed that the likelihood of infeasibility increased rapidly along the same directions. We next attempt to explain this behavior of the robust ACRP.

%combined with the previous section solidifies that this formulation is capable of solving dense and complex scenarios, but with increasing uncertainty sets, we can observe that it becomes a very complex problem. As expected, to guarantee that a larger margin for error is available for each aircraft, an instance cannot have the same amount of aircraft available as it would have in its deterministic counterpart. The main difference between changing $\Gamma$ and $\bar{\epsilon}$ is that in the former we explore how much an uncertainty set will impact each variable and component in our model based on the level of robustness. In the latter, by varying the size of the uncertainty set, we are likely to explore situations that are not expected and even unrealistic for the aircraft given that large uncertainty sets are not commonly accepted are standard margins of error. 

\subsection{Feasibility analysis}
\label{fease}

\begin{figure}[ht]
	\centering\vspace{-2pt}
	\subfloat[$\Gamma = 1$ ]{\includegraphics[width=0.45\textwidth]{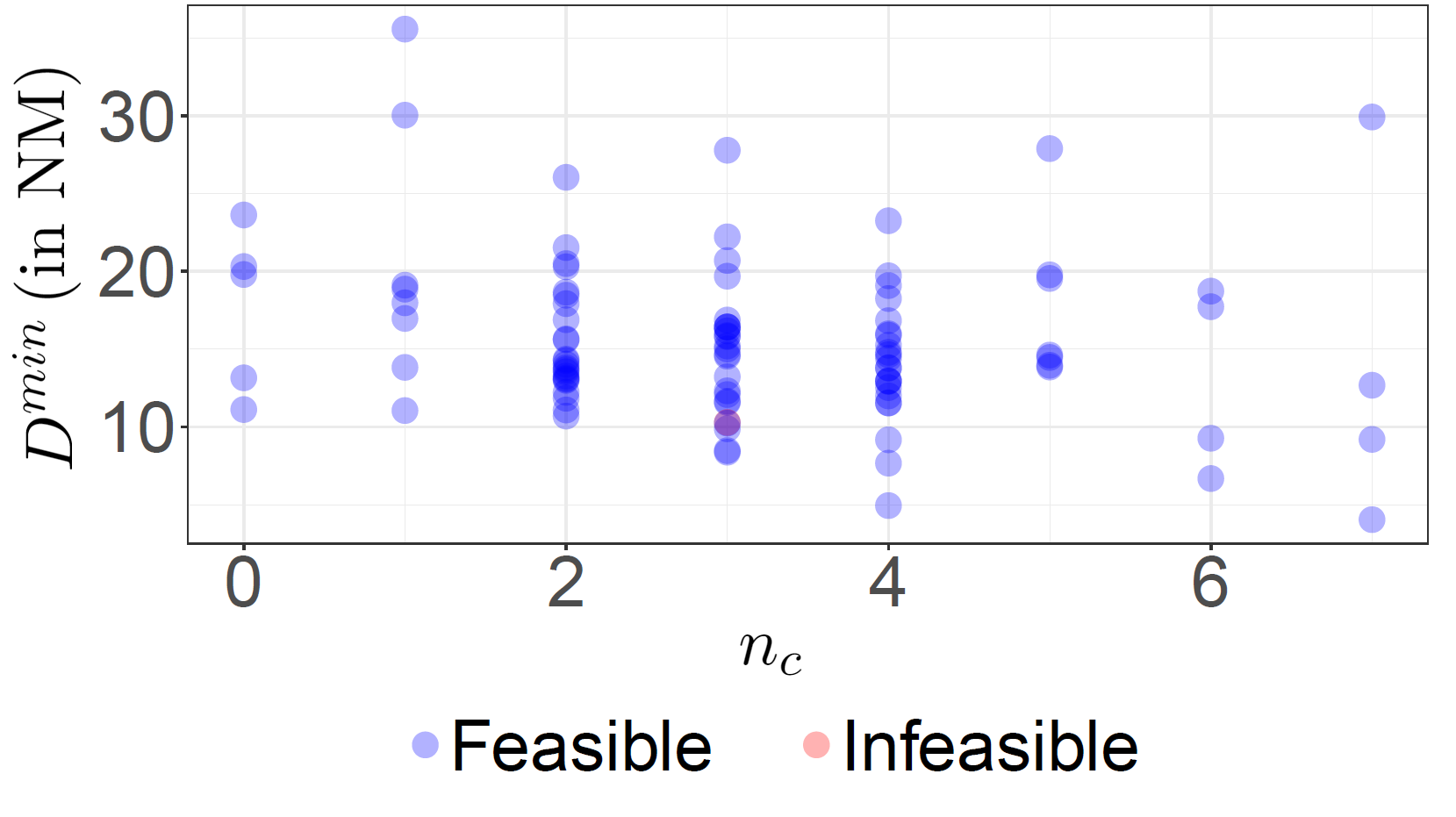}\label{10g0}} \hfill
	\subfloat[$\Gamma = 2$ ]{\includegraphics[width=0.45\textwidth]{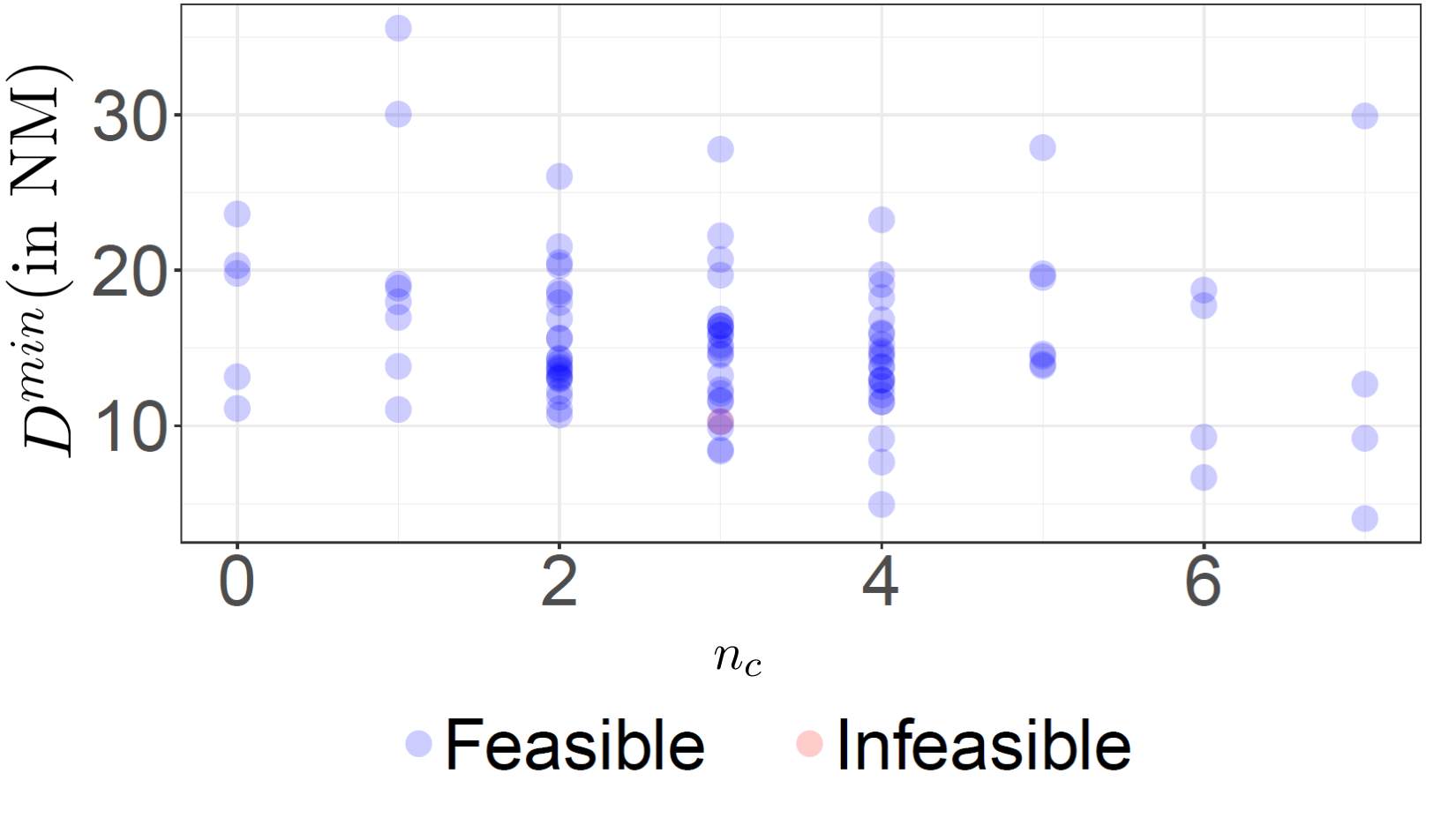}\label{10g2}}\\
	\subfloat[$\Gamma = 3$ ]{\includegraphics[width=0.45\textwidth]{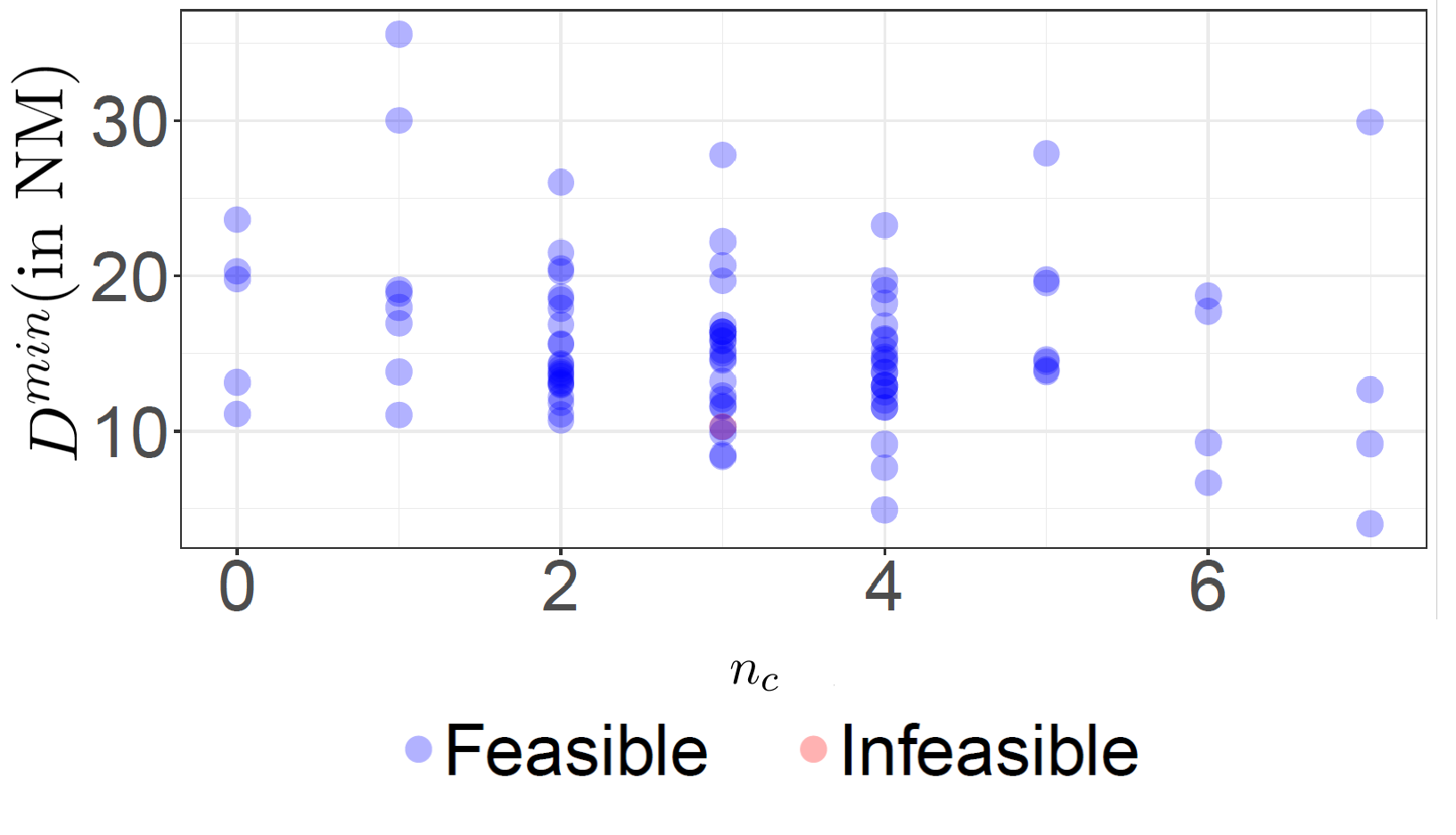}\label{10g3}} \hfill
	\subfloat[$\Gamma = 4$ ]{\includegraphics[width=0.45\textwidth]{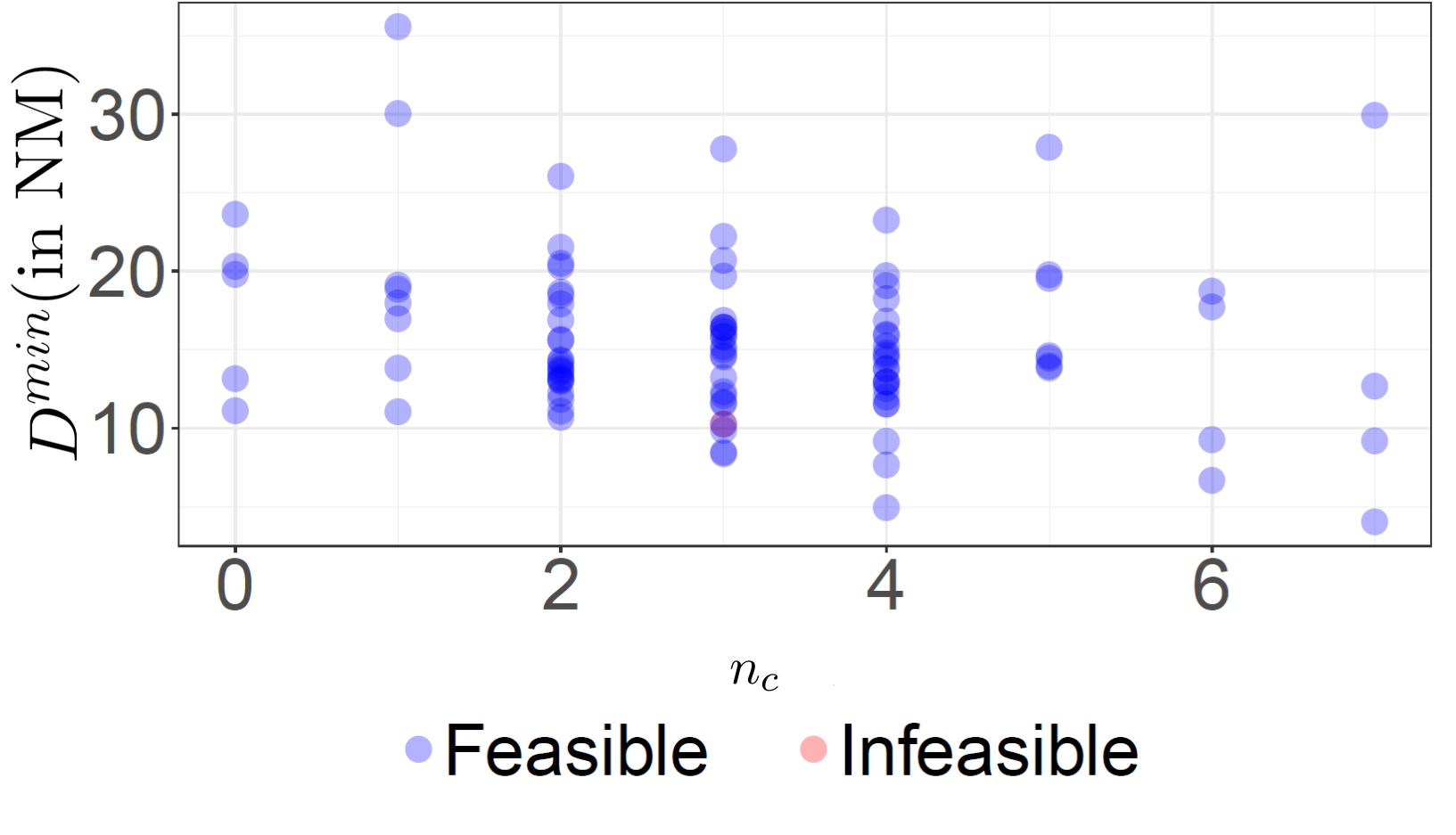}\label{10g4}}\\
	\caption{Feasibility of 100 RCP-10 instances based on the number of conflict $n_c$ and the total pairwise minimal distance between all pairs of aircraft $D^{\text{min}}$ using $\bar{\epsilon} = 5\%$.}
	\label{ilus3} 
\end{figure}

\begin{figure}[ht]
	\centering\vspace{-2pt}
	\subfloat[$\Gamma = 1$ ]{\includegraphics[width=0.45\textwidth]{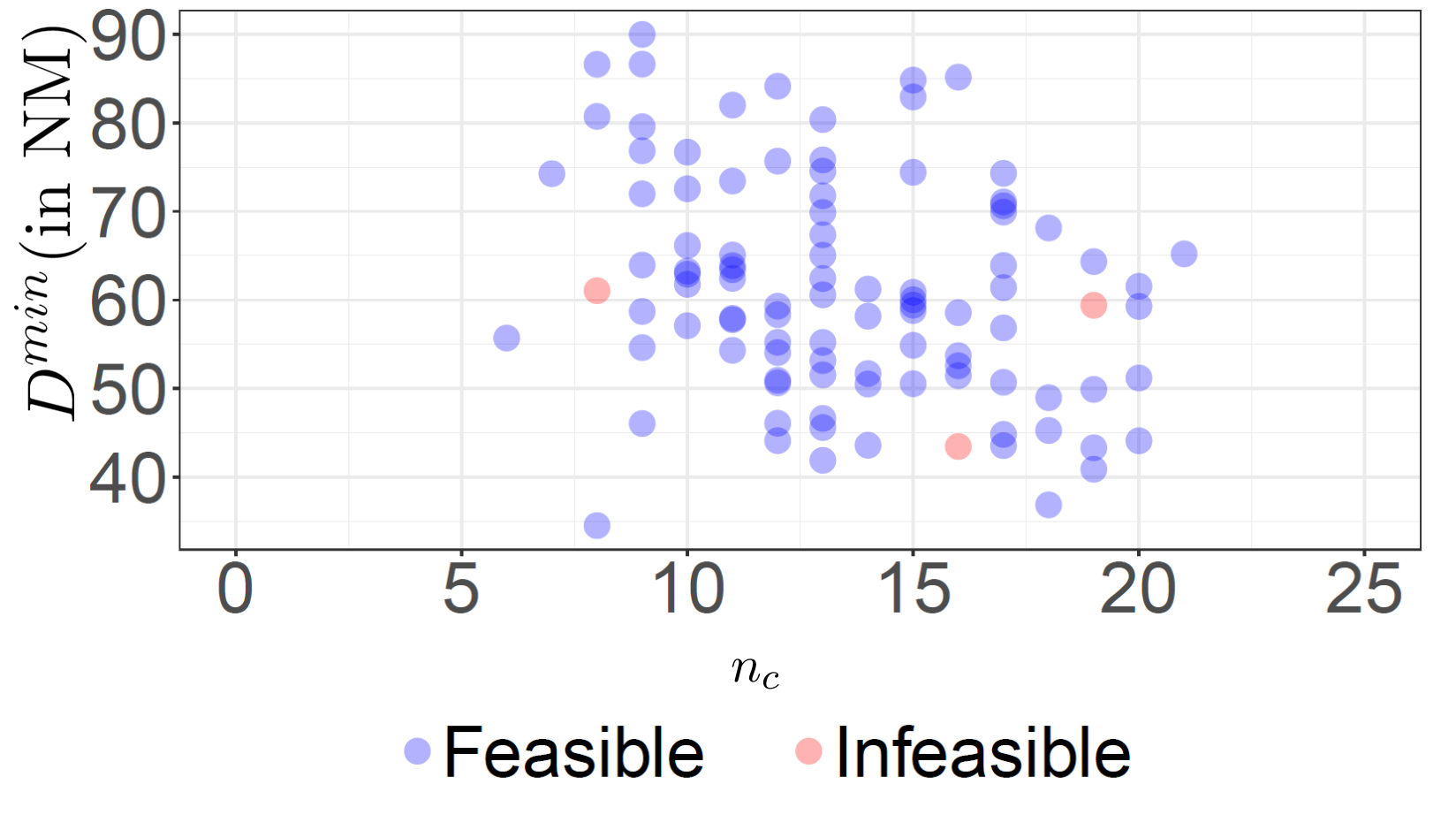}\label{20g0}} \hfill
	\subfloat[$\Gamma = 2$ ]{\includegraphics[width=0.45\textwidth]{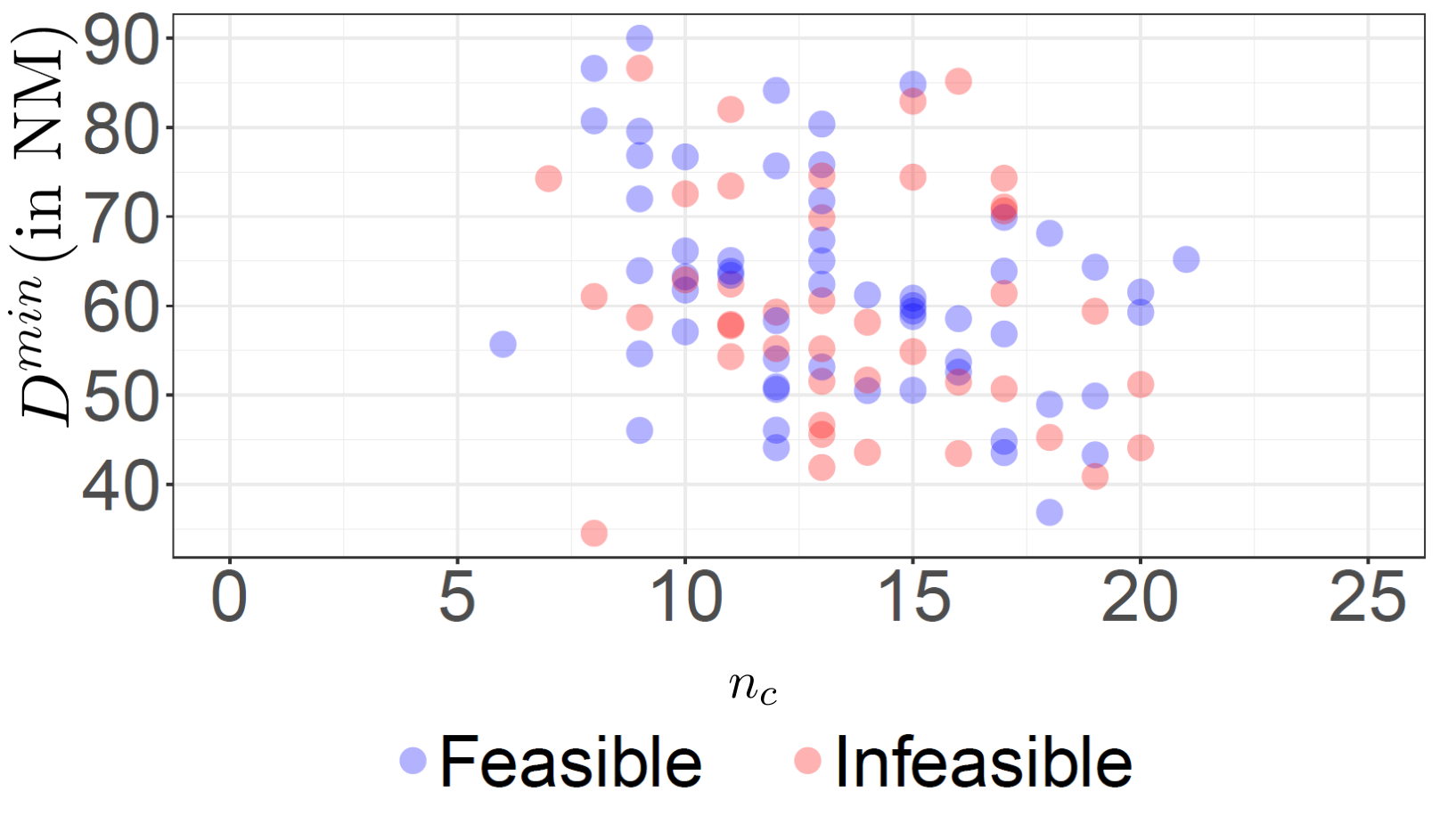}\label{20g2}}\\
	\subfloat[$\Gamma = 3$ ]{\includegraphics[width=0.45\textwidth]{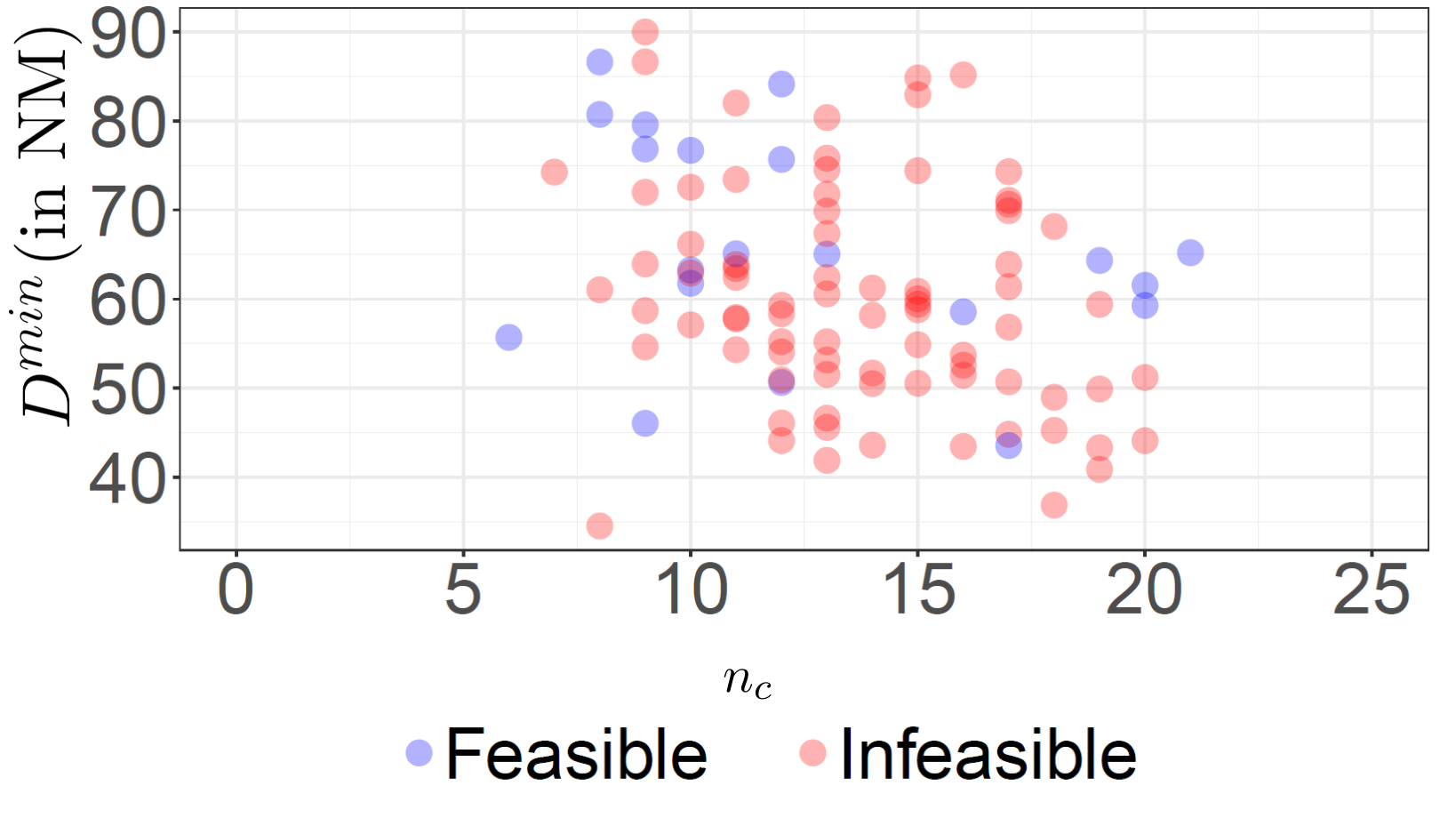}\label{20g3}} \hfill
	\subfloat[$\Gamma = 4$ ]{\includegraphics[width=0.45\textwidth]{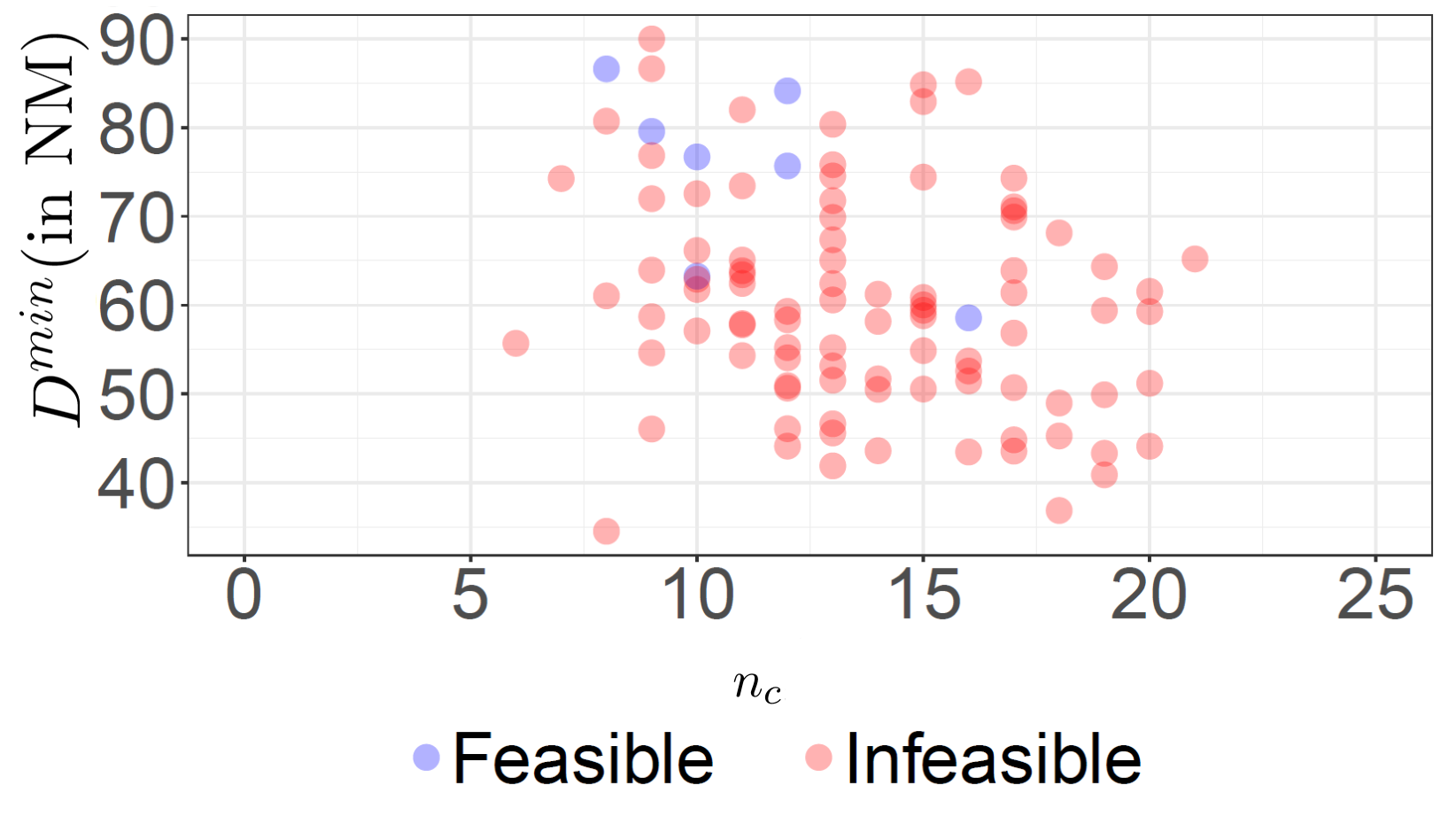}\label{20g4}}\\
	\caption{Feasibility of 100 RCP-20 instances based on the number of conflict $n_c$ and the total pairwise minimal distance between all pairs of aircraft $D^{\text{min}}$ using $\bar{\epsilon} = 5\%$.}
	\label{ilus4} 
\end{figure}

\begin{figure}[ht]
	\centering\vspace{-2pt}
	\subfloat[$\bar{\epsilon} = 2.5\%$]{\includegraphics[width=0.45\textwidth]{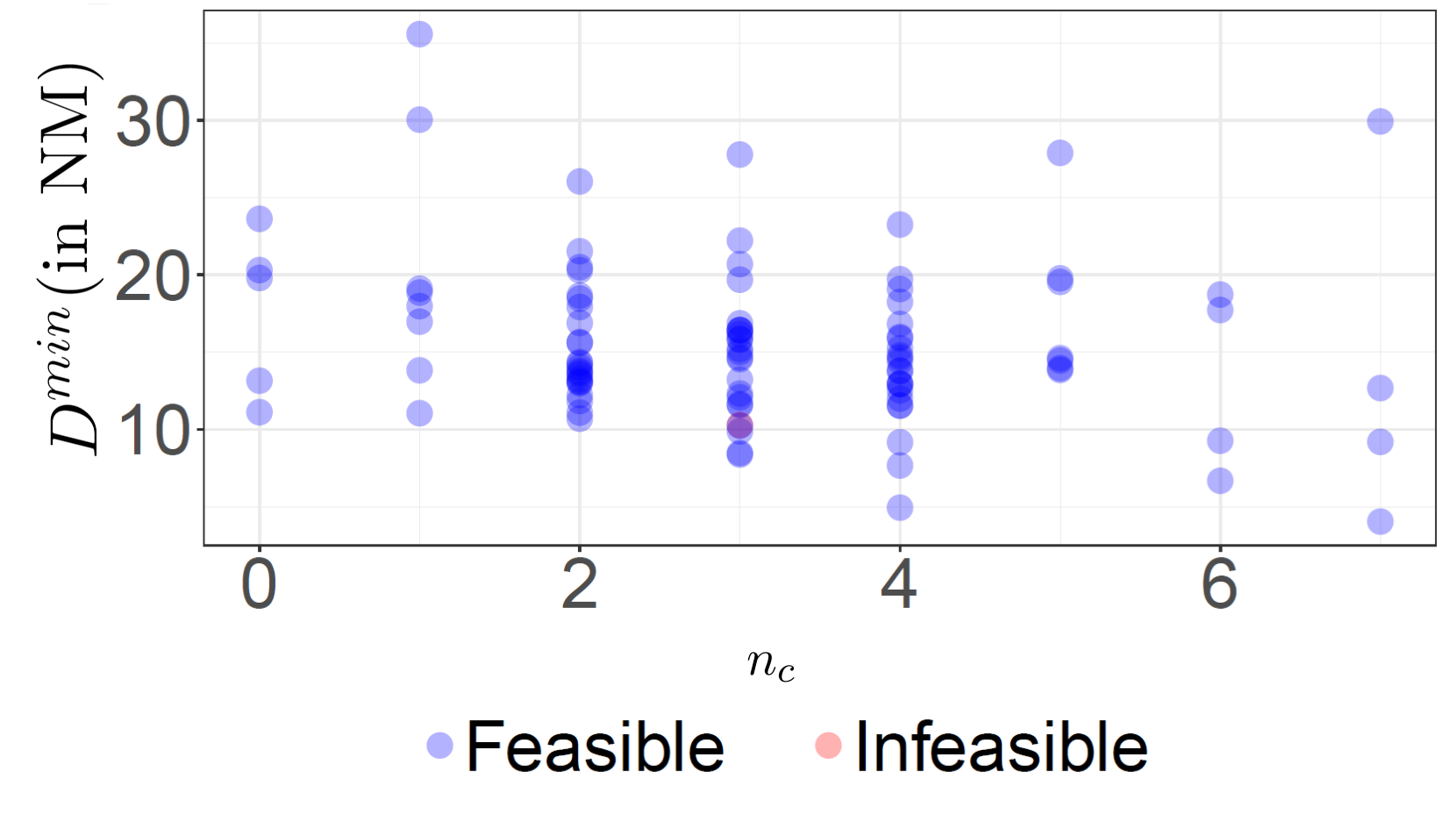}\label{10e1}} \hfill
	\subfloat[$\bar{\epsilon} = 5.0\%$]{\includegraphics[width=0.45\textwidth]{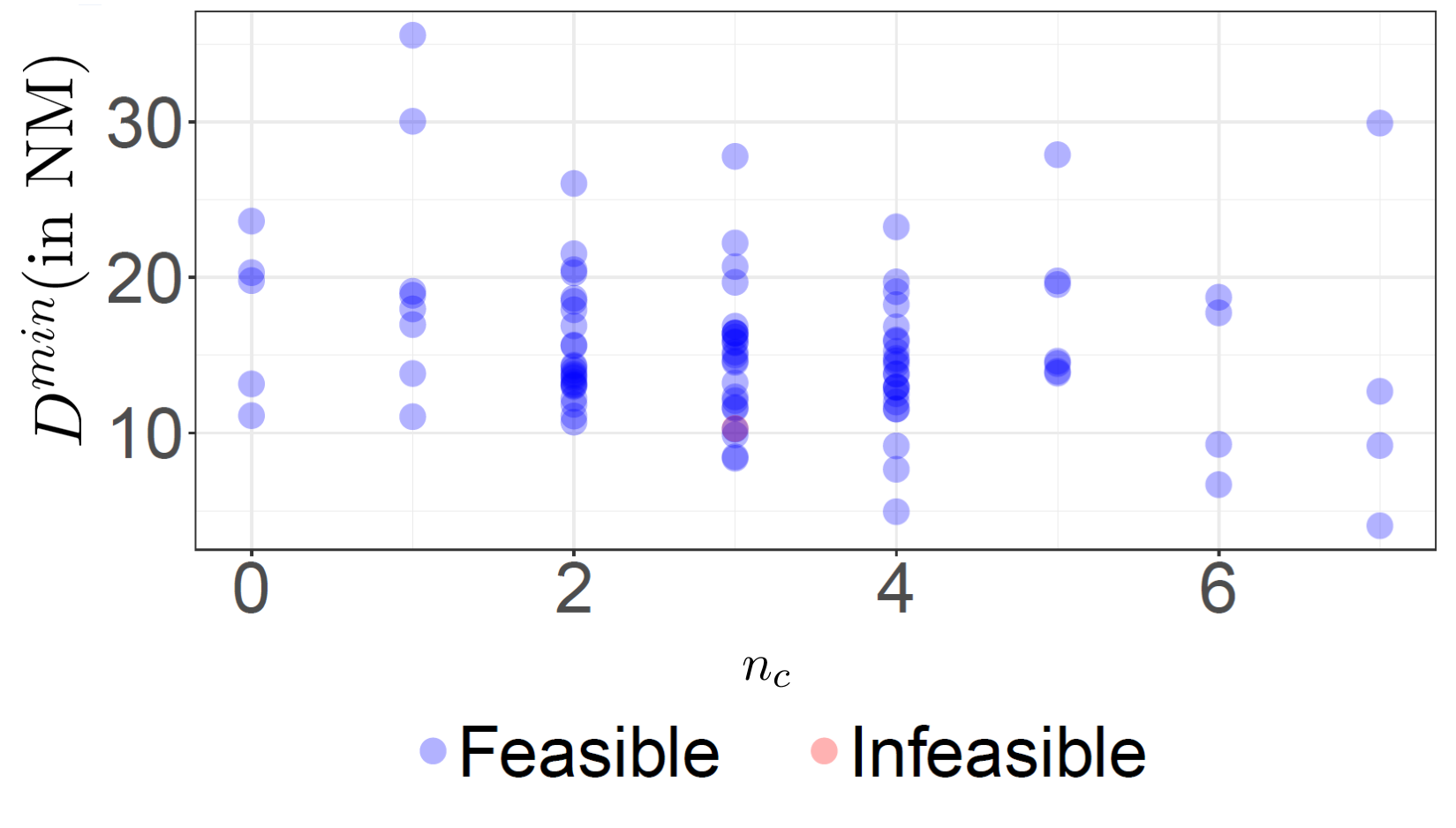}\label{10e2}}\\
	\subfloat[$\bar{\epsilon} = 7.5\%$]{\includegraphics[width=0.45\textwidth]{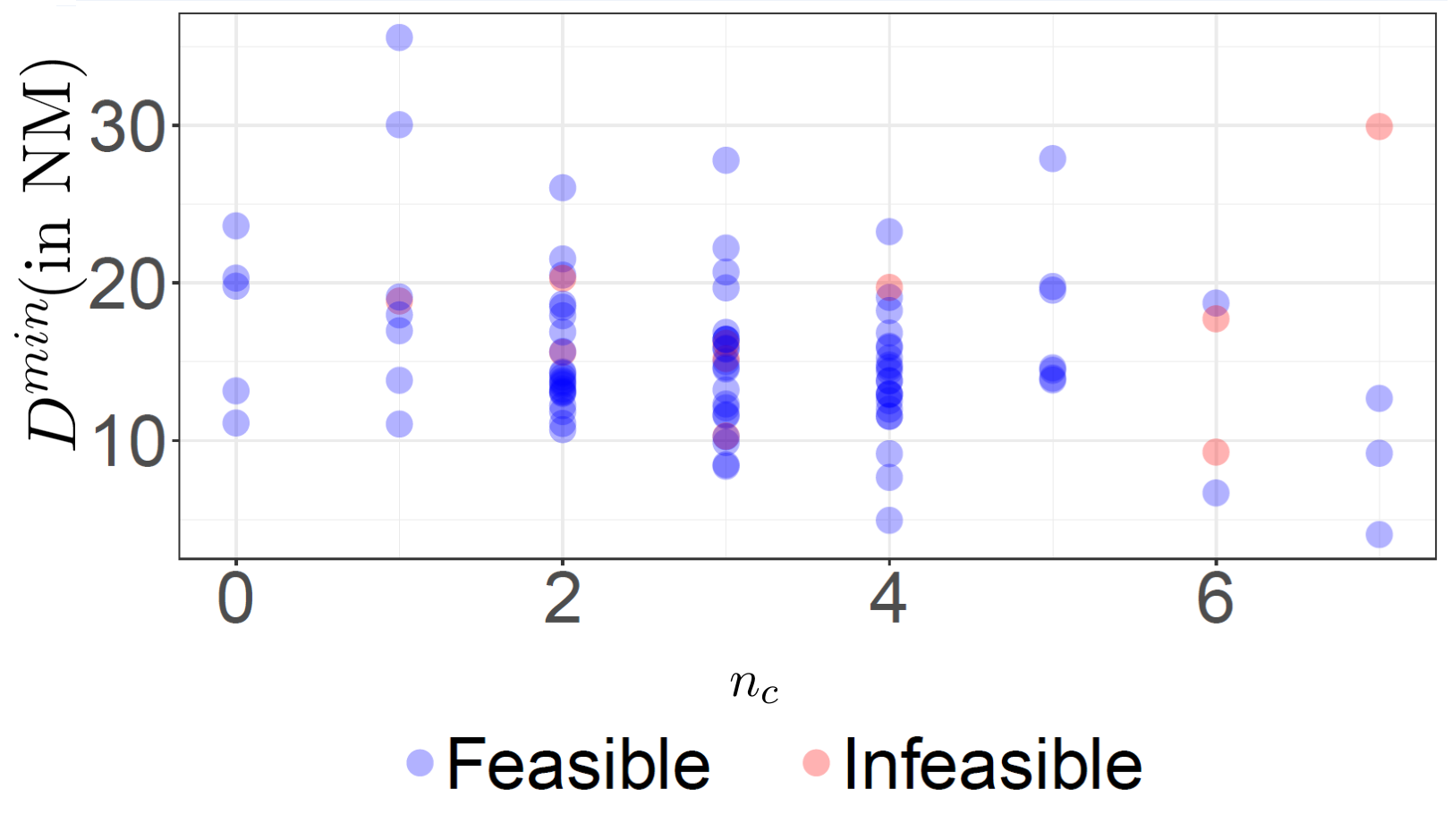}\label{10e3}} \hfill
	\subfloat[$\bar{\epsilon} = 10\%$]{\includegraphics[width=0.45\textwidth]{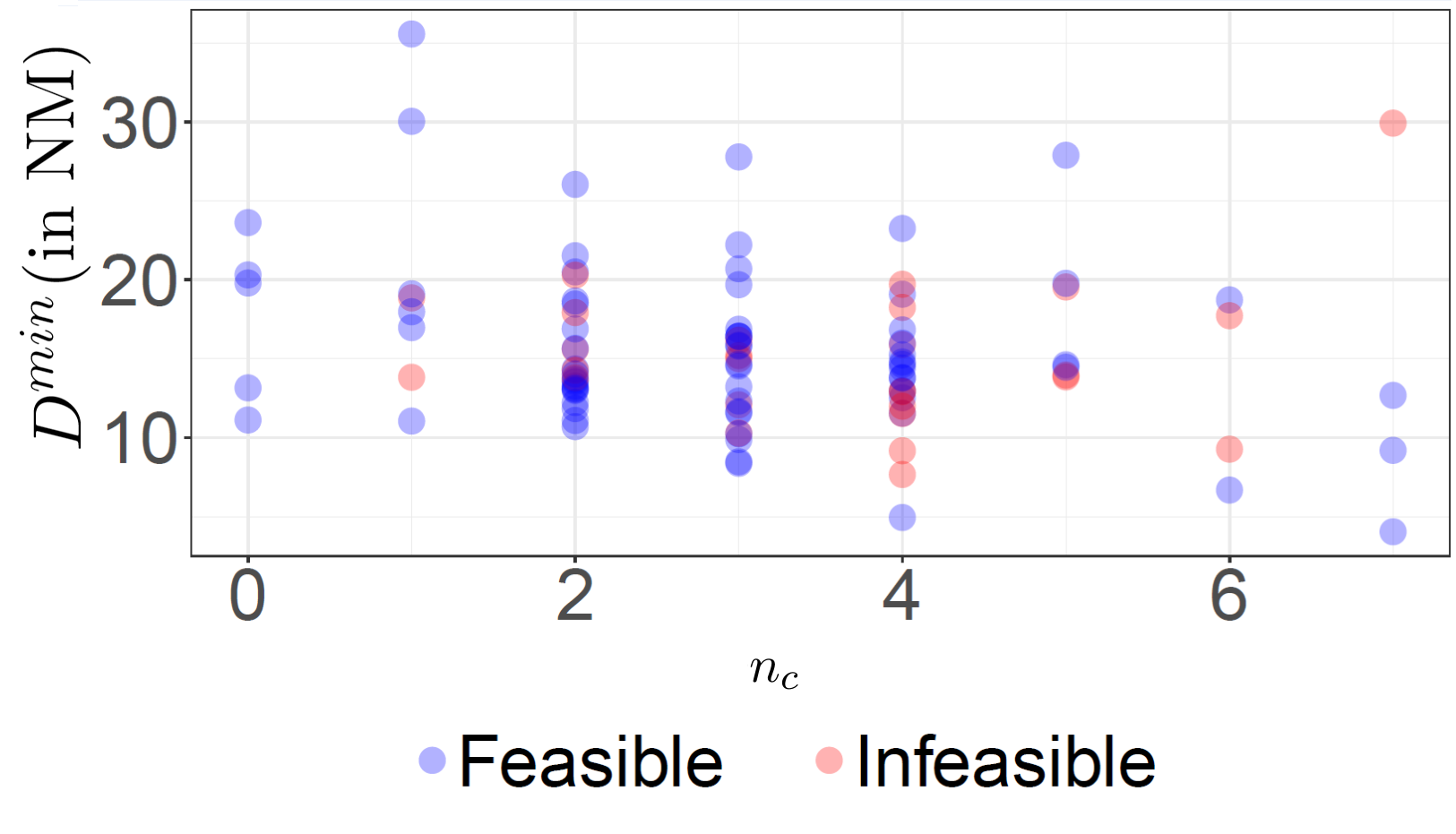}\label{10e4}}\\
	\caption{Feasibility of 100 RCP-10 instances based on the number of conflict $n_c$ and the total pairwise minimal distance between all pairs of aircraft $D^{\text{min}}$ using $\Gamma = 4$.}
	\label{ilus5} 
\end{figure}

\begin{figure}[ht]
	\centering\vspace{-2pt}
	\subfloat[$\bar{\epsilon} = 2.5\%$]{\includegraphics[width=0.45\textwidth]{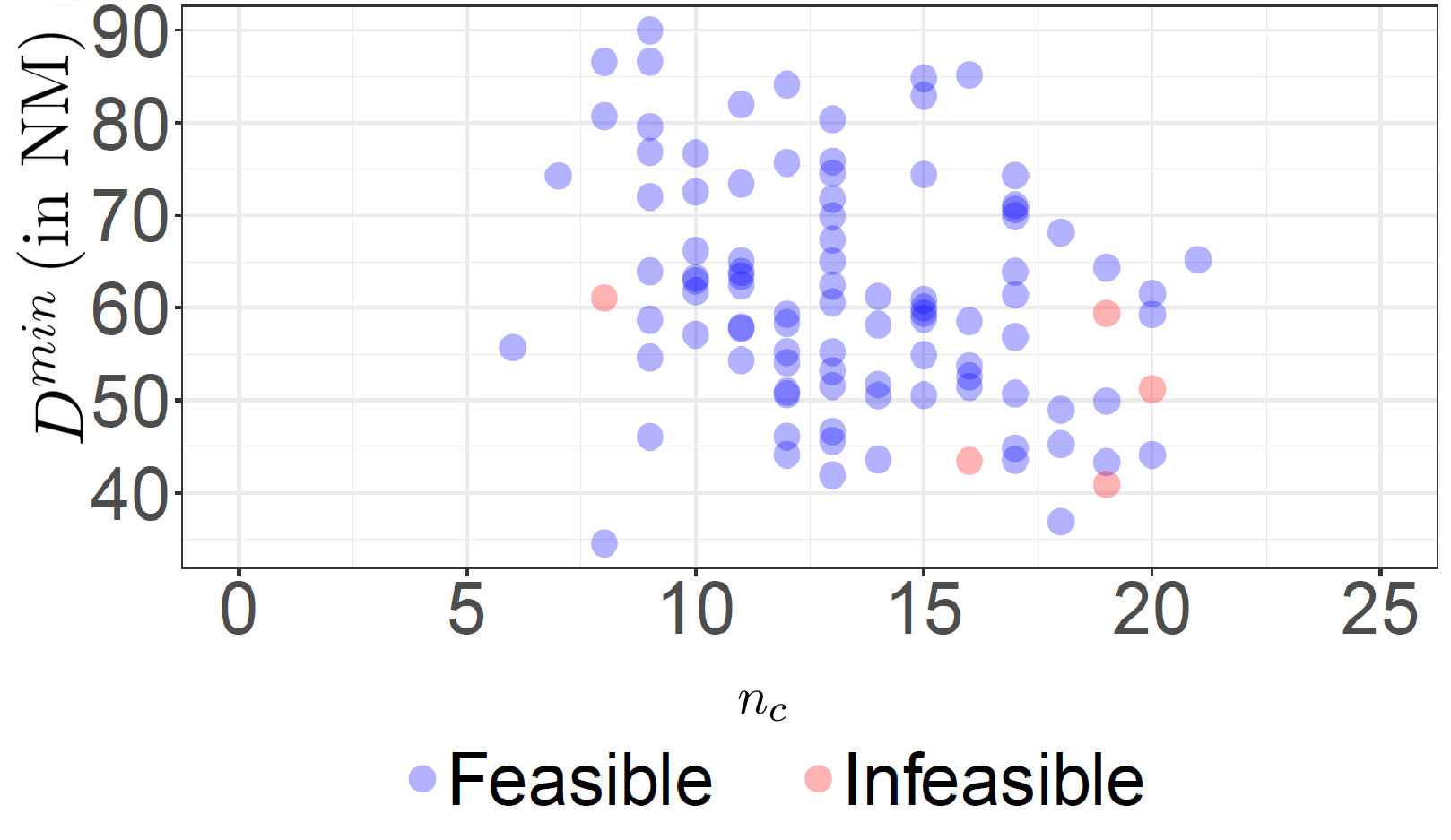}\label{20e1}} \hfill
	\subfloat[$\bar{\epsilon} = 5.0\%$]{\includegraphics[width=0.45\textwidth]{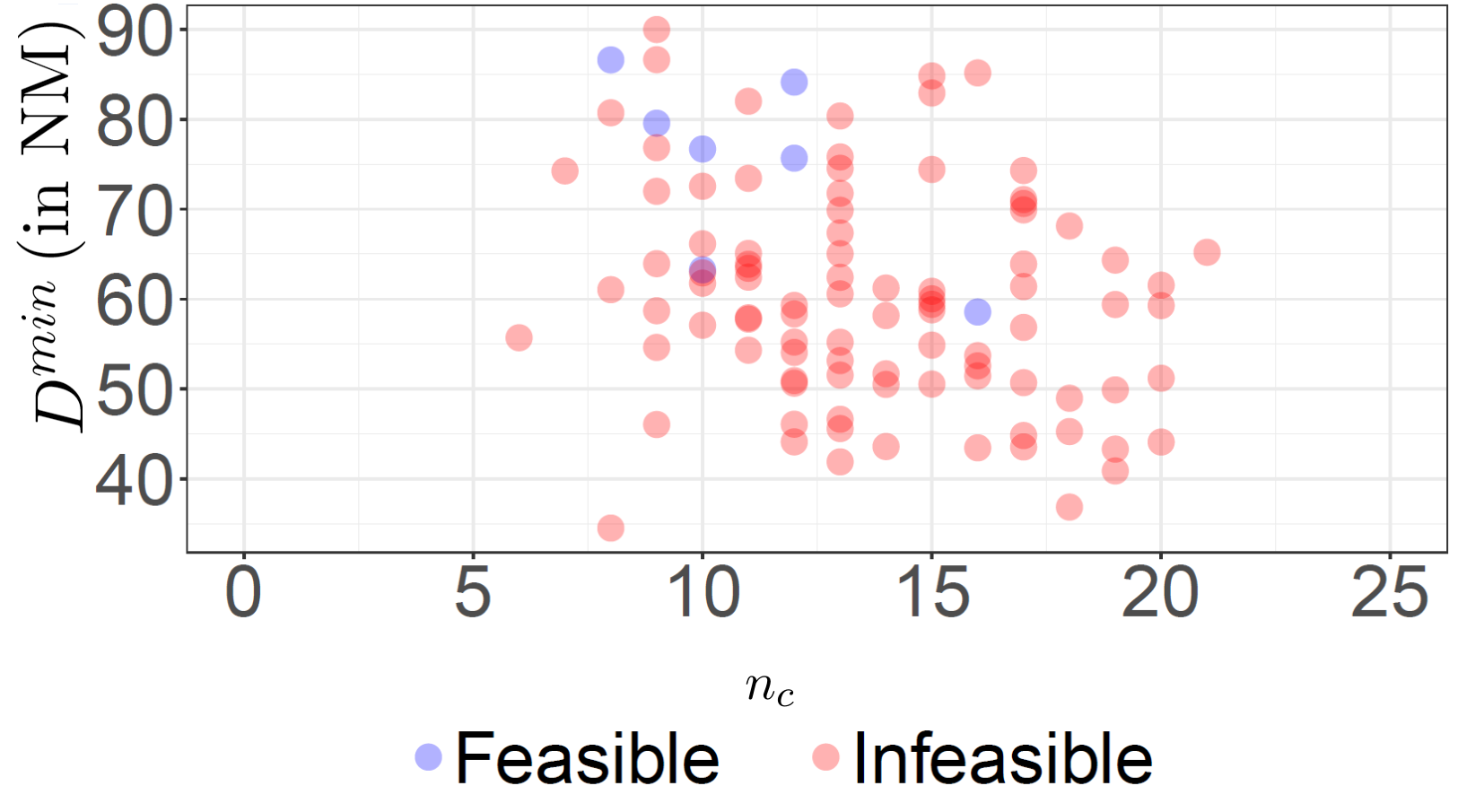}\label{20e2}}\\
	\subfloat[$\bar{\epsilon} = 7.5\%$]{\includegraphics[width=0.45\textwidth]{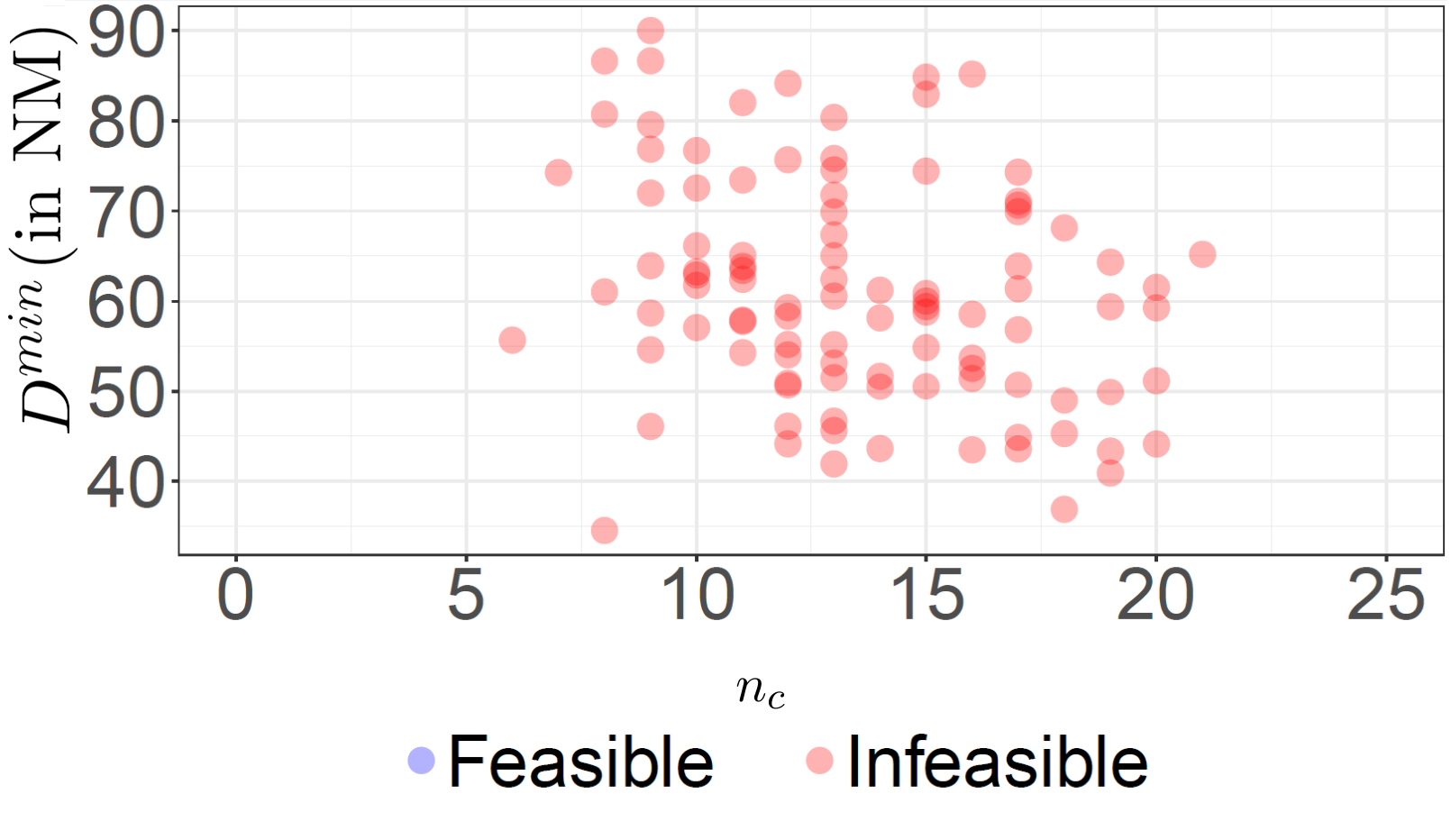}\label{20e3}} \hfill
	\subfloat[$\bar{\epsilon} = 10\%$]{\includegraphics[width=0.45\textwidth]{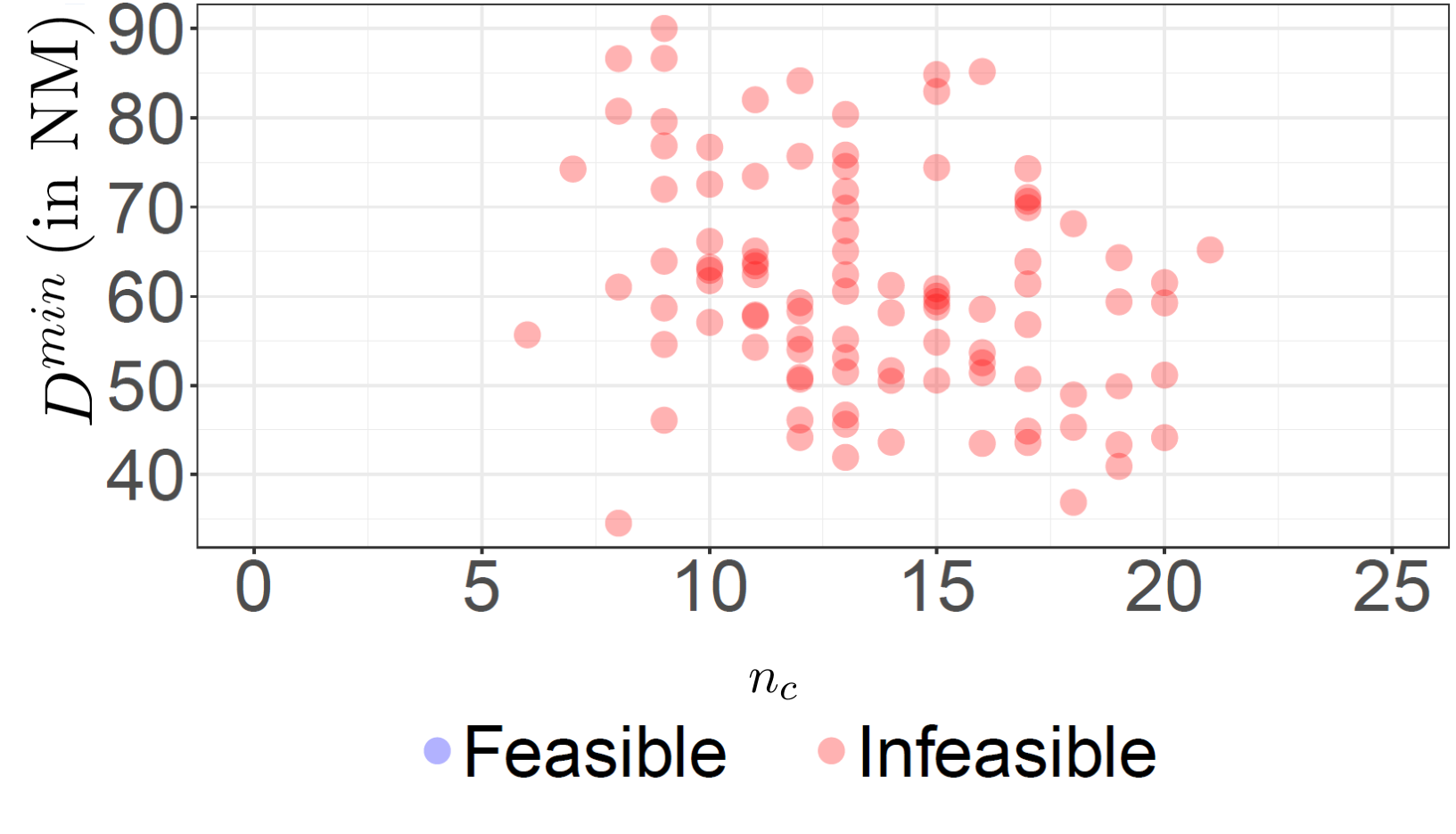}\label{20e4}}\\
	\caption{Feasibility of 100 RCP-20 instances based on the number of conflict $n_c$ and the total pairwise minimal distance between all pairs of aircraft $D^{\text{min}}$ using $\Gamma = 4$.}
	\label{ilus6} 
\end{figure}

In this section, we explore the correlation between instance feasibility and infeasibility and the characteristics of the instance. Our goal is to identify which features of conflict resolution instances can explain the existence of feasible solutions to the robust ACRP. Since CP instances are fully symmetric  we do not examine these instances. Further, since most RCP-30 instances are infeasible, we focus the analysis to RCP-10 and RCP-20 instances. For each group of 100 instances, we generate two-dimensional scatter plots where the x-axis dimension represents the number of conflicts ($n_c$) of this instance and the y-axis represents the total minimal pairwise distance $D^{\text{min}} = \sum_{(i,j) \in \P} d_{ij}(\tm)$ (in NM) of this instance. Each dot in the plots represents in of the 100 instances of the corresponding group (RCP-10 or RCP-20) and the color of the dot indicates whether this instance is feasible (blue) or infeasible (red). We seek to verify that instances with a low total minimal pairwise distance and a high number of conflicts are more likely to be infeasible when increasing the level of robustness or the maximum uncertainty in the system. The results of this experiment are reported in Figures \ref{ilus3}, \ref{ilus4}, \ref{ilus5} and \ref{ilus6}: Figures \ref{ilus3} and \ref{ilus4} illustrate the outcome for varying level of robustness for RCP-10 and RCP-20 respectively; while Figures \ref{ilus5} and \ref{ilus6} illustrate the outcome for varying maximum uncertainty for RCP-10 and RCP-20, respectively.

%As expected, for instances where the initial minimal separation per pair of aircraft is more than 5 NM, there is no expected initial conflict. On the other hand, if the pairwise distance is small or less than 5 NM as the minimal criteria, conflicts are expected. 

%Based on the configuration of the instances used for testing, when the total minimal distance is small, meaning that all aircraft are close to each other and the density of conflict is higher, it is more likely that such instances cannot be solved for all level of robustness. In Figure \ref{ilus3} and \ref{ilus4}, we observe such behaviour. 

In Figure \ref{ilus3}, we observe that for all values of $\Gamma$, all RCP-10 instances but one (out of 100) with $n_c = 3$ conflicts and $D^{\text{min}} \approx$ 10 NM can be solved. Figure \ref{ilus4} reveals that for RCP-20 instances and $\Gamma=1$ only 3 infeasible instances. Increasing $\Gamma$, we observe that the proportion of infeasible instances increase more rapidly the bottom-right quadrant of the plot, suggesting that on average instances with a high number of conflict and low total minimal pairwise distance are more likely to be infeasible. Examining RCP-10 instances, we observe in Figure \ref{ilus5} that only a single instance is infeasible for $\bar{\epsilon}=2.5\%$ and $\bar{\epsilon}=5.0\%$. Increasing the maximum uncertainty to $\bar{\epsilon}=7.5\%$ and $10\%$, reveal that more instances are infeasible and those instances exhibit a smaller total pairwise minimal distance and a higher number of conflicts. As shown in Figure \ref{ilus6}, the proportion of infeasible RCP-20 instances increases rapidly with the maximum uncertainty. Using $\bar{\epsilon}=2.5\%$, only 5 instances (out of 100) are infeasible and they mostly represent cases where the total minimal distance is small relative other instances, and the number of conflicts is above 15 which is relatively high. Using $\bar{\epsilon}=5.0\%$ value, 94\% of instances are infeasible and the feasible instances tend to have a large total minimal pairwise distance (above 70 NM) which suggests that aircraft have enough ``room'' to accommodate robust conflict-free trajectories. We find that for $\bar{\epsilon}\geq 7.5\%$ all RCP-20 instances are infeasible.

%The same analysis could be reflected when comparing the minimum separation distance. Overall, it is clear that more conflict requires, the total separation between tend to decrease. This is due to more aircraft and more conflict in a limited space causes the model to obtain solutions that only satisfy the bare minimum of separation while the lower density instances, there is more degree of freedom and therefore higher separation distance between pairs. By observing the results for RCP-10 with a fixed uncertainty set (Figure \ref{ilus8}) shows that as $\Gamma$ increases, the total distance between pairs decreases and it is even clearer in Figure \ref{ilus7} where the concentration of instances is the lower region of the graph. For different uncertainty sets as in Figures \ref{ilus9} and \ref{ilus10}, fewer instances are feasible, which explains fewer data points as $epsilon$ grows especially for instances with more conflicts. 

\section{Conclusion}
\label{con}

\subsection{Summary of Findings} 

In this study, we proposed a new formulation for the robust Aircraft Conflict Resolution Problem (ACRP) under trajectory prediction uncertainty. We introduced an aircraft trajectory prediction uncertainty model based on aircraft velocity components and showed that this approach can be incorporated into a robust optimization formulation for the ACRP. We considered the ACRP with continuous speed and heading control maneuvers and showed that robust separation constraints can be reformulated as tractable integer-linear constraints using state-of-the-art approaches in robust optimization. We adapted the complex number formulation for the ACRP and an exact algorithm to solve the resulting robust complex number formulation model to optimality. We conducted a series of numerical experiments on benchmarking instances of the literature to explore the behavior of the robust ACRP and test the computational performance of the proposed approach. A total of $1535$ benchmarking instances were used. These instances include two types of ACRPs with up to 30 aircraft per instance for different levels of randomness reflected by the level of robustness and the size of the uncertainty sets. The performance of the proposed solution method highlights the scalability of the approach compared to existing deterministic methods in the literature. We find that increasing the level of robustness and/or the maximum uncertainty in the model may lead to infeasible problems, notably when the number of aircraft and conflicts is significant. Further, pre- and post-optimization analyses reveal that the number of conflicts and the total pairwise minimal distance between aircraft can explain the behaviour of the proposed robust ACRP in terms of instance feasibility. To the best of our knowledge, this is the first exact robust optimization formulation for the ACRP.

%In our numerical experiments, we compared with the deterministic case ($\Gamma=0$). Our results reveal that increasing $\Gamma$ and the uncertainty set, the resulting model is still very competitive and we find that the runtime of the robust formulation does not increase significantly compared to the deterministic counterpart. Even though the gamma formulation requires more variables and constraint, it can solve most of the cases for lower values of $\Gamma$ and uncertainty sets. However, higher values for those parameters, it becomes more challenging to get feasible solutions within 10 minutes as the time limit and changes in $\Gamma$ are more likely to be solved instead of changes in $\bar{\epsilon}$ as shown in previous section, which means that controlling the level of robustness is easier than altering the dimension of such robustness. 

\subsection{Future Research and Perspectives}

The proposed approach for the robust ACRP relies on several assumptions which may be limiting and could be further relaxed. First, the proposed formulation assumes that uniform motion laws apply, which translate into infinite acceleration and deceleration rates. Although this is a common assumption in the literature, dedicated efforts are required to generate more accurate models. Second, aircraft trajectory recovery is not taken into account. Changing aircraft heading angles may require subsequent maneuvers to return aircraft to their target trajectory, which may translate into significant penalties for airlines. Since the determination of aircraft recovery trajectories is function of aircraft predicted trajectories, accounting for uncertainty is critical to ensure safe recovery trajectories. 

In terms of computational experiments, our analysis has focused on uniform uncertainty model across aircraft. More realistic air traffic scenarios, possibly generated from weather data, should be explored so as to gain more practical insights into the impact of robustness in aircraft conflict resolution. Unfortunately, such data is not publicly available and future research efforts may be needed to augment existing data repositories for conflict resolution problems. The coordination of aircraft conflict resolution maneuvers also presents considerable operational challenges. As in the vast majority of exact methods for the ACRP, the proposed approach assumes that all maneuvers start at the same time, which may not always be feasible in practice. Hence, there is a need to develop exact methods that are able to account for asynchronous conflict resolution maneuvers. Exploring the role of robust and stochastic optimization approaches to develop such rich formulations is a promising avenue of research in aircraft conflict resolution.

\bibliographystyle{elsarticle-harv}
\bibliography{biblio}

\end{document}